\documentclass[11pt,a4paper]{article}
\usepackage[margin=26mm]{geometry}
\usepackage{fontspec}
\usepackage{amsmath,amssymb,amsthm,mathtools}
\numberwithin{equation}{section}
\usepackage{enumitem,booktabs,array,longtable,microtype}
\usepackage[colorlinks,citecolor=blue]{hyperref}
\usepackage[numbers,sort&compress]{natbib}
\usepackage{indentfirst}
\usepackage{fancyhdr}
\usepackage{mathrsfs}
\newtheorem{theorem}{Theorem}[section]
\newtheorem{lemma}{Lemma}[section]

\newtheorem{corollary}{Corollary}[section]
\theoremstyle{definition}
\newtheorem{definition}{Definition}[section]
\newtheorem{remark}{Remark}[section]

\newcounter{namedassumptionctr}
\newcommand{\namedassumption}[2]{%
	\par\medskip
	\begingroup
	\renewcommand{\thenamedassumptionctr}{#1}%
	\refstepcounter{namedassumptionctr}%
	\label{#2}%
	\endgroup
	\noindent\textbf{(#1)}\quad\ignorespaces
}
\allowdisplaybreaks
\date{}
\begin{document}
	\date{}

	\title{Random pullback attractors for evolution equations with nonlocal perturbations and fractional Brownian motion}

	\author{Shiqi Li$^1$, Jiawei He$^{1,*}$, Jianhua Huang$^2$\\[3pt]
		\small $^1$School of Mathematics, Guangxi University, Nanning 530004, China\\
		\small $^2$College of Science, National University of Defense Technology, Changsha 410073, China\\
		\small $^*$Corresponding author: \texttt{jwhe@gxu.edu.cn}}

	\date{}
	\maketitle
	\begin{abstract}
		We study semilinear evolution equations with time-nonlocal perturbations of  \(\mu\ge0\) and additive fractional Brownian noise with Hurst parameter \(H\in(1/2,1)\). Under the  regularity assumptions and a globally Lipschitz nonlinearity, we establish global existence and uniqueness in the concept of pathwise mild solutions. We overcome the difficult  in which the solution family lacks the cocycle property, by shifting the corresponding stochastic Volterra equation on a weighted space, and then the continuous random dynamical system follows. Based on a dissipativity condition, we establish the pullback random asymptotic compactness and the compact random  pullback attractor \(\mathcal A_\mu\). 
		Finally, we prove that the family of random pullback attractors
		\(\mathcal A_\mu\) is upper semicontinuous at \(\mu=0\) almost surely.

		\medskip
		
	\end{abstract}


	\section{Introduction}
	\label{sec:introduction}
	It is known that the Fick law $\vec F=-\lambda\nabla u$ and the law of  energy balance $\partial_t \epsilon =-\operatorname{div}\vec F$ implies  the heat equation in homogeneous medium of 3D regular enough bounded domain for  the body temperature $u$ with zero initial temperature satisfying $\epsilon =\sigma u$. This model is also adapted to widely used to model the heat conduction in materials and diffusion phenomena in natural sciences and so on. However, in nonhomogeneous media, the
	constitutive equation $\epsilon =\sigma u$ needs to modify by
	$$
	\epsilon =\sigma u +\gamma\, k*u,  
	$$
	that describes a case of diffusion process  that represents the transition from sub-diffusion (as $t \to 0$) to normal diffusion (as $t \to \infty $), see e.g. \cite{Chen17,he_huang,Metzler2000}, $k$ is a suitable relaxation kernel function. The energy balance for the system is deduced by the form
	$$
	-\operatorname{div}\vec F +g(u,t)=\sigma \partial_t u +\gamma\, \partial_t (k*u) ,
	$$
	where $g(u,t)$ is the energy supply. Consequently, we arrive at the following nonlocal equation 
	$$
	\sigma \partial_t u +\gamma\, \partial_t (k*u) =\lambda \Delta u+g(u,t).
	$$
	In many practical settings, the effects of external noise and external forcing cannot be ignored. This observation naturally motivates the consideration of an extended formulation  
	$$
	g(u,t) =f(u)+G\,B^H(t),
	$$
	here \(B^H\) is a two-sided \(L_2\)-valued fractional
	Brownian motion with Hurst parameter \(H\in(1/2,1)\), $G$ is the noise coefficient. In this way,  we are interested in the following SPDE for \(\mu\geq0\), 
	\begin{equation}\label{e1.1}
		\begin{cases}
			d u
			=-\mu\,d(k*u)+(A u+f(u))dt+G\,d\dot{B}^H(t),
			& t>0,\\ 
			u(0)=x\in X, 
		\end{cases}
	\end{equation}
	where $-A$ is a strictly positive and symmetric operator with a compact inverse, and it is the generator of an analytic semigroup on $X$.
	
	As shown by Caraballo et al.~\cite{crc}, the constitutive relation for the heat conduction model can also be replaced by
	$$
	\epsilon =-c_0\nabla u+\int_{-\infty}^t \gamma(t-s)\nabla u(s,\cdot)ds,
	$$
	in a 3D regular enough bounded domain, here $\gamma$ is the heat flux relation function and $c_0>0$ is the thermal conduction. And then they established the pullback attractors for a stochastic reaction-diffusion equation with white noise. Following this work, Liu et al.~\cite{liu_caraballo_kloeden} extended this model including two memory terms and subcritical and critical nonlinearities.
	They also studied the random pullback attractors and their unpper semi-continuity. 
	
	It is noteworthy that the current model \eqref{e1.1} cannot be directly applied to the semigroup framework and the Lax-Milgram developed in~\cite{liu_caraballo_kloeden}, since the nonlocality embedded in \eqref{e1.1} is characterized by derivative forms, instead of the integral-type memory effect.
	In recent years, 
	a prominent starting point for investigating the long-term dynamics of non-local systems is to construct the theoretical framework on weighted spaces with a compact open topology, which was first pioneered by Miller and Sell \cite{miller_sell} in their research on Volterra equations. Doan and Kloeden
	\cite{doan_kloeden} followed this idea for studying the semi-dynamic of autonomous Caputo equations on finite dimension, while  Cui and Kloeden \cite{cui_kloeden} constructed a skew-product semiflow and skew-product attractors 
	for nonautonomous Caputo equations. Moreover, Cui et al.~\cite{cui_kloeden_xin} also 
	obtained a compact absorbing set and the H\"older regularity of the attractors. 
	Zhang et al.  \cite{zhw} constructed a mean-square
	skew-product semiflow and obtained a forward attracting set in the
	weak mean-square topology under the framework of the tempered fractional Brownian motion.
	On infinite dimension, Peng, Deng and He \cite{deng} studied the uniform attractor for a generalized Rayleigh--Stokes equation on
	\(L^p(\mathbb R^d)\). He and Huang \cite{he_huang} established the skew-product attractors for the time nonlocal evolution equations for a more general Hilbert space.
	
	We remark that fractional Brownian motion (fBm) is not a semimartingale for
	\(H\in(1/2,1)\), this means that the classical It\^o integral may not valid to handle problem \eqref{e1.1}. Nevertheless,
	when \(\mu=0\) in \eqref{e1.1}, Chen et al.~\cite{chen} used 
	Z\"ahle's pathwise integral to construct a random dynamical system, 
	and Gao et al. \cite{gao} proved the existence
	of random attractors under a smallness condition on the nonlinear multiplicative noise by virtue of the Young integral, see e.g. \cite{young}. Motivated by
	the aforementioned results, we use the Young integral to prove the Cocycle property.
	Moreover,
	our main results are listed as follows points. 
	since the semigroup theory cannot be direct applied to the current system driven by nonlocal perturbation, the first difficulty is to obtain the appropriate pathwise  mild solution by virtue of subordination principle.
	The second difficulty is to establish the H\"{o}lder regularity motivated by the property of Young integral, indeed we require the regularity of stochastic convolution is better than that of the solution. We then shift the stochastic convolution to complete the existence of solution and the cocycle properties in a weighted Banach space. 
	Moreover, by a  random dynamical system, we need to verify the continuity and measurability of the cocycle. The next one is that, we overcome the difficulty that the
	pullback asymptotic compactness is not direct obtained by compact embedding and Arzel\`a-Ascoli theorem in the weighted Banach space, and then
	we establish the existence of random
	pullback attractor \(\mathcal A_\mu\),  which has a tempered random bound. 
	Finally, we study the upper semi-continuity of attractors on the limits behavior as \(\mu\to 0\). As our knowledge, 
	previous research on the dynamics of stochastic infinite-dimensional complex systems has mainly
	focused on noise types such as additive Wiener noise, linear multiplicative Wiener noise, or additive fractional
	noise. This paper seeks to establish a research framework for addressing the dynamics of infinite-dimensional
	complex systems driven by linear multiplicative fractional noise.
	
	The paper is organized as follows.
	Section~\ref{sec:preliminaries} gives estimates for the kernels and
	the stochastic convolution. Section~\ref{sec:pathwise-existence-uniqueness}
	proves existence and uniqueness of global pathwise mild solutions.
	Section~\ref{sec:lifted-cocycle} constructs a continuous random dynamical
	system on \(\mathcal Y\). Section~\ref{sec:generalized-pullback-attractor}
	establishes the existence of a minimal random pullback attractor.
	Finally, Section~\ref{subsec:vanishing-memory} proves the
	convergence of random attractors as \(\mu\to0\).

	\section{Preliminaries}
	\label{sec:preliminaries}
	Let \(\Omega\subset\mathbb R^d\) be a bounded domain with smooth boundary, and set
	\(
	X:=L^2(\Omega)
	\)
	with norm \(\|\cdot\|\) and inner product
	\((\cdot,\cdot)\).
	Let \(A:D(A)\subset X\to X\) be a densely defined linear operator
	such that \(-A\) is strictly positive and self-adjoint with compact
	inverse.
	Then \(-A\) generates an analytic semigroup 	\(\{e^{-tA}\}_{t\geq0}\) on \(X\).
	There exists an orthonormal basis \(\{\varphi_n\}_{n\geq1}\) of \(X\)
	consisting of eigenvectors of \(-A\), with corresponding eigenvalues
	\[
	-A\varphi_n=\lambda_n\varphi_n,\qquad
	0<\lambda_1\leq\lambda_2\leq\cdots,\qquad
	\lambda_n\to\infty.
	\]
	For \(\delta\geq0\), let \(
	X_\delta:=D(-A^\delta)=\left\{
	v=\sum_{n=1}^{\infty}v_n\varphi_n:
	\sum_{n=1}^{\infty}
	\lambda_n^{2\delta}|v_n|^2<\infty
	\right\}\) with
	the norm \(\|v\|_\delta=\|(-A)^\delta v\|_X\).
	Since \(-A\geq\lambda_1 I\) with \(\lambda_1>0\), this norm is
	equivalent to the graph norm of \((-A)^\delta\).
	For every \(\delta>0\), the embedding
	\(
	X_\delta\hookrightarrow X
	\)
	is compact.
	Throughout this paper, \(\mathbb N=\{1,2,\ldots\}\), \(\mathbb R_+:=[0,\infty)\), $\widehat{\cdot}$ denotes the Laplace transform, and $*$ denotes convolution on \(\mathbb R_+\).  We recall the following definitions; see \cite{sc,pru}.
	A function \(g\in C^\infty(0,\infty)\) is called
	\emph{completely monotone} if
	\[
	(-1)^n g^{(n)}(t)\geq0,
	\qquad t>0,\quad n=0,1,\ldots.
	\]
	A function \(g:(0,\infty)\to[0,\infty)\) is called
	a \emph{Stieltjes function} if it admits the representation
	\[
	g(t)=\frac at+b+
	\int_{0}^{\infty}\frac1{t+s}\,\rho(ds),
	\qquad t>0,
	\]
	where \(a,b\geq0\) and \(\rho\) is a positive Borel measure
	satisfying
	\(
	\int_{0}^{\infty}\frac1{1+s}\,\rho(ds)<\infty.
	\)
	Finally, a function \(g:(0,\infty)\to[0,\infty)\)
	is a \emph{complete Bernstein function} if and only if
	\(t\mapsto g(t)/t\) is a Stieltjes function.
	
	\subsection{Operator kernel estimation}
	\label{subsec: kernel estimation}

	\namedassumption{HA}{ass:kernel}
	The kernel \(k\in L^1_{\mathrm{loc}}(\mathbb R_+)\) is nonnegative and
	completely monotone.
	
	For \(\mu\geq0\) and \(\lambda>0\), define
	\(
	\psi(\lambda)
	:=\lambda\bigl(1+\mu\widehat k(\lambda)\bigr),
	\)
	and let the kernel \(m_\mu\) through
	\begin{equation}\label{eq:mhat-def}
		\widehat m_\mu (\lambda)
		:=\frac{1}{\psi(\lambda)}
		=\frac{1}{\lambda(1+\mu\widehat k(\lambda))}.
	\end{equation}

	\begin{lemma}\label{lem:m-properties}
		Suppose that (\ref{ass:kernel}) holds, and fix \(\mu\geq0\).
		Then $\psi$ is a complete
		Bernstein function, \(\widehat m_\mu\)  is a Stieltjes function, and $m_\mu$ is
		completely monotone.  Furthermore,
		\begin{equation*}
			m_\mu+\mu k*m_\mu=1
			\qquad\text{on }\mathbb R_+,
		\end{equation*}
		and
		\begin{equation}\label{eq:m-basic-properties}
			0\leq m_\mu(t)\leq1,
			\qquad
			\lim_{t\rightarrow0}m_\mu(t)=1.
		\end{equation}
	\end{lemma}
	\begin{proof}
		Since $k$ is completely monotone, $\widehat k$ is a Stieltjes function\cite{sc}.
		Consequently, $\lambda\widehat k(\lambda)$ is a complete Bernstein
		function.  Hence $\psi(\lambda)$ is a
		complete Bernstein function. Therefore \(\widehat m_\mu\) is a
		Stieltjes function.
		Because \(k\in L^1_{\mathrm{loc}}(\mathbb R_+)\) is completely monotone,
		\(\widehat k(\lambda)\to0\) as \(\lambda\to\infty\), and
		\(\widehat m_\mu(\lambda)\leq1/\lambda\to0\). The Stieltjes representation
		therefore shows that \( m_\mu\) is  a completely
		monotone function, with no singular mass at the origin.
		In addition, by \eqref{eq:mhat-def}, we have
		\(
		\widehat m_\mu(\lambda)\bigl(1+\mu\widehat k(\lambda)\bigr)
		=\frac1\lambda.
		\)
		Using the uniqueness of the Laplace transform yields $m_\mu+\mu k*m_\mu=1$.
		Since \(m_\mu,k,\mu\ge0\),
		$m_\mu+\mu k*m_\mu=1$ also gives
		$0\leq m_\mu(t)\leq1$.  Finally,
		\begin{equation}\label{eq:m-k}
			0\leq1-m_\mu(t)=\mu(k*m_\mu)(t)
			\leq\mu\int_0^t k(s)\,ds\to0
			\qquad\text{as }t\rightarrow0,
		\end{equation}
		which proves \eqref{eq:m-basic-properties}.
	\end{proof}
	
	We extend \(m_\mu\) continuously to zero by setting \(m_\mu(0)=1\).
	Next, for \(\nu>0\), let
	\(p_{\mu,\nu}\) be the solution of
	\[
	p_{\mu,\nu}
	+\nu m_\mu*p_{\mu,\nu}
	=m_\mu.
	\]
	Equivalently,
	\begin{equation}\label{eq:pnu-laplace}
		\widehat p_{\mu,\nu}(\lambda)
		=
		\frac{\widehat m_\mu(\lambda)}
		{1+\nu\widehat m_\mu(\lambda)}
		=
		\frac{1}{\psi(\lambda)+\nu}.
	\end{equation}

	\begin{lemma}\label{lem:pnu-properties}
		For every \(\mu\geq0\) and $\nu>0$, the function $p_{\mu,\nu}$ is completely monotone and
		satisfies
		\begin{equation}\label{eq:pnu-basic-bound}
			0\leq p_{\mu,\nu}(t)\leq m_\mu(t)\leq1,
			\qquad t>0,
		\end{equation}
		and
		\(
		\lim_{t\rightarrow0}p_{\mu,\nu}(t)=1.
		\)
		Moreover, if $0<\nu_1<\nu_2$, then
		\(
		p_{\mu,\nu_1}(t)\geq p_{\mu,\nu_2}(t)
		\) for \( t>0\).
	\end{lemma}
	
	\begin{proof}
		Since $\psi$ is a complete Bernstein function,
		$1/(\psi+\nu)$ is a Stieltjes function.  Hence $p_{\mu,\nu}$ is completely
		monotone.  From 
		\(
		p_{\mu,\nu}=m_\mu-\nu m_\mu*p_{\mu,\nu},
		\)
		we have \eqref{eq:pnu-basic-bound}.  Since,
		\[
		0\leq(m_\mu*p_{\mu,\nu})(t)\leq\int_0^t m_\mu(s)\,ds\to0
		\qquad\text{as }t\rightarrow0,
		\]
		we obtain \(
		\lim_{t\rightarrow0}p_{\mu,\nu}(t)=1
		\) by \eqref{eq:m-basic-properties}.  Finally,
		from \eqref{eq:pnu-laplace}, we get
		\(
		\widehat p_{\mu,\nu_1}(\lambda)-\widehat p_{\mu,\nu_2}(\lambda)
		=
		(\nu_2-\nu_1)
		\widehat p_{\mu,\nu_1}(\lambda)
		\widehat p_{\mu,\nu_2}(\lambda)
		\) for $0<\nu_1<\nu_2$.
		Taking inverse Laplace transforms and using continuity, we obtain
		\begin{equation}\label{eq:pnu-resolvent-identity}
			p_{\mu,\nu_1}(t)-p_{\mu,\nu_2}(t)
			=(\nu_2-\nu_1)(p_{\mu,\nu_1}*p_{\mu,\nu_2})(t),
			\qquad t\geq0.
		\end{equation}
		Since both kernels are nonnegative and \(\nu_2>\nu_1\),
		it follows that \(p_{\mu,\nu_1}(t)\geq p_{\mu,\nu_2}(t)\).
	\end{proof}
	
	Consider first the deterministic linear problem
	\begin{equation*}
		\partial_tu=-\mu\partial_t(k*u)+Au+F(t),
		\qquad u(0)=x.
	\end{equation*}
	Taking the Laplace transform  gives
	\begin{equation*}
		\widehat u(\lambda)
		=
		(\psi(\lambda)I-A)^{-1}u_0
		+
		(\psi(\lambda)I-A)^{-1}\widehat F(\lambda).
	\end{equation*}
	Accordingly, we define the operator
	\begin{equation*}
		\widehat P_\mu(\lambda)
		:=
		(\psi(\lambda)I-A)^{-1},
		\qquad \lambda>0.
	\end{equation*}
	For \(v\in X\), write
	\(
	v=\sum_{n=1}^\infty v_n\varphi_n
	\) and \(
	v_n=(v,\varphi_n).
	\)
	From \(
	\widehat P_\mu(\lambda)v
	=
	\sum_{n=1}^\infty
	\frac{v_n}{\psi(\lambda)+\lambda_n}\varphi_n,
	\)
	we get
	\begin{equation}\label{eq:P-spectral-definition}
		P_\mu(t)v
		:=
		\sum_{n=1}^\infty
		p_{\mu,\lambda_n}(t)v_n\varphi_n,
		\qquad t\geq0.
	\end{equation}
	By \eqref{eq:pnu-basic-bound} and Parseval's identity, 
	the  series in  \eqref{eq:P-spectral-definition} converges in $X$.
	
	\begin{lemma}\label{prop:P-properties}
		Assume that
		(\ref{ass:kernel}) holds and fix
		\(\mu\geq0\). The family \(P_\mu\) defined by
		\eqref{eq:P-spectral-definition} satisfies
		\(P_\mu(0)=I\)
		and \(P_\mu(\cdot)x\in C([0,\infty);X)\) for every \(x\in X\).
		Moreover, for  \(t\geq0\),
		\begin{equation}\label{eq:P-exact-norm}
			\|P_\mu(t)\|_{\mathcal L(X)}
			=p_{\mu,\lambda_1}(t)
			\leq m_\mu(t)\leq1.
		\end{equation}
	\end{lemma}
	
	\begin{proof}
		By Lemma~\ref{lem:pnu-properties}, extend \(p_{\mu,\nu}\)
		continuously to \(t=0\) by setting \(p_{\mu,\nu}(0):=1\).
		Then \eqref{eq:P-spectral-definition} gives
		\(P_\mu(0)=I\).
		Fix \(v=\sum_{n=1}^{\infty}v_n\varphi_n\in X\)
		and \(t_0\geq0\). By \eqref{eq:pnu-basic-bound}, we obtain
		\[
		|p_{\mu,\lambda_n}(t)-p_{\mu,\lambda_n}(t_0)|^2|v_n|^2
		\leq |v_n|^2\quad \text{and}
		\quad
		\sum_{n=1}^{\infty}|v_n|^2=\|v\|_X^2<\infty.
		\]
		Thus Lemma~\ref{lem:pnu-properties}, Parseval's identity
		and dominated convergence yield
		\[
		\begin{aligned}
			\|P_\mu(t)v-P_\mu(t_0)v\|_X^2
			&=
			\sum_{n=1}^{\infty}
			|p_{\mu,\lambda_n}(t)-p_{\mu,\lambda_n}(t_0)|^2
			|v_n|^2
			\to0,
			\qquad\text{as }t\to t_0.
		\end{aligned}
		\]
		Hence \(P_\mu(\cdot)v\in C([0,\infty);X)\).
		
		Moreover, Lemma~\ref{lem:pnu-properties} and
		Parseval's identity give
		\(
		\|P_\mu(t)v\|_X^2
		\leq
		p_{\mu,\lambda_1}(t)^2\|v\|_X^2\)
		for \( t\geq0.
		\)
		Since
		\(P_\mu(t)\varphi_1=p_{\mu,\lambda_1}(t)\varphi_1\)
		and \(\|\varphi_1\|_X=1\),
		\[
		\begin{aligned}
			\|P_\mu(t)\|_{\mathcal L(X)}
			\geq
			\frac{\|P_\mu(t)\varphi_1\|_X}{\|\varphi_1\|_X}
			=
			\frac{p_{\mu,\lambda_1}(t)\|\varphi_1\|_X}
			{\|\varphi_1\|_X}
			=
			p_{\mu,\lambda_1}(t).
		\end{aligned}
		\]
		Therefore,
		\(
		\|P_\mu(t)\|_{\mathcal L(X)}
		=
		p_{\mu,\lambda_1}(t).
		\)	Together with \eqref{eq:pnu-basic-bound}, this proves
		\eqref{eq:P-exact-norm}.
	\end{proof}		
	
	For the pathwise treatment of the fractional Brownian noise, we impose the following additional condition on
	the kernel \(k\).
	
	\namedassumption{HB}{ass:kernel-holder}
	There exist \(\eta\in(0,1]\) and constant \(C>0\) such that
	\begin{equation*}
		\int_0^t k(r)\,dr
		\leq C
		t^\eta,
		\qquad 0<t\leq1.
	\end{equation*}

	\begin{lemma}\label{prop:P-holder}
		Assum that (\ref{ass:kernel}) and
		(\ref{ass:kernel-holder}) hold and fix \(\mu\geq0\).
		Let \(\delta\in(0,1]\) and set
		\(
		\rho:=\min\{\delta,\eta\}.
		\)
		Then there exists a constant \(C>0\) such that
		\[
		\|P_\mu(t)-P_\mu(s)\|_{\mathcal L(X_\delta,X)}
		\leq C|t-s|^\rho,
		\qquad s,t\geq0.
		\]
	\end{lemma}
	\begin{proof}
		By Lemma~\ref{lem:pnu-properties} and Bernstein's theorem,
		there exists a probability measure \(\varpi_{\mu,\nu}\)
		on \([0,\infty)\) such that
		\begin{equation}\label{eq-probability measure}
			p_{\mu,\nu}(t)
			=
			\int_{0}^{\infty}
			e^{-rt}\,\varpi_{\mu,\nu}(dr),
			\qquad t\geq0.
		\end{equation}
		For \(0\leq s<t\) and \(h=t-s\), we get
		\[
		\begin{aligned}
			0
			\leq
			p_{\mu,\nu}(s)-p_{\mu,\nu}(t)
			&\leq
			\int_{0}^{\infty}
			\bigl(1-e^{-rh}\bigr)\,\varpi_{\mu,\nu}(dr)=
			1-p_{\mu,\nu}(h).
		\end{aligned}
		\]
		By \eqref{eq:m-k}, \eqref{eq:pnu-basic-bound} and
		(\ref{ass:kernel-holder}), for \(0<h\leq1\), one has
		\[
		0\leq1-m_\mu(h)
		\leq
		\mu\int_0^h k(r)\,dr
		\leq Ch^\eta,\]
		and
		\[
		0\leq\nu(m_\mu*p_{\mu,\nu})(h)
		=
		m_\mu(h)-p_{\mu,\nu}(h)
		\leq\min\{1,\nu h\}.
		\]
		Consequently,
		\[
		\begin{aligned}
			1-p_{\mu,\nu}(h)
			&=
			1-m_\mu(h)+\nu(m_\mu*p_{\mu,\nu})(h)
			\leq Ch^\eta+\min\{1,\nu h\}.
		\end{aligned}
		\]
		Using
		\(
		\min\{1,a\}\leq a^\delta\) with
		\(a\geq0\) and \( 0<\delta\leq1,
		\)
		we obtain, for \(\nu\geq\lambda_1\) and \(0<h\leq1\),
		\[
		\begin{aligned}
			\nu^{-\delta}
			|p_{\mu,\nu}(t)-p_{\mu,\nu}(s)|
			&\leq
			C\lambda_1^{-\delta}h^\eta
			+\nu^{-\delta}(\nu h)^\delta
			=
			C\lambda_1^{-\delta}h^\eta+h^\delta
			\leq Ch^\rho.
		\end{aligned}
		\]
		For \(h\geq1\), \eqref{eq:pnu-basic-bound} gives
		\[
		\nu^{-\delta}
		|p_{\mu,\nu}(t)-p_{\mu,\nu}(s)|
		\leq\lambda_1^{-\delta}
		\leq\lambda_1^{-\delta}h^\rho.
		\]
		Thus, we have
		\begin{equation}\label{eq:scalar-holder}
			\sup_{\nu\geq\lambda_1}
			\nu^{-\delta}
			|p_{\mu,\nu}(t)-p_{\mu,\nu}(s)|
			\leq C|t-s|^\rho,
			\qquad s,t\geq0.
		\end{equation}
		
		For \(v=\sum_{n=1}^\infty v_n\varphi_n\in X_\delta\),
		\eqref{eq:P-spectral-definition} and
		\eqref{eq:scalar-holder} yield
		\[
		\begin{aligned}
			\|(P_\mu(t)-P_\mu(s))v\|_X^2
			\leq
			C|t-s|^{2\rho}
			\sum_{n=1}^{\infty}
			\lambda_n^{2\delta}|v_n|^2=
			C|t-s|^{2\rho}\|v\|_\delta^2.
		\end{aligned}
		\]
		The conclusion follows.
	\end{proof}

	\namedassumption{HC}{ass:noise-coefficient}
	Assume that there exists
	\(\delta\in(0,1]\) such that
	\(
	\rho:=\min\{\delta,\eta\}>1-H.
	\)
	Moreover,  operator \(
	G\in\mathcal L(X,X_\delta)
	\).
	
	Fix \(\mu\geq0\). Lemma~\ref{prop:P-holder} and  (\ref{ass:noise-coefficient})
	give
	\begin{equation}\label{eq:PGH}
		\begin{aligned}
			\|(P_\mu(t)-P_\mu(s))G\|_{\mathcal L(X)}
			&\leq
			\|P_\mu(t)-P_\mu(s)\|_{\mathcal L(X_\delta,X)}
			\|G\|_{\mathcal L(X,X_\delta)}
			\leq
			C|t-s|^\rho
			\|G\|_{\mathcal L(X,X_\delta)}.
		\end{aligned}
	\end{equation}
	Thus
	\(
	P_\mu(\cdot)G
	\in
	C^\rho
	([0,\infty);\mathcal L(X)).
	\)
	By (\ref{ass:noise-coefficient}), we  fix
	\(
	\gamma\in(\max\{1/2,1-\rho\},H).
	\)
	Then
	\(
	\rho+\gamma>1.
	\)

	\subsection{Fractional Brownian motion}
	\label{subsec: fractional Brownian motion}
	
	We recall that a two-sided real-valued fractional Brownian motion $\widetilde{\beta}^H(\cdot)$ with a
	Hurst index $H\in(0,1)$ is a centered Gaussian process with covariance function
	\begin{equation}\label{eq:covariance structure}
		\mathbb{E}\big(\widetilde{\beta}^H(t)\widetilde{\beta}^H(s)\big)
		=\frac12\big(|t|^{2H}+|s|^{2H}-|t-s|^{2H}\big),\quad \text{for } s,t\in\mathbb{R}.
	\end{equation}
	Let \(H\in(1/2,1)\), and let
	\(
	Q\in\mathcal L(X)
	\)
	be a nonnegative self-adjoint trace-class operator. Then there exist
	an orthonormal basis \(\{e_j\}_{j\geq1}\) of \(X\) and a sequence
	\(\{q_j\}_{j\geq1}\subset[0,\infty)\) such that
	\(
	Qe_j=q_je_j\)
	and
	\(\sum_{j=1}^{\infty}q_j
	=
	\operatorname{Tr}Q
	<
	\infty.
	\)
	Let \(\{\beta_j^H\}_{j\geq1}\) be a sequence of independent two-sided
	real-valued fractional Brownian motions with Hurst parameter \(H\).
	An \(X\)-valued \(Q\)-fractional Brownian motion is represented by
	\begin{equation}\label{eq:Q-fbm-representation}
		B^H(t)
		=
		\sum_{j=1}^{\infty}
		\sqrt{q_j}\,
		\beta_j^H(t)e_j,
		\qquad
		t\in\mathbb R.
	\end{equation}
	Since \(Q\) is trace class, the series in
	\eqref{eq:Q-fbm-representation} converges in mean square in \(X\) for each
	\(t\in\mathbb R\), and \(B^H\) admits a continuous version.
	Moreover, for every \(\gamma\in(0,H)\), its sample paths are locally
	\(\gamma\)-Hölder continuous almost surely.

	Following \cite{chen} and \cite{gao},
	let \(C_0(\mathbb R;X)\) denote the space of continuous
	\(X\)-valued paths on \(\mathbb R\) which are zero at zero.
	We endow this space with the compact-open topology.
	Let
	\(
	\mathcal F:=\mathcal B(C_0(\mathbb R;X)),
	\)
	and let \(\mathbb P\) be the law on this space of  the \(Q\)-fractional Brownian motion constructed above.
	On this probability space, define 
	\begin{equation}\label{eq:canonical-fbm}
		B^H(t,\omega):=\omega(t),
		\qquad t\in\mathbb R,\quad
		\omega\in C_0(\mathbb R;X).
	\end{equation}
	Under \(\mathbb P\), this process is an \(X\)-valued
	\(Q\)-fractional Brownian motion with Hurst parameter \(H\).
	For \(t\in\mathbb R\), define the Wiener shift by
	\begin{equation}\label{eq:Wiener shift}
		(\theta_t\omega)(r)
		:=\omega(t+r)-\omega(t),
		\qquad r\in\mathbb R.
	\end{equation}
	These shifts satisfy
	\[
	\theta_0=I,\quad 
	\theta_{t+s}=\theta_t\circ\theta_s
	\quad \text {for}\quad s,t\in\mathbb R.\]
	The map \((t,\omega)\to\theta_t\omega\) is jointly
	measurable, and the stationary increments of fractional
	Brownian motion imply that
	\(
	\mathbb P\circ\theta_t^{-1}=\mathbb P\)
	for \(t\in\mathbb R.
	\)
	Moreover,
	\(
	\bigl(C_0(\mathbb R;X),\mathcal F,\mathbb P,\theta\bigr)
	\)
	is an ergodic metric dynamical system.
	
	Recall the exponents \(\rho=\min\{\delta,\eta\}\) and
	\(\gamma\in(\max\{1/2,1-\rho\},H)\) fixed after
	(\ref{ass:noise-coefficient}).
	As in Lemma~15 of \cite{chen}, define
	\[
	\mathcal O
	:=
	\left\{
	\omega\in C_0(\mathbb R;X):
	[\omega]_{\gamma;[-N,N]}<\infty
	\text{ for every }N\in\mathbb N
	\right\},
	\]
	where
	\begin{equation}\label{eq:local Holder}
		[\omega]_{\gamma;[-N,N]}
		:=
		\sup_{-N\leq s<t\leq N}
		\frac{\|\omega(t)-\omega(s)\|_X}
		{|t-s|^\gamma}.
	\end{equation}
	By the standard sample-path regularity of fractional Brownian motion, the set \(\mathcal O\) is measurable,
	\(
	\mathbb P(\mathcal O)=1
	\)
	and 
	\(
	\theta_t\mathcal O=\mathcal O\) for \(t\in\mathbb R\).
	In what follows, all pathwise statements are understood for
	\(\omega\in\mathcal O\). 
	On \(\mathcal O\), we use the trace \(\sigma\)-algebra
	\(\mathcal F|_{\mathcal O}\) and the restriction of
	\(\mathbb P\), which we denote again by \(\mathcal F\)
	and \(\mathbb P\), respectively.
	
	Fixed \(\omega\in\mathcal O\), in view of \eqref{eq:canonical-fbm}, the stochastic
	convolution is  interpreted pathwise as
	\begin{equation*}
		\int_0^t P_\mu(t-q)G\,dB^H(q,\omega)
		:=
		\int_0^t P_\mu(t-q)G\,d\omega(q)
		\qquad t\geq0.
	\end{equation*}
	The formal noise term \(G\,dB^H(t)/dt\) in \eqref{e1.1}
	is understood through the following definition.
	
	\begin{definition}\label{def:pathwise-mild-solution}
		Let \(T>0\), \(x\in X\), and \(\omega\in\mathcal O\).
		A function
		\(
		u(\cdot,\omega,x)\in C([0,T];X)
		\)
		is called a pathwise mild solution of \eqref{e1.1} on \([0,T]\) if, for \(t\in[0,T]\),
		\begin{equation}\label{eq:mild-stochastic}
			\begin{aligned}
				u_\mu(t,\omega,x)
				={}
				P_\mu(t)x+ \int_0^t P_\mu(t-q)f(u_\mu(q,\omega,x))\,dq 
				+
				\int_0^t
				P_\mu(t-q)G\,d\omega(q)
				.
			\end{aligned}
		\end{equation}
	\end{definition}
	More generally, for a given function \(\Psi\in C([0,\infty);X)\),
	we consider the Volterra integral equation
	\begin{equation}\label{eq:generalized-profile-equation}
		\begin{aligned}
			u_\mu(t,\omega,\Psi)
			= \Psi(t)
			+ \int_0^t P_\mu(t-q)f(u_\mu(q,\omega,\Psi))\,dq 
			+ \int_0^t P_\mu(t-q)G\,d\omega(q),
			\qquad t\geq0.
		\end{aligned}
	\end{equation}
	Taking \(\Psi=P_\mu(\cdot)x\) in
	\eqref{eq:generalized-profile-equation} gives
	\eqref{eq:mild-stochastic}.

	\subsection{Young integration and the stochastic convolution}
	\label{subsec:young-integration}
	
	We recall the classical Young integration theorem.
	
	\begin{lemma}\label{lem:young-integration}\cite{young}
		Let \(E_1\) and \(E_2\) be Banach spaces, and let
		\(
		K\in C^\rho([a,b];\mathcal L(E_1,E_2))\)
		and \(	
		\omega\in C^\gamma([a,b];E_1),
		\)
		where
		\(
		\rho+\gamma>1.
		\)
		Then the Young integral
		\(
		\int_a^b K(q)\,d\omega(q)
		\)
		is well defined in \(E_2\). Moreover, there exists a constant
		\(C>0\) such that
		\begin{equation}\label{eq:young-estimate}
			\begin{aligned}
				\left\|
				\int_a^b K(q)\,d\omega(q)
				\right\|_{E_2}
				\leq
				\|K(a)\|_{\mathcal L(E_1,E_2)}
				[\omega]_{\gamma;[a,b]}
				|b-a|^\gamma
				+
				C
				[K]_{\rho;[a,b]}
				[\omega]_{\gamma;[a,b]}
				|b-a|^{\rho+\gamma},
			\end{aligned}
		\end{equation}
		where
		\[
		[K]_{\rho;[a,b]}
		:=
		\sup_{\substack{a\leq s<t\leq b}}
		\frac{
			\|K(t)-K(s)\|_{\mathcal L(E_1,E_2)}
		}{
			|t-s|^\rho
		}.
		\]
		Furthermore, 
		for \(s\in\mathbb R\) and \(t\geq0\), the Young integral is compatible with the Wiener shift
		\begin{equation}\label{eq:young-shift}
			\int_s^{s+t}K(q)\,d\omega(q)
			=
			\int_0^t
			K(s+r)\,d(\theta_s\omega)(r),
		\end{equation}
		whenever the integrals involved are well defined.
	\end{lemma}

	We now introduce a two-parameter stochastic convolution. For
	\(t,r\geq0\) and \(\omega\in\mathcal O\), define
	\begin{equation}\label{eq:two-parameter-convolution}
		Z_\mu(t,r,\omega)
		:=
		\int_0^t P_\mu(t+r-q)G\,d\omega(q).
	\end{equation}
	For \(r=0\), we write
	\(	Z_\mu(t,0,\omega)
	=
	\int_0^t P_\mu(t-q)G\,d\omega(q).\)
	\begin{lemma}\label{lem:stochastic-convolution}
		Assump that (\ref{ass:kernel}),
		(\ref{ass:kernel-holder})  and
		(\ref{ass:noise-coefficient}) hold. Let  \(\rho+\gamma>1\) and \(t,r\geq0\).
		Fixed \(\omega\in\mathcal O\), the Young integral \eqref{eq:two-parameter-convolution}
		is well defined and   the following properties hold.
		
		\smallskip
		\noindent\textup{(i)}
		For every \(t\geq0\),  
		the map \(r\mapsto Z_\mu(t,r,\omega)\) is continuous and bounded on \([0,\infty)\).
		
		\smallskip
		\noindent\textup{(ii)}
		For every \(T>0\),
		\begin{equation}\label{eq:convolution-uniform-continuity}
			\sup_{0\leq t\leq T-h}\sup_{r\geq0}
			\|Z_\mu(t+h,r,\omega)-Z_\mu(t,r,\omega)\|_X
			\to0
			\qquad\text{as }h\to 0.
		\end{equation}
		In particular, \(Z_\mu(\cdot,0,\omega)\in C([0,\infty);X)\).
		
		\smallskip
		\noindent\textup{(iii)}
		For all \(s,t,r\geq0\),
		\begin{equation}\label{eq:convolution-restart}
			Z_\mu(s+t,r,\omega)
			=
			Z_\mu(s,t+r,\omega)
			+
			Z_\mu(t,r,\theta_s\omega).
		\end{equation}
	\end{lemma}
	\begin{proof}
		Fix \(\omega\in\mathcal O\) and set
		\(Z_\mu(0,r,\omega)=0\) for \(r\geq0\).
		For \(t>0\), the embedding \(X_\delta\hookrightarrow X\),
		\eqref{eq:P-exact-norm} and \eqref{eq:PGH} give
		\begin{equation}\label{eq:kg}
			\sup_{q\in[0,t]}
			\|P_\mu(t+r-q)G\|_{\mathcal L(X)}
			+
			[P_\mu(t+r-\cdot)G]_{\rho;[0,t]}
			\leq
			C\|G\|_{\mathcal L(X,X_\delta)},
			\qquad r\geq0.
		\end{equation}
		Hence
		\(P_\mu(t+r-\cdot)G\in C^\rho([0,t];\mathcal L(X))\).
		Since \(\omega\in C^\gamma([0,t];X)\) and
		\(\rho+\gamma>1\),
		Lemma~\ref{lem:young-integration} implies that
		\eqref{eq:two-parameter-convolution} is well defined.
		
		\textup{(i)}
		Lemma~\ref{lem:young-integration} and \eqref{eq:kg} give
		\begin{equation}\label{eq:two-parameter-bound}
			\sup_{r\geq0}\|Z_\mu(t,r,\omega)\|_X
			\leq
			C\|G\|_{\mathcal L(X,X_\delta)}
			[\omega]_{\gamma;[0,t]}
			\bigl(t^\gamma+t^{\rho+\gamma}\bigr).
		\end{equation}
		Fix \(t>0\), \(r,r_0\geq0\), and
		\(\beta\in(1-\gamma,\rho)\). Define
		\(
		D_{r,r_0}(q)
		:=
		[P_\mu(t+r-q)-P_\mu(t+r_0-q)]G
		\)
		for \( q\in[0,t]\).	By \eqref{eq:PGH}, we have
		\[
		\sup_{q\in[0,t]}
		\|D_{r,r_0}(q)\|_{\mathcal L(X)}
		\leq
		C\|G\|_{\mathcal L(X,X_\delta)}
		|r-r_0|^\rho,
		\]
		and
		\[
		\begin{aligned}
			\|D_{r,r_0}(q_1)-D_{r,r_0}(q_2)\|_{\mathcal L(X)}
			&\leq
			C\|G\|_{\mathcal L(X,X_\delta)}
			\min\{|r-r_0|^\rho,|q_1-q_2|^\rho\}\\
			&\leq
			C\|G\|_{\mathcal L(X,X_\delta)}
			|r-r_0|^{\rho-\beta}|q_1-q_2|^\beta,
		\end{aligned}
		\]
		for \(0\leq q_1, q_2\leq t\).	Hence,
		\(
		[D_{r,r_0}]_{\beta;[0,t]}
		\leq
		C\|G\|_{\mathcal L(X,X_\delta)}
		|r-r_0|^{\rho-\beta}.
		\)
		Thus \(D_{r,r_0}\in C^\beta([0,t];\mathcal L(X))\).
		
		Since \(\beta+\gamma>1\),
		Lemma~\ref{lem:young-integration} yields	
		\begin{equation}\label{eq-zr}
			\begin{aligned}
				\|Z_\mu(t,r,\omega)-Z_\mu(t,r_0,\omega)\|_X
				&\leq
				C\|G\|_{\mathcal L(X,X_\delta)}
				[\omega]_{\gamma;[0,t]}
				\left(
				|r-r_0|^\rho t^\gamma
				+
				|r-r_0|^{\rho-\beta}t^{\beta+\gamma}
				\right)
				\to0,
			\end{aligned}
		\end{equation}
		as \(r\to r_0\).
		Together with \eqref{eq:two-parameter-bound} and
		\(Z_\mu(0,r,\omega)=0\), this proves \textup{(i)}.
		
		\smallskip
		\noindent\textup{(ii)}
		For \(0<h\leq T\) and \(0\leq t\leq T-h\),
		\[
		\begin{aligned}
			Z_\mu(t+h,r,\omega)-Z_\mu(t,r,\omega)
			=
			\int_t^{t+h}
			P_\mu(t+h+r-q)G\,d\omega(q)
			+
			Z_\mu(t,r+h,\omega)-Z_\mu(t,r,\omega).
		\end{aligned}
		\]
		By \eqref{eq:P-exact-norm}, \eqref{eq:PGH}
		and Lemma~\ref{lem:young-integration},
		\[
		\begin{aligned}
			&\left\|
			\int_t^{t+h}P_\mu(t+h+r-q)G\,d\omega(q)
			\right\|_X
			\leq
			C\|G\|_{\mathcal L(X,X_\delta)}
			[\omega]_{\gamma;[0,T]}
			\bigl(h^\gamma+h^{\rho+\gamma}\bigr).
		\end{aligned}
		\]
		Combining this with \eqref{eq-zr} gives
		\begin{equation*}
			\begin{aligned}
				&\sup_{0\leq t\leq T-h}\sup_{r\geq0}
				\|Z_\mu(t+h,r,\omega)-Z_\mu(t,r,\omega)\|_X\\
				&\quad\leq
				C\|G\|_{\mathcal L(X,X_\delta)}
				[\omega]_{\gamma;[0,T]}
				\left(
				h^\gamma+h^{\rho+\gamma}
				+h^\rho T^\gamma
				+h^{\rho-\beta}T^{\beta+\gamma}
				\right)
				\to0,
			\end{aligned}
		\end{equation*}
		as \(h\to0\).
		Thus, \eqref{eq:convolution-uniform-continuity} holds. Since \(T>0\) is arbitrary,
		\(Z_\mu(\cdot,0,\omega)\in C([0,\infty);X)\).

		\smallskip
		\noindent\textup{(iii)}
		By \eqref{eq:two-parameter-convolution} and
		\eqref{eq:young-shift}, for \(s,t,r\geq0\),
		\[
		\begin{aligned}
			Z_\mu(s+t,r,\omega)
			={}&
			\int_0^s P_\mu(s+t+r-q)G\,d\omega(q)
			+
			\int_s^{s+t}P_\mu(s+t+r-q)G\,d\omega(q)\\
			={}&
			Z_\mu(s,t+r,\omega)
			+
			\int_0^t P_\mu(t+r-q)G\,d(\theta_s\omega)(q)
			=
			Z_\mu(s,t+r,\omega)+Z_\mu(t,r,\theta_s\omega).
		\end{aligned}
		\]
		This proves \eqref{eq:convolution-restart}.
	\end{proof}

	\section{Existence and uniqueness of pathwise mild solutions}
	\label{sec:pathwise-existence-uniqueness}
	In this section, we establish existence and uniqueness for a
	Volterra integral equation with an arbitrary continuous given function.
	Throughout this section, \(\mu\geq0\) and 
	\(\omega\in\mathcal O\) are fixed.
	
	\namedassumption{HL}{ass:nonlinearity}
	The   \(f:X\to X\) is globally Lipschitz continuous which there exists
	a constant \(L >0\) such that
	\begin{equation}\label{eq:f-lipschitz}
		\|f(x)-f(y)\|_X\leq L \|x-y\|_X,
		\qquad x,y\in X.
	\end{equation}
	
	\begin{theorem}\label{thm:generalized-profile-solution}
		Suppose that
		Lemma~\ref{lem:stochastic-convolution} and
		(\ref{ass:nonlinearity}) hold.
		For every \(\Psi\in C([0,\infty);X)\) and
		\(\omega\in\mathcal O\), equation
		\eqref{eq:generalized-profile-equation}
		admits a unique global solution
		\[
		u_\mu (\cdot,\omega,\Psi)\in C([0,\infty);X).
		\]
	\end{theorem}
	
	\begin{proof}
		Fix \(\Psi\in C([0,\infty);X)\) and \(\omega\in\mathcal O\).
		Set
		\begin{equation}\label{eq:xi-omega}
			\phi^*(t):=\Psi(t)+Z_\mu(t,0,\omega),
			\qquad t\geq0.
		\end{equation}
		By Lemma~\ref{lem:stochastic-convolution},
		\(\phi^*\in C([0,\infty);X)\).
		Choose \(\kappa>L\) and define
		\begin{equation}\label{eq:auxiliary-weight-definition}
			\chi(t):=
			e^{\kappa t}
			\left(1+\sup_{0\leq r\leq t}\|\phi^*(r)\|_X\right),
			\qquad t\geq0.
		\end{equation}
		Then \(\chi\) is continuous and nondecreasing, with
		\(\chi(t)\geq1\).
		Consider
		\begin{equation}\label{eq:global-weighted-space}
			\mathcal X_\kappa
			:=
			\left\{
			u\in C([0,\infty);X):\|u\|_\kappa<\infty
			\right\},
			\qquad
			\|u\|_\kappa
			:=
			\sup_{t\geq0}\frac{\|u(t)\|_X}{\chi(t)}.
		\end{equation}
		
		Let  \(\{u_n\}_{n\geq1}\) be a Cauchy sequence in \(\mathcal X_\kappa\).
		For every \(T>0\),
		\[
		\sup_{0\leq t\leq T}\|u_n(t)-u_m(t)\|_X
		\leq
		\chi(T)\|u_n-u_m\|_\kappa
		\to0,\qquad
		\text{as }m,n\to\infty.
		\]
		By completeness of \(C([0,T];X)\),
		\(u_n\) converges locally uniformly to some
		\(u\in C([0,\infty);X)\).
		For every \(\varepsilon>0\), choose \(N\) such that
		\(
		\|u_n-u_m\|_\kappa\leq\varepsilon\) 
		for \( n,m\geq N.
		\)
		Letting \(m\to\infty\) pointwise gives
		\(
		\|u_n(t)-u(t)\|_X\leq\varepsilon\chi(t)\) for
		\(t\geq0\) and \(n\geq N.
		\)
		Hence
		\(
		\|u_n-u\|_\kappa\leq\varepsilon\)
		and
		\(\|u\|_\kappa
		\leq
		\|u-u_N\|_\kappa+\|u_N\|_\kappa<\infty.
		\)
		Thus \(\mathcal X_\kappa\) is a Banach space.

		By \eqref{eq:auxiliary-weight-definition},
		\[
		\frac{\|\phi^*(t)\|_X}{\chi(t)}
		=
		e^{-\kappa t}
		\frac{\|\phi^*(t)\|_X}
		{1+\sup_{0\leq r\leq t}\|\phi^*(r)\|_X}
		\leq e^{-\kappa t}\leq1.
		\]
		Consequently,
		\(\phi^*\in\mathcal X_\kappa\) and
		\(\|\phi^*\|_\kappa\leq1\).
		For \(0\leq q\leq t\),
		\[
		\frac{\chi(q)}{\chi(t)}
		=
		e^{-\kappa(t-q)}
		\frac{1+\sup_{0\leq r\leq q}\|\phi^*(r)\|_X}
		{1+\sup_{0\leq r\leq t}\|\phi^*(r)\|_X}
		\leq e^{-\kappa(t-q)}.
		\]
		Therefore,
		\begin{equation}\label{eq:weight-integral-estimate}
			\frac1{\chi(t)}\int_0^t\chi(q)\,dq
			\leq
			\int_0^t e^{-\kappa(t-q)}\,dq
			=
			\frac{1-e^{-\kappa t}}{\kappa}
			\leq\frac1\kappa.
		\end{equation}
		Since  \(\chi(q)\geq1\),
		\begin{equation}\label{eq:constant-weight-bound}
			\frac{t}{\chi(t)}
			\leq
			\frac1{\chi(t)}\int_0^t\chi(q)\,dq
			\leq\frac1\kappa.
		\end{equation}
		
		Define
		\begin{equation}\label{eq:fixed-point-operator}
			(\Phi u)(t)
			:=
			\phi^*(t)
			+
			\int_0^t P_\mu(t-q)f(u(q))\,dq,
			\qquad t\geq0.
		\end{equation}
		For \(u\in\mathcal X_\kappa\),
		\(f(u(\cdot))\in C([0,\infty);X)\).
		By Lemma~\ref{prop:P-properties} and
		\eqref{eq:P-exact-norm}, the \(P_\mu(t-q)f(u(q))\) in
		\eqref{eq:fixed-point-operator} is continuous on \([0,t]\).
		Hence the Bochner integral \(\int_0^t P_\mu(t-q)f(u(q))\,dq\) is well defined.
		For \(t>0\), the substitution \(q=ta\) gives
		\begin{equation}\label{eq:fixed-point-con}
			\int_0^t P_\mu(t-q)f(u(q))\,dq
			=
			t\int_0^1 P_\mu(t(1-a))f(u(ta))\,da.
		\end{equation}
		The \eqref{eq:fixed-point-con} also holds at \(t=0\), since both sides vanish.
		Fix \(T>0\) and \(t_2\in[0,T]\), and let \(t_1\to t_2\) in \([0,T]\).
		By Lemma~\ref{prop:P-properties},
		\eqref{eq:P-exact-norm} and \eqref{eq:f-lipschitz},
		for each \(a\in[0,1]\),
		\begin{equation*}
			\begin{aligned}
				&\|P_\mu(t_1(1-a))f(u(t_1a))
				-P_\mu(t_2(1-a))f(u(t_2a))\|_X\\
				&\quad\leq L\|u(t_1a)-u(t_2a)\|_X
				+\|[P_\mu(t_1(1-a))-P_\mu(t_2(1-a))]f(u(t_2a))\|_X
				\to 0.
			\end{aligned}
		\end{equation*}
		Moreover, \eqref{eq:P-exact-norm} gives
		\(
		\|P_\mu(t_1(1-a))f(u(t_1a))\|_X
		\leq \sup_{0\leq q\leq T}\|f(u(q))\|_X<\infty .
		\)
		By dominated convergence on \([0,1]\),
		\begin{equation*}
			\int_0^1 P_\mu(t_1(1-a))f(u(t_1a))\,da
			\to
			\int_0^1 P_\mu(t_2(1-a))f(u(t_2a))\,da
			\quad\text{in }X .
		\end{equation*}
		Since \(t_1\to t_2\), the change of variables in
		\eqref{eq:fixed-point-con} yields
		\begin{equation*}
			\int_0^{t_1} P_\mu(t_1-q)f(u(q))\,dq
			\to
			\int_0^{t_2} P_\mu(t_2-q)f(u(q))\,dq
			\quad\text{in }X .
		\end{equation*}
		Thus \(	\int_0^t P_\mu(t-q)f(u(q))\,dq\) 
		is continuous on \([0,\infty)\). Since
		\(\phi^*\in C([0,\infty);X)\), it follows that
		\(\Phi u\in C([0,\infty);X)\).

		By (\ref{ass:nonlinearity}) and
		\eqref{eq:P-exact-norm}, we get
		\[
		\begin{aligned}
			\frac{\|(\Phi u)(t)\|_X}{\chi(t)}
			&\leq
			1+\frac{t\|f(0)\|_X}{\chi(t)}
			+
			L\|u\|_\kappa
			\frac1{\chi(t)}\int_0^t\chi(q)\,dq.
		\end{aligned}
		\]
		Using \eqref{eq:weight-integral-estimate} and
		\eqref{eq:constant-weight-bound}, we have
		\begin{equation}\label{eq:self-map-weighted-bound}
			\|\Phi u\|_\kappa
			\leq
			1+\frac{\|f(0)\|_X}{\kappa}
			+\frac L\kappa\|u\|_\kappa<\infty.
		\end{equation}
		Thus \(\Phi:\mathcal X_\kappa\to\mathcal X_\kappa\).
		
		For \(u_1,u_2\in\mathcal X_\kappa\),
		(\ref{ass:nonlinearity}),
		\eqref{eq:P-exact-norm} and
		\eqref{eq:weight-integral-estimate} give
		\[
		\begin{aligned}
			\frac{\|(\Phi u_1)(t)-(\Phi u_2)(t)\|_X}{\chi(t)}
			&\leq
			L\|u_1-u_2\|_\kappa
			\frac1{\chi(t)}\int_0^t\chi(q)\,dq
			\leq
			\frac L\kappa\|u_1-u_2\|_\kappa.
		\end{aligned}
		\]
		Taking the supremum over \(t\geq0\), we obtain
		\begin{equation}\label{eq:global-weighted-contraction}
			\|\Phi u_1-\Phi u_2\|_\kappa
			\leq
			\frac L\kappa\|u_1-u_2\|_\kappa,
			\qquad \frac L\kappa<1.
		\end{equation}
		By the Banach fixed point theorem, \(\Phi\) has a unique
		fixed point in \(\mathcal X_\kappa\), denoted by
		\(u_\mu (\cdot,\omega,\Psi)\). Thus
		\(
		u_\mu (\cdot,\omega,\Psi)
		=
		\Phi\bigl(u_\mu (\cdot,\omega,\Psi)\bigr),
		\)
		which is a global solution of
		\eqref{eq:generalized-profile-equation}.

		To prove uniqueness in \(C([0,\infty);X)\), let
		\(\widetilde u_1,\widetilde u_2\) be two solutions corresponding
		to the same \(\Psi\) and \(\omega\).
		By \eqref{eq:P-exact-norm} and
		(\ref{ass:nonlinearity}),
		\[
		\|\widetilde u_1(t)-\widetilde u_2(t)\|_X
		\leq \int_0^t 
		\bigr\|P_\mu(t-q)\bigl[f(\widetilde u_1(q))-f(\widetilde u_2(q))\bigr]\bigr\|_X\,dq\leq 
		L\int_0^t
		\|\widetilde u_1(q)-\widetilde u_2(q)\|_X\,dq.
		\]
		Gronwall's inequality on each finite interval \([0,T]\)
		yields \(\widetilde u_1=\widetilde u_2\) on \([0,\infty)\).
	\end{proof}				
	
	\section{The cocycle and random dynamical system}
	\label{sec:lifted-cocycle}	
	For \(\mu>0\), the family \(P_\mu(t)\) does not in general satisfy the
	semigroup property. For  \(\mu\geq0\), we introduce a suitable function space
	and construct a continuous random dynamical system on this
	space.  We recall the definition of a random dynamical system;
	see \cite[Definition~1.3(i)]{ge}.
	
	\begin{definition}\label{def:random-dynamical-system}
		Let \((\mathcal O,\mathcal F,\mathbb P,\theta)\) be the
		metric dynamical system introduced above  and let
		\(\mathcal Y\) be a separable Banach space.
		A mapping
		\(
		\varphi:[0,\infty)\times\mathcal O\times\mathcal Y
		\to \mathcal Y
		\)
		is called a random dynamical system on \(\mathcal Y\)
		over \(\theta\) if the following conditions hold.
		
		\smallskip
		\noindent\textup{(i)}
		The mapping \(\varphi\) is measurable with respect to
		\(
		\mathcal B([0,\infty))\otimes\mathcal F
		\otimes\mathcal B(\mathcal Y)\)
		and \(
		\mathcal B(\mathcal Y).
		\)
		
		\smallskip
		\noindent\textup{(ii)}
		For every \(\omega\in\mathcal O\) and \(\Psi\in\mathcal Y\),
		\(
		\varphi(0,\omega,\Psi)=\Psi.
		\)
		
		\smallskip
		\noindent\textup{(iii)}
		For all \(s,t\geq0\), \(\omega\in\mathcal O\), and
		\(\Psi\in\mathcal Y\), the cocycle identity holds:
		\[
		\varphi(s+t,\omega,\Psi)
		=
		\varphi\bigl(t,\theta_s\omega,\varphi(s,\omega,\Psi)\bigr).
		\]
		The random dynamical system is called continuous if,
		for every \(t\geq0\) and \(\omega\in\mathcal O\),
		the mapping
		\(
		\Psi\to \varphi(t,\omega,\Psi)
		\)
		is continuous on \(\mathcal Y\).
	\end{definition}
	
	\subsection{The phase space}
	\label{subsec:profile-space}
	
	For \(\Psi\in C([0,\infty);X)\), we introduce the weighted profile space
	\begin{equation*}
		\mathcal Y
		:=
		\left\{
		\Psi\in C([0,\infty);X):
		\|\Psi\|_{\mathcal Y}<\infty
		\right\},
	\end{equation*}
	where
	\begin{equation}\label{eq:profile-norm}
		\|\Psi\|_{\mathcal Y}
		:=
		\sum_{N=0}^{\infty}
		\frac{1}{2^N}\sup_{r\in[0,N]}\|\Psi(r)\|_X.
	\end{equation}
	For \(t\geq0\), define 
	\begin{equation}\label{eq:tauPsi}
		(\tau_t\Psi)(r):=\Psi(t+r),
		\qquad r\geq0.
	\end{equation}
	
	\begin{lemma}\label{lem:profile-space-estimates}
		The space \((\mathcal Y,\|\cdot\|_{\mathcal Y})\) is a separable Banach
		space with the following properties.
		
		\smallskip
		\noindent\textup{(i)}
		For every \(T\geq0\) and \(N\in\mathbb N\) with \(N\geq T\), we have 
		\begin{equation}\label{eq:profile-local-control}
			\sup_{0\leq r\leq T}\|\Psi(r)\|_X
			\leq
			2^{N}\|\Psi\|_{\mathcal Y},
			\qquad
			\Psi\in\mathcal Y.
		\end{equation}
		
		\smallskip
		\noindent\textup{(ii)}
		For every \(T\geq0\) and \(m\in\mathbb N\) with \(m\geq T\),
		\begin{equation}\label{eq:translation-bound}
			\sup_{0\leq t\leq T}
			\|\tau_t\Psi\|_{\mathcal Y}
			\leq
			2^{m} \|\Psi\|_{\mathcal Y},
			\qquad
			\Psi\in\mathcal Y.
		\end{equation}
		In particular, \(\tau_t\) is a bounded linear operator on
		\(\mathcal Y\) for \(t\geq0.\)

		\smallskip
		\noindent\textup{(iii)}
		If \(\Psi\in C([0,\infty);X)\) is bounded, then 
		\(\Psi \in\mathcal Y\) and
		\begin{equation}\label{eq:bounded estimate}
			\|\Psi\|_{\mathcal Y}
			\leq
			2\sup_{r\geq0}\|\Psi(r)\|_X.
		\end{equation}
	\end{lemma}

	\begin{proof}
		Let \(\{\Psi_j\}\) be a Cauchy
		sequence in \(\mathcal Y\). By \eqref{eq:profile-norm},
		for  \(N\in\mathbb N\),
		\[
		\sup_{0\leq r\leq N}
		\|\Psi_j(r)-\Psi_k(r)\|_X
		\leq
		2^N\|\Psi_j-\Psi_k\|_{\mathcal Y}
		\to0
		\qquad\text{as }j,k\to\infty.
		\]
		By completeness of \(C([0,N];X)\), there exists
		\(\Psi\in C([0,\infty);X)\) such that
		\(\Psi_j\to\Psi\) uniformly on every bounded interval.
		By Fatou's lemma applied to the
		series in \eqref{eq:profile-norm}, we have
		\[
		\|\Psi\|_{\mathcal Y}
		\leq
		\lim_{n\to\infty}
		\bigr(\inf_{j\geq n}\|\Psi_j\|_{\mathcal Y}\bigr)<\infty.
		\]
		Thus \(\Psi\in\mathcal Y\).
		Given \(\varepsilon>0\), choose \(J\) such that
		\(\|\Psi_j-\Psi_k\|_{\mathcal Y}\leq\varepsilon\)
		for all \(j,k\geq J\).
		For each fixed \(j\geq J\), the same application of Fatou's lemma gives
		\[
		\|\Psi_j-\Psi\|_{\mathcal Y}
		\leq
		\lim_{n\to\infty}
		\left(\inf_{k\geq n}\|\Psi_j-\Psi_k\|_{\mathcal Y}\right)
		\leq\varepsilon.
		\]
		Hence \(\Psi_j\to\Psi\) in \(\mathcal Y\), and \(\mathcal Y\) is complete.

		To prove separability, for \(m\in\mathbb N\), choose
		\(\chi_m\in C([0,\infty);[0,1])\) satisfying
		\(
		\chi_m=1 \) on \([0,m]\)
		and
		\(	\chi_m=0\) on \([m+1,\infty).
		\)
		By \eqref{eq:profile-norm}, for every \(\Psi\in\mathcal Y\),
		\begin{equation}\label{eq:profile-cutoff-approximation}
			\|\Psi-\chi_m\Psi\|_{\mathcal Y}
			\leq
			\sum_{N=m+1}^{\infty}2^{-N}
			\sup_{0\leq r\leq N}\|\Psi(r)\|_X
			\to0
			\qquad\text{as }m\to\infty.
		\end{equation}
		Set
		\(
		E_m:=
		\left\{
		\Psi\in\mathcal Y:
		\Psi(r)=0 \text{ for }r\geq m
		\right\}.
		\)
		For \(\Psi\in E_m\), \eqref{eq:profile-norm} gives
		\[
		2^{-m}\sup_{0\leq r\leq m}\|\Psi(r)\|_X
		\leq
		\|\Psi\|_{\mathcal Y}
		\leq
		2\sup_{0\leq r\leq m}\|\Psi(r)\|_X.
		\]
		Thus \(E_m\), viewed as the closed subspace
		\(
		\{u\in C([0,m];X):u(m)=0\},
		\)
		is separable because \(X\) is separable.
		Since \(\chi_m\Psi\in E_{m+1}\),
		\eqref{eq:profile-cutoff-approximation} shows that
		\(\bigcup_{m=1}^{\infty}E_m\) is dense in \(\mathcal Y\).
		Hence \(\mathcal Y\) is separable.
		
		For \textup{(i)}, fix \(T\geq0\) and choose \(N\in\mathbb N\)
		with \(N\geq T\). Then
		\[
		\sup_{0\leq r\leq T}\|\Psi(r)\|_X
		\leq\sup_{0\leq r\leq N}\|\Psi(r)\|_X
		\leq2^N\|\Psi\|_{\mathcal Y}.
		\]
		
		For \textup{(ii)}, choose \(m\in\mathbb N\) with \(m\geq T\).
		For every \(t\in[0,T]\), \eqref{eq:profile-norm} gives
		\[
		\begin{aligned}
			\|\tau_t\Psi\|_{\mathcal Y}
			\leq
			\sum_{N=0}^{\infty}2^{-N}
			\sup_{0\leq r\leq N+m}\|\Psi(r)\|_X
			=
			2^m\sum_{j=m}^{\infty}2^{-j}
			\sup_{0\leq r\leq j}\|\Psi(r)\|_X
			\leq
			2^m\|\Psi\|_{\mathcal Y}.
		\end{aligned}
		\]
		Taking the supremum over \(t\in[0,T]\) proves
		\eqref{eq:translation-bound}.

		Finally, if \(\Psi\in C([0,\infty);X)\) is bounded, then
		\[
		\|\Psi\|_{\mathcal Y}
		\leq
		\left(\sum_{N=0}^{\infty}2^{-N}\right)
		\sup_{r\geq0}\|\Psi(r)\|_X
		=2\sup_{r\geq0}\|\Psi(r)\|_X<\infty.
		\]
		Therefore \(\Psi\in\mathcal Y\), and
		\eqref{eq:bounded estimate} follows.
	\end{proof}
	
	Moreover, the following statement is a direct consequence of
	Lemma~\ref{lem:stochastic-convolution}.
	
	\begin{lemma}\label{lem:two-parameter-convolution}
		Under the assumptions of Lemma~\ref{lem:stochastic-convolution}, for
		each \(\omega\in\mathcal O\) the map
		\(t\mapsto Z_\mu(t,\cdot,\omega)\) is continuous from
		\([0,\infty)\) into \(\mathcal Y\).
	\end{lemma}
	\begin{proof}
		Fix \(\omega\in\mathcal O\).
		By Lemma~\ref{lem:stochastic-convolution} \textup{(i)},
		for every \(t\geq0\), the function
		\(r\to Z_\mu(t,r,\omega)\) is continuous and bounded.
		Hence \eqref{eq:bounded estimate} gives
		\begin{equation*}
			Z_\mu(t,\cdot,\omega)\in\mathcal Y,
			\qquad
			\|Z_\mu(t,\cdot,\omega)\|_{\mathcal Y}
			\leq
			2\sup_{r\geq0}\|Z_\mu(t,r,\omega)\|_X
			<\infty.
		\end{equation*}
		For \(s,t\geq0\), Lemma~\ref{lem:stochastic-convolution}\textup{(i)}
		gives
		\(
		Z_\mu(t,\cdot,\omega)-Z_\mu(s,\cdot,\omega)
		\in C_b([0,\infty);X).
		\)
		By \eqref{eq:bounded estimate} and
		\eqref{eq:convolution-uniform-continuity},
		\[
		\begin{aligned}
			\|Z_\mu(t,\cdot,\omega)-Z_\mu(s,\cdot,\omega)\|_{\mathcal Y}
			&\leq
			2\sup_{r\geq0}
			\|Z_\mu(t,r,\omega)-Z_\mu(s,r,\omega)\|_X
			\to0
			\qquad\text{as }t\to s.
		\end{aligned}
		\]
	\end{proof}

	\subsection{The cocycle property}
	\label{subsec:cocycle-property}
	
	Fix \(\mu\geq0\). For each \(\Psi\in\mathcal Y\) and
	\(\omega\in\mathcal O\), let \(u_\mu(\cdot,\omega,\Psi)\) denote the
	unique global solution of \eqref{eq:generalized-profile-equation}
	given by Theorem~\ref{thm:generalized-profile-solution}.
	Since \(Z_\mu(t,\omega)=Z_\mu(t,0,\omega)\), we have
	\begin{equation}\label{eq:uPsi}
		\begin{aligned}
			u_\mu(t,\omega,\Psi)
			=\Psi(t)
			+\int_0^t P_\mu(t-q)f\bigl(u_\mu(q,\omega,\Psi)\bigr)\,dq
			+Z_\mu(t,0,\omega),\qquad t\geq0.
		\end{aligned}
	\end{equation}
	For \(t\geq0\), define the function
	\(\Gamma_\mu(t,\omega,\Psi)\) by
	\begin{equation}\label{eq:lifted-map}
		\begin{aligned}
			\Gamma_\mu(t,\omega,\Psi)(r)
			:=\Psi(t+r)
			+\int_0^t P_\mu(t+r-q)f\bigl(u_\mu(q,\omega,\Psi)\bigr)\,dq
			+Z_\mu(t,r,\omega),\qquad r\geq0.
		\end{aligned}
	\end{equation}
	
	We verify that \(\Gamma_\mu(t,\omega,\Psi)\in\mathcal Y\).
	By \eqref{eq:tauPsi} and
	Lemma~\ref{lem:profile-space-estimates}\textup{(ii)},
	we derive \(\tau_t\Psi\in\mathcal Y\).
	Let \(r_n\to r\) in \([0,\infty)\).
	By Lemma~\ref{prop:P-properties}, \eqref{eq:P-exact-norm}
	and dominated convergence,
	\[
	\begin{aligned}
		&\left\|\int_0^t
		[P_\mu(t+r_n-q)-P_\mu(t+r-q)]
		f(u_\mu(q,\omega,\Psi))\,dq\right\|_X
		\to0,
	\end{aligned}
	\]
	Moreover, \eqref{eq:P-exact-norm} gives
	\[
	\sup_{r\geq0}
	\left\|\int_0^tP_\mu(t+r-q)f(u_\mu(q,\omega,\Psi))\,dq\right\|_X
	\leq
	\int_0^t\|f(u_\mu(q,\omega,\Psi))\|_X\,dq
	<\infty.
	\]
	Thus the integral in \eqref{eq:lifted-map}, as a function of \(r\),
	belongs to \(\mathcal Y\) by
	Lemma~\ref{lem:profile-space-estimates}\textup{(iii)}.
	Since \(Z_\mu(t,\cdot,\omega)\in\mathcal Y\) by
	Lemma~\ref{lem:two-parameter-convolution}, \eqref{eq:lifted-map}
	defines a mapping
	\(
	\Gamma_\mu:[0,\infty)\times\mathcal O\times\mathcal Y
	\to\mathcal Y.
	\)
	
	Using \eqref{eq:lifted-map} at \(r=0\) and 
	\eqref{eq:uPsi}, we obtain
	\begin{equation}\label{eq:gamma-evaluation-at-zero}
		\Gamma_\mu(t,\omega,\Psi)(0)=u_\mu(t,\omega,\Psi),
		\qquad t\geq0.
	\end{equation}
	For  \(x\in X\), 
	set
	\begin{equation}\label{eq:Psi-x}
		\Psi_{\mu,x}(r):=P_\mu(r)x,
		\qquad r\geq0.
	\end{equation}	
	Let \(\Psi_0:=\Psi_{\mu,0}\equiv0\).
	By Lemma~\ref{prop:P-properties} and \eqref{eq:P-exact-norm},
	\(\Psi_{\mu,x}\) is continuous and
	\begin{equation}\label{eq:Psi-x-bound} \Psi_{\mu,x}(0)=x,  \quad \text {and}\quad 
		\sup_{r\geq0}\|\Psi_{\mu,x}(r)\|_X\leq\|x\|_X.
	\end{equation}
	Hence Lemma~\ref{lem:profile-space-estimates}\textup{(iii)} gives
	\begin{equation}\label{eq:Psi-xy}
		\Psi_{\mu,x}\in\mathcal Y, \quad \text{and}\quad 	\|\Psi_{\mu,x}\|_{\mathcal Y}\leq2\|x\|_X.
	\end{equation}
	Taking \(\Psi=\Psi_{\mu,x}\) in \eqref{eq:uPsi} gives
	\eqref{eq:mild-stochastic} with initial value \(x\). By uniqueness and
	\eqref{eq:gamma-evaluation-at-zero},
	\[
	\Gamma_\mu(t,\omega,\Psi_{\mu,x})(0)
	=u_\mu(t,\omega,\Psi_{\mu,x})
	=u_\mu(t,\omega,x),
	\qquad t\geq0,\quad\omega\in\mathcal O.
	\]

	The following Lemmas are the main step in proving the cocycle property.
	\begin{lemma}\label{lem:solution-restart}
		Let \(\Psi\in\mathcal Y\), \(\omega\in\mathcal O\),
		and \(s,t\geq0\). Then,
		\begin{equation}\label{eq:solution-restart}
			u_{\Gamma_\mu(s,\omega,\Psi)}(t,\theta_s\omega)
			=u_\mu(s+t,\omega,\Psi).
		\end{equation}
	\end{lemma}
	\begin{proof}
		Fix \(s\geq0\), and set
		\(\Psi^{(s)}:=\Gamma_\mu(s,\omega,\Psi)\in\mathcal Y\).
		By \eqref{eq:uPsi},
		\begin{equation}\label{eq:uPsi-s}
			u_{\Psi^{(s)}}(t,\theta_s\omega)
			=
			\Psi^{(s)}(t)
			+\int_0^t P_\mu(t-q)
			f(u_{\Psi^{(s)}}(q,\theta_s\omega))\,dq
			+Z_\mu(t,0,\theta_s\omega).
		\end{equation}
		By \eqref{eq:uPsi}, \eqref{eq:lifted-map} and
		\eqref{eq:convolution-restart}, with \(q=s+\xi\),
		\begin{equation}\label{eq:uPsi-x}
			\begin{aligned}
				u_\mu(s+t,\omega,\Psi)
				={}&\Psi(s+t)
				+\int_0^s P_\mu(s+t-q)
				f(u_\mu(q,\omega,\Psi))\,dq\\
				&+\int_0^t P_\mu(t-\xi)
				f(u_{\Psi}(s+\xi,\omega))\,d\xi
				+Z_\mu(s+t,0,\omega)\\
				={}&\Psi^{(s)}(t)
				+\int_0^t P_\mu(t-\xi)
				f(u_{\Psi}(s+\xi,\omega))\,d\xi
				+Z_\mu(t,0,\theta_s\omega).
			\end{aligned}
		\end{equation}
		Since \(\theta_s\omega\in\mathcal O\), comparing
		\eqref{eq:uPsi-s} and \eqref{eq:uPsi-x}, uniqueness in
		Theorem~\ref{thm:generalized-profile-solution} gives
		\begin{equation*}
			u_\mu(s+t,\omega,\Psi)
			=
			u_\mu(t,\theta_s\omega,\Psi^{(s)}),
			\qquad t\geq0.
		\end{equation*}
		This proves \eqref{eq:solution-restart}.
	\end{proof}

	\begin{theorem}\label{thm:lifted-cocycle}
		For	\(\omega\in\mathcal O\) and
		\(\Psi\in\mathcal Y\),	the mappings \(\Gamma_\mu(t,\omega,\cdot)\)
		satisfy
		\(	\Gamma_\mu(0,\omega,\Psi)=\Psi\)
		and
		\begin{equation}\label{eq:lifted-cocycle}
			\Gamma_\mu(s+t,\omega,\Psi)
			=
			\Gamma_\mu\bigl(t,\theta_s\omega,\Gamma_\mu(s,\omega,\Psi)\bigr), \quad \text{for}\quad s,t\geq0.
		\end{equation}
	\end{theorem}
	
	\begin{proof}
		By \eqref{eq:lifted-map} and \(Z_\mu(0,r,\omega)=0\),
		we have \(\Gamma_\mu(0,\omega,\Psi)=\Psi\).
		Fix \(s,t\geq0\), \(\omega\in\mathcal O\), and
		\(\Psi\in\mathcal Y\). Set
		\(\Psi^{(s)}:=\Gamma_\mu(s,\omega,\Psi)\in\mathcal Y\).
		By \eqref{eq:solution-restart}, for every \(r\geq0\),
		\begin{equation}\label{eq:cocycle-inf}
			\begin{aligned}
				\int_0^t P_\mu(t+r-q)
				f\bigl(u_\mu(q,\theta_s\omega,\Psi^{(s)})\bigr)\,dq
				&=
				\int_0^t P_\mu(t+r-q)
				f\bigl(u_\mu(s+q,\omega,\Psi)\bigr)\,dq\\
				&
				=\int_s^{s+t}P_\mu(s+t+r-\xi)
				f\bigl(u_\mu(\xi,\omega,\Psi)\bigr)\,d\xi.
			\end{aligned}
		\end{equation}
		Moreover, \eqref{eq:convolution-restart} gives
		\begin{equation}\label{eq:cocycle-Z}
			Z_\mu(s,t+r,\omega)+Z_\mu(t,r,\theta_s\omega)
			=
			Z_\mu(s+t,r,\omega).
		\end{equation}
		
		Combining \eqref{eq:cocycle-inf} and \eqref{eq:cocycle-Z} with
		\eqref{eq:lifted-map}, we obtain
		\[
		\begin{aligned}
			&\Gamma_\mu(t,\theta_s\omega,\Psi^{(s)})(r)\\
			&=\Psi(s+t+r)
			+\int_0^s P_\mu(s+t+r-q)
			f\bigl(u_\mu(q,\omega,\Psi)\bigr)\,dq\\
			&\qquad+\int_s^{s+t}P_\mu(s+t+r-q)
			f\bigl(u_\mu(q,\omega,\Psi)\bigr)\,dq
			+Z_\mu(s+t,r,\omega)\\
			&=\Psi(s+t+r)
			+\int_0^{s+t}P_\mu(s+t+r-q)
			f\bigl(u_\mu(q,\omega,\Psi)\bigr)\,dq
			+Z_\mu(s+t,r,\omega)
			=\Gamma_\mu(s+t,\omega,\Psi)(r).
		\end{aligned}
		\]
		Since \(r\geq0\) is arbitrary, \eqref{eq:lifted-cocycle} follows.
	\end{proof}

	\subsection{Continuity and measurability of the   cocycle}
	\label{subsec:continuity-measurability}
	
	In this subsection, we establish the continuity and
	measurability of the  cocycle  mappings \(\Gamma_\mu\).
	
	\begin{lemma}\label{lem:tr-continuou}
		Fix \(\Psi\in\mathcal Y\) and \(\omega\in\mathcal O\), and define
		\[
		\mathscr F(t,r):=
		\int_0^t P_\mu(t+r-q)f(u_\mu(q,\omega,\Psi))\,dq,
		\qquad t,r\geq0.
		\]
		Then \((t,r)\mapsto \mathscr F(t,r)\) is  continuous from
		\([0,\infty)\times[0,\infty)\) into \(X\).
		Moreover, \(\mathscr F(t,\cdot)\in\mathcal Y\) 
		and \(t\mapsto \mathscr F(t,\cdot)\) is continuous   in  \(\mathcal Y\) for every \(t\geq0\).
	\end{lemma}

	\begin{proof}
		Fix \(T>0\). By Theorem~\ref{thm:generalized-profile-solution}
		and (\ref{ass:nonlinearity}), we get
		\(
		f(u_\mu(\cdot,\omega,\Psi))\in C([0,T];X)\)
		and 
		\(\sup_{0\leq q\leq T}\|f(u_\mu(q,\omega,\Psi))\|_X<\infty.
		\)
		By \eqref{eq:P-exact-norm}, for \(0\leq t\leq T\) and \(r\geq0\),
		\begin{equation}\label{eq:deterministic-integral-bound}
			\begin{aligned}
				\|\mathscr F(t,r)\|_X
				&\leq \int_0^t \|f\bigl(u_\mu(q,\omega,\Psi)\bigr)\|_X\,dq
				\leq T\sup_{0\leq q\leq T}\|f\bigl(u_\mu(q,\omega,\Psi)\bigr)\|_X
				<\infty .
			\end{aligned}
		\end{equation}
		The change of variables \(q=ta\) gives
		\begin{equation}\label{eq:deterministic-fixed-interval}
			\mathscr F(t,r)
			=t\int_0^1 P_\mu(t(1-a)+r)
			f(u_\mu(ta,\omega,\Psi))\,da,
			\qquad t,r\geq0.
		\end{equation}
		
		Let \((t_2,r_2)\to(t_1,r_1)\) in
		\([0,T]\times[0,\infty)\). By Lemma~\ref{prop:P-properties}
		and (\ref{ass:nonlinearity}), for \(a\in[0,1]\),
		\[
		\begin{aligned}
			&\bigl\|
			P_\mu(t_2(1-a)+r_2)f(u_\mu(t_2a,\omega,\Psi))
			-P_\mu(t_1(1-a)+r_1)f(u_\mu(t_1a,\omega,\Psi))
			\bigr\|_X\\
			&\quad\leq
			L\|u_\mu(t_2a,\omega,\Psi)-u_\mu(t_1a,\omega,\Psi)\|_X\\
			&\qquad+
			\bigl\|
			[P_\mu(t_2(1-a)+r_2)-P_\mu(t_1(1-a)+r_1)]
			f(u_\mu(t_1a,\omega,\Psi))
			\bigr\|_X
			\to0.
		\end{aligned}
		\]
		Moreover, \eqref{eq:P-exact-norm} gives
		\[
		\|P_\mu(t(1-a)+r)f(u_\mu(ta,\omega,\Psi))\|_X
		\leq\sup_{0\leq q\leq T}\|f(u_\mu(q,\omega,\Psi))\|_X.
		\]
		Dominated convergence in \eqref{eq:deterministic-fixed-interval}
		therefore yields
		\(\mathscr F(t_2,r_2)\to\mathscr F(t_1,r_1)\) in \(X\).
		Together with \eqref{eq:deterministic-integral-bound} and
		Lemma~\ref{lem:profile-space-estimates}\textup{(iii)}, this gives
		\(\mathscr F(t,\cdot)\in\mathcal Y\).
		
		For each integer \(N\geq0\), the joint continuity of
		\(\mathscr F\) implies uniform continuity on
		\([0,T]\times[0,N]\). Hence, as \(t_2\to t_1\) in \([0,T]\),
		we have	\(
		\sup_{0\leq r\leq N}
		\|\mathscr F(t_2,r)-\mathscr F(t_1,r)\|_X\to0
		.
		\)
		By \eqref{eq:deterministic-integral-bound}, 
		\[
		2^{-N}\sup_{0\leq r\leq N}
		\|\mathscr F(t_2,r)-\mathscr F(t_1,r)\|_X
		\leq
		2^{1-N}T\sup_{0\leq q\leq T}
		\|f(u_\mu(q,\omega,\Psi))\|_X\]
		and 
		\[
		\sum_{N=0}^{\infty}2^{1-N}T
		\sup_{0\leq q\leq T}\|f(u_\mu(q,\omega,\Psi))\|_X
		=
		4T\sup_{0\leq q\leq T}
		\|f(u_\mu(q,\omega,\Psi))\|_X<\infty.
		\]
		Thus, by \eqref{eq:profile-norm} and dominated convergence for series,
		\begin{equation}\label{eq:deterministic-weighted-continuity}
			\begin{aligned}
				\|\mathscr F(t_2,\cdot)-\mathscr F(t_1,\cdot)\|_{\mathcal Y}
				&=
				\sum_{N=0}^{\infty}2^{-N}
				\sup_{0\leq r\leq N}
				\|\mathscr F(t_2,r)-\mathscr F(t_1,r)\|_X
				\to0.
			\end{aligned}
		\end{equation}
		Since \(T>0\) is arbitrary, the proof is complete.
	\end{proof}

	\begin{lemma}\label{prop:lifted-map-continuity}
		Let \(T>0\), and choose \(m\in\mathbb N\) with \(m\geq T\).
		Then, for every \(\omega\in\mathcal O\) and any two
		functions \(\Psi^1,\Psi^2\in\mathcal Y\), we have
		\begin{equation}\label{eq:lifted-map-lipschitz}
			\begin{aligned}
				\sup_{0\leq t\leq T}
				\|\Gamma_\mu(t,\omega,\Psi^1)
				-\Gamma_\mu(t,\omega,\Psi^2)\|_{\mathcal Y}
				\leq
				2^{m+1}\bigl(1+LTe^{LT}\bigr)
				\|\Psi^1-\Psi^2\|_{\mathcal Y}.
			\end{aligned}
		\end{equation}
		Moreover,  the mapping
		\(
		(t,\Psi)\mapsto\Gamma_\mu(t,\omega,\Psi)
		\)
		is continuous from \([0,\infty)\times\mathcal Y\)
		into \(\mathcal Y\).
	\end{lemma}
	
	\begin{proof}
		Fix \(T>0\), \(\omega\in\mathcal O\), and
		\(m\in\mathbb N\) with \(m\geq T\).
		By \eqref{eq:uPsi}, \eqref{eq:P-exact-norm} and
		(\ref{ass:nonlinearity}),
		\[
		\begin{aligned}
			&\|u_\mu(q,\omega,\Psi^1)-u_\mu(q,\omega,\Psi^2)\|_X
			\leq\|\Psi^1(q)-\Psi^2(q)\|_X
			+L\int_0^q
			\|u_\mu(\xi,\omega,\Psi^1)-u_\mu(\xi,\omega,\Psi^2)\|_X\,d\xi,
		\end{aligned}
		\]
		where \(0\leq q\leq T\). Then,	Gronwall's inequality and \eqref{eq:profile-local-control} yield
		\begin{equation}\label{eq:generalized-solution-continuity}
			\begin{aligned}
				\sup_{0\leq q\leq T}
				\|u_\mu(q,\omega,\Psi^1)-u_\mu(q,\omega,\Psi^2)\|_X
				\leq e^{LT}
				\sup_{0\leq q\leq T}\|\Psi^1(q)-\Psi^2(q)\|_X
				\leq2^me^{LT}\|\Psi^1-\Psi^2\|_{\mathcal Y}.
			\end{aligned}
		\end{equation}
		
		By \eqref{eq:lifted-map}, \eqref{eq:translation-bound},
		\eqref{eq:generalized-solution-continuity},
		Lemma~\ref{lem:tr-continuou} and
		Lemma~\ref{lem:profile-space-estimates}\textup{(iii)},
		for \(0\leq t\leq T\),
		\[
		\begin{aligned}
			\|\Gamma_\mu(t,\omega,\Psi^1)
			-\Gamma_\mu(t,\omega,\Psi^2)\|_{\mathcal Y}
			&\leq
			2^{m}\|\Psi^1-\Psi^2\|_{\mathcal Y}
			+2L\int_0^t
			\|u_\mu(q,\omega,\Psi^1)-u_\mu(q,\omega,\Psi^2)\|_X\,dq
			\\&\leq
			2^{m}\bigl(1+LTe^{LT}\bigr)
			\|\Psi^1-\Psi^2\|_{\mathcal Y}.
		\end{aligned}
		\]
		This proves \eqref{eq:lifted-map-lipschitz}.
		
		Fix \(\Psi\in\mathcal Y\), and let \(t_2\to t_1\) in \([0,T]\).
		For each integer \(N\geq0\), uniform continuity of \(\Psi\)
		on \([0,N+m]\) gives
		\[
		\sup_{0\leq r\leq N}
		\|\Psi(t_2+r)-\Psi(t_1+r)\|_X\to0.
		\]
		Moreover,
		\[
		2^{-N}\sup_{0\leq r\leq N}
		\|\Psi(t_2+r)-\Psi(t_1+r)\|_X
		\leq2^{1-N}\sup_{0\leq r\leq N+m}\|\Psi(r)\|_X,
		\]
		and
		\[
		\begin{aligned}
			\sum_{N=0}^{\infty}2^{1-N}
			\sup_{0\leq r\leq N+m}\|\Psi(r)\|_X
			&=
			2^{m}\sum_{j=m}^{\infty}2^{-j}
			\sup_{0\leq r\leq j}\|\Psi(r)\|_X
			\leq2^{m}\|\Psi\|_{\mathcal Y}<\infty.
		\end{aligned}
		\]
		By \eqref{eq:profile-norm} and dominated convergence
		for series,
		\begin{equation}\label{eq:tau}
			\|\tau_{t_2}\Psi-\tau_{t_1}\Psi\|_{\mathcal Y}
			=
			\sum_{N=0}^{\infty}2^{-N}
			\sup_{0\leq r\leq N}
			\|\Psi(t_2+r)-\Psi(t_1+r)\|_X
			\to0.
		\end{equation}
		By \eqref{eq:lifted-map}, \eqref{eq:tau},
		\eqref{eq:deterministic-weighted-continuity},
		\eqref{eq:bounded estimate} and
		Lemma~\ref{lem:two-parameter-convolution},
		\begin{equation}\label{eq:con-gam}
			\begin{aligned}
				\|\Gamma_\mu(t_2,\omega,\Psi)
				-\Gamma_\mu(t_1,\omega,\Psi)\|_{\mathcal Y}
				&\leq
				\|\tau_{t_2}\Psi-\tau_{t_1}\Psi\|_{\mathcal Y}
				+\|\mathscr F(t_2,\cdot)-\mathscr F(t_1,\cdot)\|_{\mathcal Y}\\
				&\qquad+
				2\sup_{r\geq0}
				\|Z_\mu(t_2,r,\omega)-Z_\mu(t_1,r,\omega)\|_X
				\to0.
			\end{aligned}
		\end{equation}
		
		Finally, let \((t_2,\Psi^2)\to(t_1,\Psi^1)\).
		Choose \(T>t_1\) and \(m\in\mathbb N\) with \(m\geq T\).
		For \(t_2\in[0,T]\), \eqref{eq:lifted-map-lipschitz}
		and \eqref{eq:con-gam} give
		\[
		\begin{aligned}
			&\|\Gamma_\mu(t_2,\omega,\Psi^2)
			-\Gamma_\mu(t_1,\omega,\Psi^1)\|_{\mathcal Y}\\
			&\quad\leq
			2^{m+1}\bigl(1+LTe^{LT}\bigr)
			\|\Psi^2-\Psi^1\|_{\mathcal Y}
			+\|\Gamma_\mu(t_2,\omega,\Psi^1)
			-\Gamma_\mu(t_1,\omega,\Psi^1)\|_{\mathcal Y}
			\to0.
		\end{aligned}
		\]
	\end{proof}
	
	We next establish the joint measurability.  We use the path space
	\((\mathcal O,\mathcal F,\mathbb P)\) introduced above and the
	topology of uniform convergence on compact intervals on \(\mathcal O\).
	In particular, \(\mathcal B(\mathcal O)\subseteq\mathcal F\). Moreover, the notation \(\otimes\) 
	denotes the product \(\sigma\)-algebra.
	
	\begin{lemma}\label{lem:joint-measurability-convolution}
		For every   \(r\geq0\), the mapping
		\(	(t,\omega)\to Z_\mu (t,r,\omega)\)
		is measurable from
		\((\mathbb R_+\times\mathcal O,
		\mathcal B(\mathbb R_+)\otimes\mathcal F)\)
		into \((X,\mathcal B(X))\). 
	\end{lemma}
	
	\begin{proof}
		Fix \(r\geq0\). Since \(\mathcal O\) is a subspace of the separable
		metric space \(C_0(\mathbb R;X)\),
		\begin{equation}\label{eq:path-product-borel}
			\mathcal B(\mathbb R_+\times\mathcal O)
			=\mathcal B(\mathbb R_+)\otimes\mathcal B(\mathcal O)
			\subseteq\mathcal B(\mathbb R_+)\otimes\mathcal F.
		\end{equation}
		
		Let \((a_k,\omega_k)\to(a,\omega)\) in
		\(\mathbb R_+\times\mathcal O\), and choose \(M>0\)
		such that \(a,a_k\in[0,M]\) for every \(k\).
		Then
		\begin{equation}\label{eq:path-evaluation-continuity}
			\begin{aligned}
				\|\omega_k(a_k)-\omega(a)\|_X
				&\leq
				\sup_{0\leq q\leq M}
				\|\omega_k(q)-\omega(q)\|_X
				+\|\omega(a_k)-\omega(a)\|_X
				\to0.
			\end{aligned}
		\end{equation}
		Thus \((a,\omega)\mapsto\omega(a)\) is jointly continuous.
		
		For \(n\in\mathbb N\) and \(t\geq0\), let
		\(
		q_j:=\frac{jt}{2^n}\)
		with \(j=0,\ldots,2^n
		\)
		and define
		\begin{equation}\label{eq:young-riemann-sums}
			S_n(t,r,\omega)
			:=
			\sum_{j=0}^{2^n-1}
			P_\mu(t+r-q_j)G
			[\omega(q_{j+1})-\omega(q_j)],
			\qquad t\geq0.
		\end{equation}
		By \eqref{eq:PGH} and \eqref{eq:path-evaluation-continuity},
		each finite sum \(S_n(\cdot,r,\cdot)\) is jointly continuous,
		and therefore jointly measurable by
		\eqref{eq:path-product-borel}.
		By \eqref{eq:young-riemann-sums},
		Lemma~\ref{lem:young-integration} and
		\eqref{eq:two-parameter-convolution},
		\[
		Z_\mu(t,r,\omega)
		=\lim_{n\to\infty}S_n(t,r,\omega)
		\quad\text{in }X,
		\qquad t\geq0,\quad\omega\in\mathcal O.
		\]
		Hence \(Z_\mu(\cdot,r,\cdot)\) is jointly measurable.
	\end{proof}
	
	\begin{lemma}\label{prop:generalized-solution-measurability}
		The generalized solution mapping
		\((t,\omega,\Psi)\mapsto u_\mu(t,\omega,\Psi)\)
		is measurable from \(\mathbb R_+\times\mathcal O\times\mathcal Y\), equipped with
		\(\mathcal B(\mathbb R_+)\otimes\mathcal F\otimes\mathcal B(\mathcal Y)\),
		into \((X,\mathcal B(X))\).
	\end{lemma}

	\begin{proof}
		Set \(u_\mu^{(0)}(t,\omega,\Psi):=\Psi(t)+Z_\mu(t,0,\omega)\) and,
		for \(n\geq0\), define
		\begin{equation}\label{eq:profile-Picard-iteration}
			u_\mu^{(n+1)}(t,\omega,\Psi)
			:=u_\mu^{(0)}(t,\omega,\Psi)
			+\int_0^t P_\mu(t-q)f(u_\mu^{(n)}(q,\omega,\Psi))\,dq.
		\end{equation}
		
		For each fixed \((\omega,\Psi)\), take \(\kappa>L\) and use
		\(\chi,\mathcal X_\kappa,\Phi\) from the proof of
		Theorem~\ref{thm:generalized-profile-solution}.
		By \eqref{eq:xi-omega}, \eqref{eq:auxiliary-weight-definition}
		and \eqref{eq:fixed-point-operator}, we have
		\[
		u_\mu^{(0)}(\cdot,\omega,\Psi)=\phi^*\in\mathcal X_\kappa,
		\qquad
		u_\mu^{(n+1)}(\cdot,\omega,\Psi)
		=\Phi\bigl(u_\mu^{(n)}(\cdot,\omega,\Psi)\bigr).
		\]
		Since \(\Phi(\mathcal X_\kappa)\subseteq\mathcal X_\kappa\)
		by \eqref{eq:self-map-weighted-bound}, induction gives
		\begin{equation}\label{eq:picard-path-continuity}
			u_\mu^{(n)}(\cdot,\omega,\Psi)\in\mathcal X_\kappa
			\subset C([0,\infty);X),
			\qquad n\geq0.
		\end{equation}

		We next prove joint measurability of the iterates.
		For \((t_2,\Psi^2)\to(t_1,\Psi^1)\), choose \(T>t_1\)
		and \(m\in\mathbb N\) with \(m\geq T\).
		By \eqref{eq:profile-local-control}, for \(t_2\in[0,T]\),
		\[
		\|\Psi^2(t_2)-\Psi^1(t_1)\|_X
		\leq
		2^m\|\Psi^2-\Psi^1\|_{\mathcal Y}
		+\|\Psi^1(t_2)-\Psi^1(t_1)\|_X
		\to 0.
		\]
		Thus the mapping
		\begin{equation}\label{eq:joint-evaluation-map}
			(t,\Psi)\mapsto\Psi(t)
		\end{equation}
		is continuous and jointly measurable, since
		\(\mathcal Y\) is separable by
		Lemma~\ref{lem:profile-space-estimates}. Together with Lemma~\ref{lem:joint-measurability-convolution},
		this proves joint measurability of \(u_\mu^{(0)}\).
		
		Moreover, as \((a_2,y_2)\to(a_1,y_1)\), Lemma~\ref{prop:P-properties} and
		\eqref{eq:P-exact-norm} give
		\begin{equation}\label{eq:operator-action-continuity}
			\begin{aligned}
				\|P_\mu(a_2)y_2-P_\mu(a_1)y_1\|_X
				&\leq
				\|y_2-y_1\|_X
				+\|[P_\mu(a_2)-P_\mu(a_1)]y_1\|_X
				\to 0.
			\end{aligned}
		\end{equation}
		Assume that \(u_\mu^{(n)}\) is jointly measurable.
		By \eqref{eq:operator-action-continuity} and
		(\ref{ass:nonlinearity}), the mappings
		\[
		R_{n,k}(t,\omega,\Psi)
		:=
		\frac{t}{2^k}\sum_{j=0}^{2^k-1}
		P_\mu\left(t-\frac{jt}{2^k}\right)
		f\left(
		u_\mu^{(n)}\left(\frac{jt}{2^k},\omega,\Psi\right)
		\right)
		\]
		are jointly measurable.
		By \eqref{eq:picard-path-continuity} and
		\eqref{eq:operator-action-continuity},
		\begin{equation}\label{eq:picard-integral-approximation}
			R_{n,k}(t,\omega,\Psi)
			\to 
			\int_0^tP_\mu(t-q)
			f\bigl(u_\mu^{(n)}(q,\omega,\Psi)\bigr)\,dq
			\quad\text{in }X, \quad \text{as}\quad 
			k\to\infty.
		\end{equation}
		Equation	\eqref{eq:picard-integral-approximation} holds also at \(t=0\).
		Thus \eqref{eq:profile-Picard-iteration} implies joint
		measurability of \(u_\mu^{(n+1)}\).
		
		Finally, for each fixed \((\omega,\Psi)\),
		\(u_\mu(\cdot,\omega,\Psi)\) is the fixed point of \(\Phi\) in
		\(\mathcal X_\kappa\). By \eqref{eq:global-weighted-space}
		and \eqref{eq:global-weighted-contraction}, for every \(t\geq0\),
		\[
		\|u_\mu^{(n)}(t,\omega,\Psi)-u_\mu(t,\omega,\Psi)\|_X
		\leq \chi(t)\left(\frac{L}{\kappa}\right)^n
		\|u_\mu^{(0)}(\cdot,\omega,\Psi)-u_\mu(\cdot,\omega,\Psi)\|_\kappa
		\to0 .
		\]
		Thus \((t,\omega,\Psi)\mapsto u_\mu(t,\omega,\Psi)\) is jointly
		measurable as a pointwise limit of jointly measurable maps.
	\end{proof}

	\begin{lemma}\label{prop:lifted-map-measurability}
		The  map 
		\(
		\Gamma_\mu:\mathbb R_+\times\mathcal O\times\mathcal Y
		\to \mathcal Y
		\)
		is measurable with respect to
		\(\mathcal B(\mathbb R_+)\otimes\mathcal F
		\otimes\mathcal B(\mathcal Y)\)
		and \(\mathcal B(\mathcal Y)\).
	\end{lemma}
	
	\begin{proof}
		Fix \(r\geq0\). For \(k\in\mathbb N\), define
		\begin{equation}\label{eq:gamma-riemann-sums}
			\begin{aligned}
				R_k(t,r,\omega,\Psi)
				&:=
				\frac{t}{2^k}\sum_{j=0}^{2^k-1}
				P_\mu\left(t+r-\frac{jt}{2^k}\right)
				f\left(
				u_\mu\left(\frac{jt}{2^k},\omega,\Psi\right)
				\right),
				\qquad t\geq0.
			\end{aligned}
		\end{equation}
		By Lemma~\ref{prop:generalized-solution-measurability},
		\eqref{eq:operator-action-continuity} and
		(\ref{ass:nonlinearity}),
		\(
		(t,\omega,\Psi)\mapsto R_k(t,r,\omega,\Psi)
		\)
		is jointly measurable.
		
		For each fixed \((t,\omega,\Psi)\),
		Theorem~\ref{thm:generalized-profile-solution},
		(\ref{ass:nonlinearity}) and
		\eqref{eq:operator-action-continuity} imply that the mapping
		\begin{equation}\label{eq:gamma-integrand}
			q\mapsto P_\mu(t+r-q)f\bigl(u_\mu(q,\omega,\Psi)\bigr),
			\qquad q\in[0,t],
		\end{equation}
		is continuous into \(X\).
		Hence \eqref{eq:gamma-riemann-sums} and
		\eqref{eq:gamma-integrand} give
		\begin{equation}\label{eq:gamma-integral-approximation}
			R_k(t,r,\omega,\Psi)
			\to
			\int_0^t P_\mu(t+r-q)f\bigl(u_\mu(q,\omega,\Psi)\bigr)\,dq
			\quad\text{in }X,
			\qquad\text{as }k\to\infty.
		\end{equation}
		Therefore, the integral in
		\eqref{eq:gamma-integral-approximation}
		is jointly measurable. 
		Together with \eqref{eq:joint-evaluation-map},
		Lemma~\ref{lem:joint-measurability-convolution}, and
		\eqref{eq:lifted-map}, this proves that
		\begin{equation}\label{eq:gam-r-mea}
			(t,\omega,\Psi)\mapsto\Gamma_\mu(t,\omega,\Psi)(r)
		\end{equation}
		is jointly measurable for every fixed \(r\geq0\).
		
		For
		\(\Xi\in\mathcal Y\),
		continuity in \(r\) and \eqref{eq:profile-norm} gives
		\begin{equation}\label{eq:gam}
			\begin{aligned}
				&\|\Gamma_\mu(t,\omega,\Psi)-\Xi\|_{\mathcal Y}
				=
				\sum_{N=0}^{\infty}2^{-N}
				\sup_{r\in\mathbb Q\cap[0,N]}
				\|\Gamma_\mu(t,\omega,\Psi)(r)-\Xi(r)\|_X.
			\end{aligned}
		\end{equation}
		By \eqref{eq:gam-r-mea}, the right-hand side of \eqref{eq:gam} is measurable
		as a limit of sums of countable suprema.
		Consequently, for every \(\varepsilon>0\),
		\[
		\begin{aligned}
			&\Gamma_\mu^{-1}\bigl(B_{\mathcal Y}(\Xi,\varepsilon)\bigr)
			=
			\bigl\{(t,\omega,\Psi):
			\|\Gamma_\mu(t,\omega,\Psi)-\Xi\|_{\mathcal Y}
			<\varepsilon\bigr\}
		\end{aligned}
		\]
		is measurable.
		By Lemma~\ref{lem:profile-space-estimates},
		\(\mathcal Y\) is separable and has a countable base of
		open balls. Hence \(\Gamma_\mu\) is measurable as a
		\(\mathcal Y\)-valued mapping.
	\end{proof}

	\begin{theorem}\label{thm:lifted-random-dynamical-system}
		Suppose that (\ref{ass:kernel}),
		(\ref{ass:kernel-holder},) (\ref{ass:noise-coefficient})
		and (\ref{ass:nonlinearity}) hold and fix \(\mu\geq0\).
		Then the mapping
		\[
		\Gamma_\mu:\mathbb R_+\times\mathcal O\times\mathcal Y
		\to \mathcal Y
		\]
		defined by \eqref{eq:lifted-map} is a continuous random
		dynamical system on \(\mathcal Y\) over \(\theta\).
		Moreover, for every fixed \(\omega\in\mathcal O\),
		\[
		(t,\Psi)\mapsto\Gamma_\mu(t,\omega,\Psi)
		\]
		is continuous from \(\mathbb R_+\times\mathcal Y\)
		into \(\mathcal Y\).
	\end{theorem}
	
	\begin{proof}
		By Lemma~\ref{lem:profile-space-estimates},
		\(\mathcal Y\) is a separable Banach space.
		For \(\omega\in\mathcal O\) and \(\Psi\in\mathcal Y\), by Theorem~\ref{thm:lifted-cocycle}, 
		\(
		\Gamma_\mu(0,\omega,\Psi)=\Psi
		\)
		and the cocycle property \eqref{eq:lifted-cocycle} holds.
		Joint measurability follows from
		Lemma~\ref{prop:lifted-map-measurability}, while
		Lemma~\ref{prop:lifted-map-continuity} gives joint
		continuity in \((t,\Psi)\) for every fixed \(\omega\).
		The conclusion follows from
		Definition~\ref{def:random-dynamical-system}.
	\end{proof}

	\section{Existence of a random pullback attractor}
	\label{sec:generalized-pullback-attractor}
	
	In this section, we establish the existence of a minimal random
	pullback attractor for the lifted cocycle $\Gamma_\mu$.
	Fix \(\mu\geq0\), assume that the Lipschitz constant satisfies
	\(L<\lambda_1\), and set \(\varrho:=\lambda_1-L>0\).
	We first establish the kernel estimates needed for the pullback analysis. 
	\begin{lemma}\label{lem:attractor-kernel-estimates}
		Suppose that (\ref{ass:kernel}) holds and \(L<\lambda_1\).
		Let \(p_{\mu,\nu}\) be as in Lemma~\ref{lem:pnu-properties},
		and let \(P_\mu(t)\) be defined by
		\eqref{eq:P-spectral-definition}.
		Then the following properties hold for every \(\mu\geq0\).
		
		\smallskip
		\noindent\textup{(i)}
		For every \(\nu>0\) and \(t>0\),
		\begin{equation}\label{eq:attractor-scalar-integrability}
			\int_0^\infty p_{\mu,\nu}(s)\,ds\leq\frac1\nu,
			\qquad
			0\leq p_{\mu,\nu}(t)\leq
			\min\left\{1,\frac1{\nu t}\right\}.
		\end{equation}
		In particular, as \(t\to\infty\), we have
		\(
		p_{\mu,\nu}(t)\rightarrow0\) and \(t\,p_{\mu,\nu}(t)\to 0\).
		
		\smallskip
		\noindent\textup{(ii)}
		Let \(\varrho=\lambda_1-L>0\). We obtain
		\[	p_{\mu,\varrho}(t)
		=p_{\mu,\lambda_1}(t)
		+
		L\int_0^t
		p_{\mu,\lambda_1}(t-q)
		p_{\mu,\varrho}(q)\,dq,
		\qquad t\geq0.\]
		
		\smallskip
		\noindent\textup{(iii)}
		The operator family \(P_\mu\) belongs to
		\(C^1((0,\infty);\mathcal L(X))\).	Define
		\[
		b^*(t):=\min\left\{
		\frac1{et},\frac4{e\lambda_1t^2}
		\right\}, \qquad d^*(t):=\min\left\{1,\frac1{\lambda_1t}\right\},
		\qquad t>0.
		\]
		Then, for every \(t>0\) and \(r\geq0\),
		\begin{equation}\label{eq:attractor-operator-derivative}
			\|P_\mu(t+r)\|_{\mathcal L(X)}
			\leq d^*(t),
			\qquad
			\|P_\mu'(t+r)\|_{\mathcal L(X)}
			\leq b^*(t).
		\end{equation}
		
		\smallskip
		\noindent\textup{(iv)}
		For every \(t>0\) and \(r,h\geq0\),
		
		\[	\|P_\mu(t+r+h)-P_\mu(t+r)\|_{\mathcal L(X)}
		\leq
		\|P_\mu(t+h)-P_\mu(t)\|_{\mathcal L(X)}
		\leq p_{\mu,\lambda_1}(t).\]
	\end{lemma}
	\begin{proof}
		\textup{(i)}
		By \eqref{eq:pnu-laplace}, Lemma~\ref{lem:pnu-properties} and the monotone convergence theorem, 
		\[
		\begin{aligned}
			\int_0^\infty p_{\mu,\nu}(s)\,ds
			=\lim_{\lambda\to0}
			\int_0^\infty e^{-\lambda s}p_{\mu,\nu}(s)\,ds
			=\lim_{\lambda\to0}
			\frac1{\psi(\lambda)+\nu}
			\leq\frac1\nu.
		\end{aligned}
		\]
		By Lemma~\ref{lem:pnu-properties},
		\[
		0\leq p_{\mu,\nu}(t)
		\leq
		\min\left\{1,\frac1t\int_0^t p_{\mu,\nu}(s)\,ds\right\}
		\leq
		\min\left\{1,\frac1{\nu t}\right\},
		\qquad t>0,
		\]
		which proves \eqref{eq:attractor-scalar-integrability}.
		Moreover,
		\[
		0\leq\frac t2p_{\mu,\nu}(t)
		\leq\int_{t/2}^t p_{\mu,\nu}(s)\,ds
		\to 0
		\qquad\text{as }t\to\infty.
		\]
		Hence, as \(t\to\infty\), we get
		\(
		t\,p_{\mu,\nu}(t)\to 0,\)
		and
		\(p_{\mu,\nu}(t)\to 0.
		\)

		\textup{(ii)}
		Applying \eqref{eq:pnu-resolvent-identity} with
		\((\nu_1,\nu_2)=(\varrho,\lambda_1)\), and using
		\(\lambda_1-\varrho=L\) and 
		the convolution theorem, one has
		\[
		p_{\mu,\varrho}(t)-p_{\mu,\lambda_1}(t)
		=
		L\int_0^t
		p_{\mu,\lambda_1}(t-q)p_{\mu,\varrho}(q)\,dq,
		\qquad t\geq0.
		\]
		
		\textup{(iii)}
		Since \(\varpi_{\mu,\nu}\) is a probability measure,
		differentiating \eqref{eq-probability measure} under
		the integral sign gives
		\[
		0\leq-p_{\mu,\nu}'(t)
		=
		\int_{0}^{\infty}ze^{-tz}\,\varpi_{\mu,\nu}(dz)
		\leq
		\sup_{z\geq0}ze^{-tz}
		\leq\frac1{et}.
		\]
		Moreover, we have \[-p_{\mu,\nu}'(t)
		=
		\int_0^\infty
		(ze^{-tz/2})e^{-tz/2}\,\varpi_{\mu,\nu}(dz)
		\leq\frac2{et}p_{\mu,\nu}(t/2),\quad t>0.\]				Hence, for \(\nu\geq\lambda_1\), \eqref{eq:attractor-scalar-integrability} yields
		\(	|p_{\mu,\nu}'(t)|
		\leq
		\min\left\{\frac1{et},\frac4{e\lambda_1t^2}\right\}
		=b^*(t).\)
		Similarly,
		\begin{equation}\label{eq:two-derivative}
			0\leq p_{\mu,\nu}''(t)
			=
			\int_0^\infty z^2e^{-tz}\,\varpi_{\mu,\nu}(dz)
			\leq
			\sup_{z\geq0}z^2e^{-tz}
			=
			\frac4{e^2t^2}.
		\end{equation}
		
		Fix \(t>0\). By 
		Parseval's identity, we have 
		\(
		\sum_{n=1}^{\infty}
		|p_{\mu,\lambda_n}'(t)|^2
		|v_n|^2
		\leq
		\bigl(b^*(t)\bigr)^2\|v\|_X^2\) for 
		\(v\in X\).
		Thus the series
		\(\sum_{n=1}^{\infty}p_{\mu,\lambda_n}'(t)v_n\varphi_n\)
		converges in \(X\) and defines a bounded linear operator.
		For \(0<|h|<t/2\) and each \(n\geq1\), Taylor's theorem and
		\eqref{eq:two-derivative} give
		\[
		\left|
		\frac{p_{\mu,\lambda_n}(t+h)-p_{\mu,\lambda_n}(t)}{h}
		-p_{\mu,\lambda_n}'(t)
		\right|
		\leq \frac{|h|}{2}
		\sup_{q\in[t/2,3t/2]}|p_{\mu,\lambda_n}''(q)|
		\leq \frac{8|h|}{e^2t^2}.
		\]
		Equation \eqref{eq:P-spectral-definition} and Parseval's identity yield
		\[
		\sup_{\|v\|_X\leq1}
		\left\|
		\frac{P_\mu(t+h)-P_\mu(t)}{h}v
		-\sum_{n=1}^{\infty}p_{\mu,\lambda_n}'(t)v_n\varphi_n
		\right\|_X
		\leq \frac{8|h|}{e^2t^2}
		\to0
		\qquad\text{as }h\to0.
		\]		
		Consequently,
		\begin{equation}\label{eq:P-derivative-spectral}
			P_\mu'(t)x
			=
			\sum_{n=1}^{\infty}
			p_{\mu,\lambda_n}'(t)v_n\varphi_n
			\quad\text{and}\quad 
			\|P_\mu'(t)\|_{\mathcal L(X)}
			=
			\sup_{n\geq1}|p_{\mu,\lambda_n}'(t)|
			\leq b^*(t).
		\end{equation}
		By \eqref{eq:P-derivative-spectral},
		\eqref{eq:two-derivative} and the mean value theorem,
		for \(s,t\geq\varepsilon>0\),
		\[
		\begin{aligned}
			\|P_\mu'(t)-P_\mu'(s)\|_{\mathcal L(X)}
			&=
			\sup_{n\geq1}
			|p_{\mu,\lambda_n}'(t)-p_{\mu,\lambda_n}'(s)|
			\leq
			\frac4{e^2\varepsilon^2}|t-s|.
		\end{aligned}
		\]
		Hence \(P_\mu'\) is locally Lipschitz continuous and
		\(
		P_\mu\in C^1((0,\infty);\mathcal L(X)).
		\)
		
		Finally, \eqref{eq:P-exact-norm},
		\eqref{eq:attractor-scalar-integrability},
		\eqref{eq:P-derivative-spectral} and the monotonicity
		of \(d^*\) and \(b^*\) yield
		\[
		\|P_\mu(t+r)\|_{\mathcal L(X)}
		\leq d^*(t+r)\leq d^*(t) \quad \text{and }\quad 
		\|P_\mu'(t+r)\|_{\mathcal L(X)}
		\leq b^*(t+r)\leq b^*(t),
		\]
		for \(t>0\) and \(r\geq0\).	This proves \eqref{eq:attractor-operator-derivative}.				
		
		\textup{(iv)}
		For \(t,r,h\geq0\), \eqref{eq-probability measure} gives
		\[
		\begin{aligned}
			0\leq
			p_{\mu,\nu}(t+r)-p_{\mu,\nu}(t+r+h)
			&\leq
			\int_0^\infty
			e^{-tz}(1-e^{-hz})\,\varpi_{\mu,\nu}(dz)
			=
			p_{\mu,\nu}(t)-p_{\mu,\nu}(t+h).
		\end{aligned}
		\]
		Taking the supremum over \(\nu=\lambda_n\) and using
		\(
		0\leq
		p_{\mu,\lambda_n}(t)-p_{\mu,\lambda_n}(t+h)
		\leq p_{\mu,\lambda_1}(t),
		\)
		we obtain
		\[
		\begin{aligned}
			\|P_\mu(t+r+h)-P_\mu(t+r)\|_{\mathcal L(X)}
			&\leq
			\|P_\mu(t+h)-P_\mu(t)\|_{\mathcal L(X)}
			\leq p_{\mu,\lambda_1}(t).
		\end{aligned}
		\]
	\end{proof}
	
	For nonempty subsets \( \mathcal V, \mathcal Q\subset\mathcal Y\), define
	\begin{equation*}
		\Gamma_\mu(t,\omega,\mathcal V)
		:=\{\Gamma_\mu(t,\omega,\Psi):\Psi\in \mathcal V\},
	\end{equation*}
	and  the Hausdorff semidistance by
	\begin{equation*}
		\operatorname{dist}_{\mathcal Y}(\mathcal V,\mathcal Q)
		:=
		\sup_{\Psi\in \mathcal V}\inf_{\widetilde\Psi\in \mathcal Q}
		\|\Psi-\widetilde\Psi\|_{\mathcal Y}.
	\end{equation*}
	Next, we recall the definitions used below, following  \cite[Sections~2 and~3]{Cra} and \cite[Appendix]{cao}.
	
	\begin{definition}\label{def:attractor-random-closed-set}
		A family
		\(\mathcal C=\{\mathcal C(\omega)\}_{\omega\in\mathcal O}\)
		of nonempty closed subsets of \(\mathcal Y\) is called
		a random closed set in \(\mathcal Y\) if, for every
		\(\Psi\in\mathcal Y\), the map
		\[
		\omega\mapsto
		\inf_{\widetilde\Psi\in\mathcal C(\omega)}
		\|\Psi-\widetilde\Psi\|_{\mathcal Y}
		\]
		is measurable.
		A random closed set \(\mathcal C(\omega)\) is called a random
		compact set if \(\mathcal C(\omega)\) is compact  for every \(\omega\in\mathcal O\).
	\end{definition}
	
	\begin{definition}\label{def:random-pullback-attractor}
		Let \(\mathfrak D\) be a nonempty collection of families
		\(D=\{D(\omega)\}_{\omega\in\mathcal O}\) of nonempty subsets of
		\(\mathcal Y\). A family
		\(\mathcal A_\mu=\{\mathcal A_\mu(\omega)\}_{\omega\in\mathcal O}\)
		is a \emph{random pullback attractor} for \(\Gamma_\mu\) relative to
		\(\mathfrak D\) if the following conditions hold.
		
		\smallskip
		\noindent\textup{(i)}
		\(\mathcal A_\mu\) is a random compact set in \(\mathcal Y\).
		
		\smallskip
		\noindent\textup{(ii)}
		For every \(t\geq0\),
		\[
		\Gamma_\mu(t,\omega,\mathcal A_\mu(\omega))
		=
		\mathcal A_\mu(\theta_t\omega)
		\qquad
		\mathbb P\text{-almost surely}.
		\]

		\smallskip
		\noindent\textup{(iii)}
		For every \(D\in\mathfrak D\),
		\[
		\lim_{t\to\infty}
		\operatorname{dist}_{\mathcal Y}
		\bigl(
		\Gamma_\mu(t,\theta_{-t}\omega,D(\theta_{-t}\omega)),
		\mathcal A_\mu(\omega)
		\bigr)
		=0
		\qquad
		\mathbb P\text{-almost surely}.
		\]
		Such an attractor \(\mathcal A_\mu\) is called minimal if,
		for every random pullback attractor \(\mathcal K\)
		for \(\Gamma_\mu\) relative to the same family \(\mathfrak D\),
		\[
		\mathcal A_\mu(\omega)\subseteq\mathcal K(\omega)
		\qquad
		\mathbb P\text{-almost surely}.
		\]
	\end{definition}
	For a nonempty closed bounded random set \(B=\{B(\omega)\}_{\omega\in\mathcal O}\) in \(X\),
	define
	\(
	r_B(\omega)
	:=
	\sup_{x\in B(\omega)}\|x\|_X.
	\)
	Denote by \(\mathfrak D_{\mu,\varrho}\)
	the collection of all such random sets satisfying
	\begin{equation}\label{eq:admissible-random-initial-set}
		\lim_{T\to\infty}
		p_{\mu,\varrho}(T)\,r_B(\theta_{-T}\omega)=0
		\qquad\mathbb P\text{-almost surely}.
	\end{equation}
	For \(B\in\mathfrak D_{\mu,\varrho}\), define 
	\begin{equation}\label{eq:widetilde B}
		\widetilde B(\omega):=\{\Psi_{\mu,x}:x\in B(\omega)\},
		\qquad \Psi_{\mu,x}(r)=P_\mu(r)x.
	\end{equation}
	Since \(B(\omega)\) is complete, 	\eqref{eq:Psi-x-bound},
	\eqref{eq:profile-norm}, \eqref{eq:Psi-xy} and the equation
	\(\Psi_{\mu,x}-\Psi_{\mu,y}=\Psi_{\mu,x-y}\) with \(x,y\in X\)  imply that
	\(\widetilde B(\omega)\) is closed in \(\mathcal Y\). Moreover,
	\(\sup_{\Psi\in\widetilde B(\omega)}\|\Psi\|_{\mathcal Y}
	\leq 2r_B(\omega)<\infty.\)
	Consequently, 
	a measurable dense selection representation of \(B\)
	shows that \(r_B\) is measurable and that
	\(\widetilde B\) is a random closed set.
	We consider the family
	\begin{equation}\label{eq:attractor-initial-family}
		\mathfrak D_{\widetilde B}
		:=
		\{\widetilde B:B\in\mathfrak D_{\mu,\varrho}\}.
	\end{equation}
	
	\subsection{Bounds and convergence of the stochastic convolution}
	\label{subsec:attractor-noise}
	Fix \(\zeta\in(H,1)\) and retain the exponent
	\(\gamma\in(\max\{1/2,1-\rho\},H)\) used in the Young integral.
	To control the convolution over the infinite past, set
	\begin{equation}\label{eq:attractor-polynomial-path-set}
		\mathcal O_0
		:=
		\left\{
		\omega\in\mathcal O:
		\sup_{s\in\mathbb R}
		\frac{\|\omega(s)\|_X}{(1+|s|)^\zeta}<\infty
		\right\}.
	\end{equation}
	
	\begin{lemma}\label{prop:fbm-path-moments}
		For the two-sided \(X\)-valued \(Q\)-fractional Brownian motion \(B^H\),
		\begin{equation}\label{eq:fbm-local-holder-moment}
			\mathbb E\bigl([B^H]_{\gamma;[-M,M]}\bigr)<\infty,
			\qquad M>0,
		\end{equation}
		and
		\begin{equation}\label{eq:fbm-weighted-path-moment}
			\mathbb E\sup_{t\in\mathbb R}
			\frac{\|B^H(t)\|_X}{(1+|t|)^\zeta}<\infty.
		\end{equation}
	\end{lemma}
	\begin{proof}
		By the Gaussianity of \(B^H\), \eqref{eq:covariance structure} and
		\eqref{eq:Q-fbm-representation}, for every finite \(q\geq1\),
		\begin{equation*}
			\mathbb E\|B^H(t)-B^H(s)\|_X^q
			\leq C_q|t-s|^{qH},
			\qquad s,t\in\mathbb R.
		\end{equation*}
		Choose \(q>\max\{1,(H-\gamma)^{-1}\}\), so that
		\(\gamma<H-1/q\).
		Applying the quantitative Kolmogorov estimate
		\cite[Theorem~1.1]{Kra}  on \([-M,M]\),
		and then using H\"older's inequality, gives
		\[
		\mathbb E\!\left([B^H]_{\gamma;[-M,M]}\right)
		\leq
		\left(
		\mathbb E\!\left([B^H]_{\gamma;[-M,M]}^q\right)
		\right)^{1/q}
		<\infty,
		\qquad M>0.
		\]
		This proves \eqref{eq:fbm-local-holder-moment}.
		
		By \(B^H(0)=0\) and
		\eqref{eq:fbm-local-holder-moment}, we obtain
		\(\mathbb E\sup_{0\leq t\leq1}\|B^H(t)\|_X
		\leq
		\mathbb E\!\left([B^H]_{\gamma;[-1,1]}\right)
		<\infty.\)
		Self-similarity of fBm, equality in law of
		\(t\to B^H(-t)\) and
		\(t\to B^H(t)\) yield
		\begin{equation*}
			\begin{aligned}
				\mathbb E\sup_{|t|\leq T}\|B^H(t)\|_X
				&\leq
				\mathbb E\sup_{0\leq t\leq T}\|B^H(t)\|_X
				+
				\mathbb E\sup_{0\leq t\leq T}\|B^H(-t)\|_X
				=
				2T^H\mathbb E\sup_{0\leq t\leq1}\|B^H(t)\|_X
				\leq CT^H
				.
			\end{aligned}
		\end{equation*}
		For \(j\geq0\) and \(2^j\leq|t|\leq2^{j+1}\),
		\(
		\frac{\|B^H(t)\|_X}{(1+|t|)^\zeta}
		\leq
		2^{-j\zeta}
		\sup_{|s|\leq2^{j+1}}\|B^H(s)\|_X.
		\)
		Therefore,
		\begin{equation*}
			\sup_{t\in\mathbb R}
			\frac{\|B^H(t)\|_X}{(1+|t|)^\zeta}
			\leq
			\sup_{|t|\leq1}\|B^H(t)\|_X
			+
			\sum_{j=0}^{\infty}
			2^{-j\zeta}
			\sup_{|t|\leq2^{j+1}}\|B^H(t)\|_X.
		\end{equation*}
		Since \(\zeta>H\),	from
		Tonelli's theorem, we derive
		\[
		\begin{aligned}
			\mathbb E\sup_{t\in\mathbb R}
			\frac{\|B^H(t)\|_X}{(1+|t|)^\zeta}
			&\leq
			\mathbb E\sup_{|t|\leq1}\|B^H(t)\|_X
			+
			\sum_{j=0}^{\infty}
			2^{-j\zeta}
			\mathbb E\sup_{|t|\leq2^{j+1}}\|B^H(t)\|_X\\
			&\leq
			C\left(
			1+\sum_{j=0}^{\infty}2^{-j(\zeta-H)}
			\right)
			<\infty.
		\end{aligned}
		\]  This proves \eqref{eq:fbm-weighted-path-moment}
	\end{proof}

	\begin{lemma}\label{lem:O_0}
		The set \(\mathcal O_0\) is measurable, has \(
		\mathbb P(\mathcal O_0)=1\)
		and satisfies
		\(\theta_s\mathcal O_0=\mathcal O_0\)
		for  \(s\in\mathbb R\).	
	\end{lemma}
	\begin{proof}
		By continuity of the paths, the supremum in
		\eqref{eq:attractor-polynomial-path-set} can be taken over
		rational times. It is therefore measurable, which proves
		measurability of \(\mathcal O_0\).
		By \eqref{eq:fbm-weighted-path-moment}, the supremum in
		\eqref{eq:attractor-polynomial-path-set} is finite
		almost surely. Together with \(\mathbb P(\mathcal O)=1\),
		this gives \(\mathbb P(\mathcal O_0)=1\).
		
		For \(\omega\in\mathcal O_0\) and  \(t\in\mathbb R\), we have
		\begin{equation}\label{eq:bound-attractor-path-set}
			\|\omega(t)\|_X\leq C(1+|t|)^\zeta
		\end{equation}
		from \eqref{eq:attractor-polynomial-path-set}. Using
		\(1+|s+r|\leq(1+|s|)(1+|r|)\) and  \eqref{eq:Wiener shift}, for every fixed \(s\in\mathbb R\), we get
		\[
		\begin{aligned}
			\|(\theta_s\omega)(r)\|_X
			&\leq C\bigl((1+|s+r|)^\zeta+(1+|s|)^\zeta\bigr)
			\leq
			2C(1+|s|)^\zeta(1+|r|)^\zeta.
		\end{aligned}
		\]
		Thus, \(
		\sup_{r\in\mathbb R}
		\frac{\|(\theta_s\omega)(r)\|_X}{(1+|r|)^\zeta}
		<\infty.
		\)
		Together with the shift invariance of \(\mathcal O\),
		this gives \(\theta_s\mathcal O_0\subseteq\mathcal O_0\).
		Applying the same argument to \(\theta_{-s}\), and using
		\(\theta_s\theta_{-s}=I\), proves
		\(\theta_s\mathcal O_0=\mathcal O_0\).
	\end{proof}
	
	We use \(d^*\) and \(b^*\) from
	Lemma~\ref{lem:attractor-kernel-estimates}\textup{(iii)}, with \(d^*(0):=1\). Then, to obtain a uniform pathwise bound for the stochastic
	convolution over past time intervals, we introduce the
	following random variable.
	\begin{lemma}\label{lem:attractor-noise-majorant}
		For \(\omega\in\mathcal O_0\), define
		\begin{equation}\label{eq:attractor-noise-majorant}
			\mathcal N(\omega)
			:=\|G\|_{\mathcal L(X)}
			\left[
			\sup_{t\geq0}d^*(t)\|\omega(-t)\|_X
			+\int_0^\infty b^*(t)\|\omega(-t)\|_X\,dt
			\right],
		\end{equation}
		and  \(\mathcal N(\omega):=0\) for
		\(\omega\in\mathcal O\setminus\mathcal O_0\).
		Then \(\mathcal N\) is a finite nonnegative measurable
		random variable on \((\mathcal O,\mathcal F,\mathbb P)\).
		For every \(\omega\in\mathcal O_0\), the map
		\(s\to\mathcal N(\theta_s\omega)\) is continuous
		on \(\mathbb R\). Moreover, for every \(N\in\mathbb N\),
		\begin{equation}\label{eq:attractor-noise-N}
			\mathcal N_N(\omega)
			:=\sup_{|s|\leq N}\mathcal N(\theta_s\omega)
		\end{equation}
		is measurable and satisfies
		\(\mathbb E\mathcal N_N<\infty\).
	\end{lemma}
	
	\begin{proof}
		Fix \(\omega\in\mathcal O_0\) and \(N\in\mathbb N\).
		For  \( |s|\leq N\) and \(t\geq1\),	we have
		\(
		1+|s-t|\leq(N+2)t\)
		and
		\(	1+|s|\leq(N+1)t.
		\)
		By	\eqref{eq:local Holder} and
		\eqref{eq:attractor-polynomial-path-set},	we obtain
		\begin{equation}\label{eq:shifted-noise-increment-bound}
			\sup_{|s|\leq N}\|\omega(s-t)-\omega(s)\|_X
			\leq
			\begin{cases}
				[\omega]_{\gamma;[-N-1,N+1]}t^\gamma,
				&0<t\leq1,\\[1ex]
				\displaystyle
				Ct^\zeta
				\sup_{v\in\mathbb R}
				\frac{\|\omega(v)\|_X}{(1+|v|)^\zeta},
				&t\geq1.
			\end{cases}
		\end{equation}
		By \eqref{eq:shifted-noise-increment-bound} and the
		definitions of \(d^*,b^*\) in
		Lemma~\ref{lem:attractor-kernel-estimates}\textup{(iii)},
		for \(|s|\leq N\),
		\begin{equation*}
			\sup_{t\geq0}d^*(t)
			\|\omega(s-t)-\omega(s)\|_X
			\leq
			[\omega]_{\gamma;[-N-1,N+1]}
			+C\sup_{v\in\mathbb R}
			\frac{\|\omega(v)\|_X}{(1+|v|)^\zeta},
		\end{equation*}
		and
		\begin{equation}\label{eq:integral}
			\begin{aligned}
				&\int_0^\infty b^*(t)
				\|\omega(s-t)-\omega(s)\|_X\,dt\\
				&\quad\leq
				\frac{[\omega]_{\gamma;[-N-1,N+1]}}{e}
				\int_0^1t^{\gamma-1}\,dt
				+
				\frac{4C}{e\lambda_1}
				\left(
				\sup_{v\in\mathbb R}
				\frac{\|\omega(v)\|_X}{(1+|v|)^\zeta}
				\right)
				\int_1^\infty t^{\zeta-2}\,dt
				<\infty.
			\end{aligned}
		\end{equation}
		Combining the two bounds, \eqref{eq:Wiener shift} and
		\eqref{eq:attractor-noise-majorant} yields
		\begin{equation*}
			\mathcal N_N(\omega)
			\leq C\|G\|_{\mathcal L(X)}
			\left[
			[\omega]_{\gamma;[-N-1,N+1]}
			+\sup_{v\in\mathbb R}
			\frac{\|\omega(v)\|_X}{(1+|v|)^\zeta}
			\right]<\infty.
		\end{equation*}
		In particular, \(0\leq\mathcal N(\omega)<\infty\).
		On \(\mathcal O\setminus\mathcal O_0\),
		\(\mathcal N(\omega)=0\).
		
		For every \(\omega\in\mathcal O\), continuity of the
		paths and the monotone convergence theorem give
		\begin{equation*}
			\sup_{t\geq0}d^*(t)\|\omega(-t)\|_X
			=
			\sup_{t\in\mathbb Q\cap[0,\infty)}
			d^*(t)\|\omega(-t)\|_X
		\end{equation*}
		and \begin{equation*}
			\int_0^\infty b^*(t)\|\omega(-t)\|_X\,dt
			=
			\lim_{m\to\infty}
			\int_{1/m}^{m}b^*(t)\|\omega(-t)\|_X\,dt.
		\end{equation*}
		Each fixed-time evaluation and each truncated integral
		are continuous in the compact-open topology.
		Thus both expressions are measurable.
		Lemma~\ref{lem:O_0} and
		\eqref{eq:attractor-noise-majorant} imply that the zero
		extension \(\mathcal N\) is measurable.
		
		For the fixed \(\omega\in\mathcal O_0\), let
		\(s_k\to s\) in \([-N,N]\).
		For every \(R\geq1\), uniform continuity of \(\omega\)
		on \([-N-R,N]\) gives
		\[\sup_{0\leq t\leq R}
		\|\omega(s_k-t)-\omega(s_k)
		-\omega(s-t)+\omega(s)\|_X
		\to0.\]
		Moreover, for \(t\geq R\), \eqref{eq:shifted-noise-increment-bound} and
		\(d^*(t)\leq1/(\lambda_1t)\) yield
		\begin{equation}\label{eq:noise-shift-sup-tail}
			\sup_{\substack{|a|\leq N}}
			d^*(t)\|\omega(a-t)-\omega(a)\|_X
			\leq C_NR^{\zeta-1}
			\sup_{v\in\mathbb R}
			\frac{\|\omega(v)\|_X}{(1+|v|)^\zeta}
			\to0
			\quad\text{as }R\to\infty.
		\end{equation}
		Splitting the supremum at \(R\), we obtain
		\begin{equation}\label{eq:convergence of the supremum}
			\begin{aligned}
				&\sup_{t\geq0}d^*(t)
				\|\omega(s_k-t)-\omega(s_k)\|_X
				\to
				\sup_{t\geq0}d^*(t)
				\|\omega(s-t)-\omega(s)\|_X.
			\end{aligned}
		\end{equation}
		By \eqref{eq:shifted-noise-increment-bound}, \eqref{eq:integral} 
		and the dominated convergence,
		\begin{equation}\label{eq:convergence of the integral}
			\begin{aligned}
				&\int_0^\infty b^*(t)
				\|\omega(s_k-t)-\omega(s_k)\|_X\,dt
				\to
				\int_0^\infty b^*(t)
				\|\omega(s-t)-\omega(s)\|_X\,dt.
			\end{aligned}
		\end{equation}
		By \eqref{eq:convergence of the supremum},
		\eqref{eq:convergence of the integral},
		\eqref{eq:Wiener shift} and
		\eqref{eq:attractor-noise-majorant}, we get
		\(
		\mathcal N(\theta_{s_k}\omega)
		\to
		\mathcal N(\theta_s\omega).
		\)
		Since \(N\) is arbitrary, the mapping
		\(s\mapsto\mathcal N(\theta_s\omega)\)
		is continuous on \(\mathbb R\).
		
		By Lemma~\ref{lem:O_0}, the complement
		\(\mathcal O\setminus\mathcal O_0\) is also shift
		invariant, and
		\(\mathcal N(\theta_s\omega)=0\) there for every
		\(s\in\mathbb R\).
		Consequently,
		\(
		\mathcal N_N(\omega)
		=
		\sup_{s\in\mathbb Q\cap[-N,N]}
		\mathcal N(\theta_s\omega)\) for 
		\(\omega\in\mathcal O.
		\)
		Measurability of \(\mathcal N\) and \(\theta_s\)
		proves measurability of \(\mathcal N_N\).
		Finally, 
		\eqref{eq:fbm-local-holder-moment},
		\eqref{eq:fbm-weighted-path-moment} and
		\(\mathbb P(\mathcal O_0)=1\) give
		\[
		\begin{aligned}
			\mathbb E\mathcal N_N
			&\leq C_N\|G\|_{\mathcal L(X)}
			\left[
			\mathbb E\!\left(
			[B^H]_{\gamma;[-N-1,N+1]}
			\right)
			+
			\mathbb E\sup_{v\in\mathbb R}
			\frac{\|B^H(v)\|_X}{(1+|v|)^\zeta}
			\right]
			<\infty.
		\end{aligned}
		\]  The proof is complete.
	\end{proof}

	For the uniform estimates and convergence of the stochastic
	convolution below, let
	\[
	C_b([0,\infty);X)
	:=
	\left\{
	\Psi\in C([0,\infty);X):
	\|\Psi\|_\infty:=\sup_{r\geq0}\|\Psi(r)\|_X<\infty
	\right\}.
	\]
	Moreover, Lemma~\ref{lem:profile-space-estimates}\textup{(iii)}
	gives
	\begin{equation}\label{eq:C_b}
		C_b([0,\infty);X)\subset\mathcal Y,
		\qquad
		\|\Psi\|_{\mathcal Y}\leq2\|\Psi\|_\infty,
		\quad\Psi\in C_b([0,\infty);X).
	\end{equation}
	
	\begin{lemma}\label{lem:attractor-past-noise-convolution}
		Fix \(\mu\geq0\). For \(T,r\geq0\) and
		\(\omega\in\mathcal O_0\), define
		\begin{equation}\label{eq:attractor-past-noise}
			S_{\mu,T}(\omega)(r)
			:=\int_{-T}^0P_\mu(r-s)G\,d\omega(s)
			=Z_\mu(T,r,\theta_{-T}\omega).
		\end{equation}
		Then the following statements hold.
		
		\smallskip
		\noindent\textup{(i)}
		For every \(\omega\in\mathcal O_0\) and \(T\geq0\),
		\(S_{\mu,T}(\omega)\in C_b([0,\infty);X)\), and
		\begin{equation}\label{eq:attractor-past-noise-bound}
			\sup_{T\geq0}\|S_{\mu,T}(\omega)\|_\infty
			\leq\mathcal N(\omega).
		\end{equation}
		
		\smallskip
		\noindent\textup{(ii)}
		For every \(\omega\in\mathcal O_0\), the map
		\(T\mapsto S_{\mu,T}(\omega)\) is continuous from \([0,\infty)\)
		into \(C_b([0,\infty);X)\).
		
		\smallskip
		\noindent\textup{(iii)}
		There exists a measurable mapping
		\(W_\mu:\mathcal O\to\mathcal Y\), with
		\(W_\mu(\omega)=0\) for
		\(\omega\in\mathcal O\setminus\mathcal O_0\), such that,
		for every \(\omega\in\mathcal O_0\),
		\(W_\mu(\omega)\in C_b([0,\infty);X)\) and
		\begin{equation}\label{eq:attractor-past-noise-limit}
			\lim_{T\to\infty}S_{\mu,T}(\omega)=W_\mu(\omega)
			\quad\text{in }C_b([0,\infty);X).
		\end{equation}
		By \eqref{eq:C_b}, this convergence also holds in \(\mathcal Y\).
	\end{lemma}
	
	\begin{proof}
		Fix \(\omega\in\mathcal O_0\).
		By Lemma~\ref{lem:O_0},
		\eqref{eq:young-shift} and \eqref{eq:two-parameter-convolution},
		\[
		\int_{-T}^0P_\mu(r-s)G\,d\omega(s)
		=\int_0^TP_\mu(r+T-q)G\,d(\theta_{-T}\omega)(q)
		=Z_\mu(T,r,\theta_{-T}\omega).
		\]
		Thus \eqref{eq:attractor-past-noise} is well defined.
		
		\textup{(i)}
		Let \(T>0\) and \(0<\varepsilon<\min\{T,1\}\).
		By Lemma~\ref{lem:attractor-kernel-estimates}\textup{(iii)} gives
		\begin{equation}\label{eq:S_T}
			\begin{aligned}
				S_{\mu,T}(\omega)(r)
				=&
				\int_{-T}^{-\varepsilon}P_\mu(r-s)G\,d\omega(s)
				+\int_{-\varepsilon}^0P_\mu(r-s)G\,d\omega(s)\\
				=&-P_\mu(r+T)G\omega(-T)
				+P_\mu(r+\varepsilon)G\omega(-\varepsilon)\\
				&+\int_\varepsilon^TP_\mu'(r+t)G\omega(-t)\,dt
				+\int_{-\varepsilon}^0P_\mu(r-s)G\,d\omega(s).
			\end{aligned}
		\end{equation}
		By \eqref{eq:attractor-operator-derivative} and
		\eqref{eq:integral} at \(s=0\), we get
		\(
		\int_0^T\|P_\mu'(r+t)G\omega(-t)\|_X\,dt
		<\infty.
		\)
		Since \(t\mapsto P_\mu'(r+t)G\omega(-t)\) is continuous
		for \(t>0\), the Bochner integral
		\(\int_0^TP_\mu'(r+t)G\omega(-t)\,dt\) is well defined.
		By \eqref{eq:young-estimate}, \eqref{eq:PGH},
		\eqref{eq:P-exact-norm} and \eqref{eq:local Holder},
		\begin{equation}\label{eq:past-noise-young-endpoint}
			\sup_{r\geq0}
			\left\|\int_{-\varepsilon}^0P_\mu(r-s)G\,d\omega(s)\right\|_X
			\leq C\|G\|_{\mathcal L(X,X_\delta)}
			[\omega]_{\gamma;[-1,0]}
			\bigl(\varepsilon^\gamma+\varepsilon^{\rho+\gamma}\bigr)
			\to0.
		\end{equation}
		Moreover, as \(\varepsilon\to0\), \eqref{eq:P-exact-norm},
		\eqref{eq:attractor-operator-derivative},
		\eqref{eq:local Holder} and \(b^*(t)\leq1/(et)\) yield
		\begin{equation}\label{eq:past-noise-bochner-endpoint}
			\begin{aligned}
				&\sup_{r\geq0}\|P_\mu(r+\varepsilon)G\omega(-\varepsilon)\|_X
				+\sup_{r\geq0}
				\left\|\int_0^\varepsilon P_\mu'(r+t)G\omega(-t)\,dt\right\|_X\\
				&\quad\leq
				\|G\|_{\mathcal L(X)}[\omega]_{\gamma;[-1,0]}
				\left(\varepsilon^\gamma+
				\frac1e\int_0^\varepsilon t^{\gamma-1}\,dt\right)
				=
				\|G\|_{\mathcal L(X)}[\omega]_{\gamma;[-1,0]}
				\left(1+\frac1{e\gamma}\right)\varepsilon^\gamma
				\to0 .
			\end{aligned}
		\end{equation}
		Letting \(\varepsilon\to0\) in \eqref{eq:S_T}, using
		\eqref{eq:past-noise-young-endpoint} and
		\eqref{eq:past-noise-bochner-endpoint}, gives
		\begin{equation}\label{eq:attractor-noise-ibp}
			S_{\mu,T}(\omega)(r)
			=-P_\mu(r+T)G\omega(-T)
			+\int_0^TP_\mu'(r+t)G\omega(-t)\,dt,
			\qquad r\geq0.
		\end{equation}
		By \eqref{eq:attractor-noise-ibp},
		\eqref{eq:attractor-operator-derivative} and
		\eqref{eq:attractor-noise-majorant},
		\[
		\|S_{\mu,T}(\omega)\|_\infty
		\leq\|G\|_{\mathcal L(X)}
		\left[d^*(T)\|\omega(-T)\|_X
		+\int_0^Tb^*(t)\|\omega(-t)\|_X\,dt\right]
		\leq\mathcal N(\omega).
		\]
		Together with \(S_{\mu,0}(\omega)=0\), this proves
		\eqref{eq:attractor-past-noise-bound}.
		Continuity in \(r\) follows from \eqref{eq:attractor-past-noise}
		and Lemma~\ref{lem:two-parameter-convolution}.
		Hence \(S_{\mu,T}(\omega)\in C_b([0,\infty);X)\).
		
		\textup{(ii)}
		Let \(0\leq T_1<T_2\leq M\). By
		\eqref{eq:attractor-past-noise},
		\[
		S_{\mu,T_2}(\omega)(r)-S_{\mu,T_1}(\omega)(r)
		=\int_{-T_2}^{-T_1}P_\mu(r-s)G\,d\omega(s).
		\]
		Thus \eqref{eq:young-estimate}, \eqref{eq:PGH} and
		\eqref{eq:P-exact-norm} give
		\begin{equation*}
			\begin{aligned}
				&\|S_{\mu,T_2}(\omega)-S_{\mu,T_1}(\omega)\|_\infty
				\leq C\|G\|_{\mathcal L(X,X_\delta)}
				[\omega]_{\gamma;[-M,0]}
				\left(|T_2-T_1|^\gamma
				+|T_2-T_1|^{\rho+\gamma}\right).
			\end{aligned}
		\end{equation*}
		Consequently, we have 
		\(
		\|S_{\mu,T_2}(\omega)-S_{\mu,T_1}(\omega)\|_\infty
		\to 0
		\) as \(T_2\to T_1\).
		
		\textup{(iii)}
		By \eqref{eq:attractor-noise-ibp},
		\eqref{eq:attractor-operator-derivative} and
		\eqref{eq:bound-attractor-path-set}, for \(T\geq M\geq1\),
		\begin{equation*}
			\begin{aligned}
				&\|S_{\mu,T}(\omega)-S_{\mu,M}(\omega)\|_\infty\\
				&\quad\leq\|G\|_{\mathcal L(X)}
				\left[d^*(M)\|\omega(-M)\|_X
				+d^*(T)\|\omega(-T)\|_X
				+\int_M^Tb^*(t)\|\omega(-t)\|_X\,dt\right]\\
				&\quad\leq C(\omega)
				\left(M^{\zeta-1}+T^{\zeta-1}
				+\int_M^Tt^{\zeta-2}\,dt\right)
				\leq C(\omega)M^{\zeta-1}
				\to0
				\quad\text{as }M\to\infty.
			\end{aligned}
		\end{equation*}
		Completeness of \(C_b([0,\infty);X)\) gives
		\(W_\mu(\omega)\in C_b([0,\infty);X)\) such that
		\(
		\lim_{T\to\infty}
		\|S_{\mu,T}(\omega)-W_\mu(\omega)\|_\infty=0.
		\)
		This proves
		\eqref{eq:attractor-past-noise-limit}. By \eqref{eq:C_b},
		\begin{equation}\label{eq:past-noise-weighted-convergence}
			\|S_{\mu,T}(\omega)-W_\mu(\omega)\|_{\mathcal Y}
			\leq2\|S_{\mu,T}(\omega)-W_\mu(\omega)\|_\infty
			\to0, \quad\text{as }T\to\infty.
		\end{equation}
		
		It remains to prove measurability of \(W_\mu\).	For fixed \(n\in\mathbb N\) and \(r\geq0\),
		Lemma~\ref{lem:joint-measurability-convolution} and
		measurability of \(\theta_{-n}\) show that
		\(\omega\mapsto Z_\mu(n,r,\theta_{-n}\omega)\) is measurable.
		By Lemma~\ref{lem:two-parameter-convolution} and \eqref{eq:C_b},
		\(Z_\mu(n,\cdot,\theta_{-n}\omega)\in\mathcal Y\)
		for every \(\omega\in\mathcal O\).
		For \(\Xi\in\mathcal Y\), continuity in \(r\) and
		\eqref{eq:profile-norm} give
		\[
		\begin{aligned}
			&\|Z_\mu(n,\cdot,\theta_{-n}\omega)-\Xi\|_{\mathcal Y}
			=
			\sum_{N=0}^{\infty}2^{-N}
			\sup_{r\in\mathbb Q\cap[0,N]}
			\|Z_\mu(n,r,\theta_{-n}\omega)-\Xi(r)\|_X.
		\end{aligned}
		\]
		The right-hand side is measurable. By separability of
		\(\mathcal Y\) in Lemma~\ref{lem:profile-space-estimates},
		\(\omega\mapsto Z_\mu(n,\cdot,\theta_{-n}\omega)\)
		is therefore \(\mathcal Y\)-valued measurable.
		Set \(W_\mu(\omega):=0\) on
		\(\mathcal O\setminus\mathcal O_0\).
		By Lemma~\ref{lem:O_0}, \eqref{eq:attractor-past-noise} and
		\eqref{eq:past-noise-weighted-convergence},
		\[
		W_\mu(\omega)=\lim_{n\to\infty}
		\begin{cases}
			Z_\mu(n,\cdot,\theta_{-n}\omega),
			&\omega\in\mathcal O_0,\\
			0,&\omega\in\mathcal O\setminus\mathcal O_0,
		\end{cases}
		\quad\text{in }\mathcal Y.
		\]
		Hence \(W_\mu:\mathcal O\to\mathcal Y\) is measurable.
	\end{proof}

	\subsection{Pullback estimates and absorption}
	\label{subsec:attractor-pullback-bounds}			
	Throughout this subsection, fix \(\mu\geq0\) and set
	\(\varrho:=\lambda_1-L>0\), as in
	Lemma~\ref{lem:attractor-kernel-estimates}\textup{(ii)}.
	Define
	\[
	\mathcal H_\mu
	:=
	\left\{
	p_{\mu,\varrho},
	p_{\mu,\lambda_1},
	p_{\mu,\lambda_1}*p_{\mu,\varrho}
	\right\}.
	\]
	All three are nonnegative.	By \eqref{eq:attractor-scalar-integrability} and
	Fubini--Tonelli theorem,
	\begin{equation}\label{eq:attractor-kernels-L1}
		\begin{aligned}
			\|p_{\mu,\varrho}\|_{L^1(\mathbb R_+)}
			\leq\frac1{\varrho},
			\quad
			\|p_{\mu,\lambda_1}\|_{L^1(\mathbb R_+)}
			\leq\frac1{\lambda_1}, \quad \text{and}\quad 
			\|p_{\mu,\lambda_1}*p_{\mu,\varrho}\|_{L^1(\mathbb R_+)}
			\leq\frac1{\lambda_1\varrho}.
		\end{aligned}
	\end{equation}
	Thus every \(h\in\mathcal H_\mu\) is nonnegative and
	belongs to \(L^1(\mathbb R_+)\). 	Define
	\begin{equation}\label{eq:attractor-good-set-preliminary}
		\mathcal O_{1,\mu}
		:=
		\left\{
		\omega\in\mathcal O_0:
		\int_0^\infty h(t)\mathcal N_N(\theta_{-t}\omega)\,dt<\infty
		\text{ for every }h\in\mathcal H_\mu,\ N\in\mathbb N
		\right\}.
	\end{equation}				
	
	\begin{lemma}\label{lem:attractor-invariant-good-set}
		The set \(\mathcal O_{1,\mu}\) is measurable, has probability one,
		and satisfies \(\theta_s\mathcal O_{1,\mu}=\mathcal O_{1,\mu}\)
		for every \(s\in\mathbb R\).
		For \(\omega\in\mathcal O_{1,\mu}\), define
		\begin{align}
			R_{0,\mu}(\omega)
			&:=\frac{\|f(0)\|_X}{\varrho}+\mathcal N(\omega)
			+L\int_0^\infty p_{\mu,\varrho}(t)\mathcal N(\theta_{-t}\omega)\,dt,
			\label{eq:attractor-state-radius}\\
			R_{1,\mu}(\omega)
			&:=\frac{\|f(0)\|_X}{\lambda_1}+\mathcal N(\omega)
			+L\int_0^\infty p_{\mu,\lambda_1}(t)R_{0,\mu}(\theta_{-t}\omega)\,dt,
			\label{eq:attractor-profile-radius}
		\end{align}
		and set both functions equal to zero outside \(\mathcal O_{1,\mu}\).
		Then \(R_{0,\mu}\) and \(R_{1,\mu}\) are finite measurable random variables,
		and \(0\leq R_{1,\mu}\leq R_{0,\mu}\).
		For \(i=0,1\) and every \(\omega\in\mathcal O\),
		the map \(s\mapsto R_{i,\mu}(\theta_s\omega)\) is continuous on \(\mathbb R\).
		Moreover, for each \(N\in\mathbb N\),
		\begin{equation}\label{eq:attractor-radius-local-moments}
			\mathbb E\sup_{|s|\leq N}R_{i,\mu}(\theta_s\omega)<\infty,
			\qquad i=0,1.
		\end{equation}
	\end{lemma}

	\begin{proof}
		For \(h\in\mathcal H_\mu\) and \(N\in\mathbb N\),
		Lemma~\ref{lem:attractor-noise-majorant} and joint
		measurability of \(\theta\) imply that
		\(
		\omega\to
		\int_0^\infty
		h(t)\mathcal N_N(\theta_{-t}\omega)\,dt
		\)
		is measurable.
		By Fubini--Tonelli theorem, \eqref{eq:attractor-kernels-L1} and
		\(\mathbb E\mathcal N_N<\infty\) in Lemma \ref{lem:attractor-noise-majorant}, one has
		\begin{equation}\label{eq:radius-weighted-noise-moment}
			\begin{aligned}
				\mathbb E\int_0^\infty
				h(t)\mathcal N_N(\theta_{-t}\omega)\,dt
				=
				\int_0^\infty
				h(t)\mathbb E\bigl[
				\mathcal N_N(\theta_{-t}\omega)
				\bigr]\,dt
				=
				\|h\|_{L^1(\mathbb R_+)}
				\mathbb E\mathcal N_N
				<\infty.
			\end{aligned}
		\end{equation}
		Since \(\mathcal H_\mu\) is finite,
		\eqref{eq:attractor-good-set-preliminary},
		\eqref{eq:radius-weighted-noise-moment} and
		Lemma~\ref{lem:O_0} give
		\(	\mathbb P(\mathcal O_{1,\mu})=1.\)

		For \(s\in\mathbb R\),
		choose \(m\in\mathbb N\) with \(m\geq|s|\).
		By \eqref{eq:attractor-noise-majorant}, we obtain
		\[
		\begin{aligned}
			\mathcal N_N(\theta_{-t}\theta_s\omega)
			=
			\sup_{|r|\leq N}
			\mathcal N(\theta_{r+s-t}\omega)
			\leq
			\sup_{|q|\leq N+m}
			\mathcal N(\theta_{q-t}\omega)
			=
			\mathcal N_{N+m}(\theta_{-t}\omega).
		\end{aligned}
		\]
		Thus \eqref{eq:attractor-good-set-preliminary} and
		Lemma~\ref{lem:O_0} imply
		\(\theta_s\mathcal O_{1,\mu}\subseteq\mathcal O_{1,\mu}\).
		Applying the same inclusion with \(-s\) yields
		\(
		\theta_s\mathcal O_{1,\mu}=\mathcal O_{1,\mu}.
		\)
		
		Since \(\mathcal N\leq\mathcal N_1\) from \eqref{eq:attractor-noise-N},
		\eqref{eq:attractor-good-set-preliminary},
		\eqref{eq:attractor-state-radius} and
		Lemma~\ref{lem:attractor-noise-majorant} give
		\(0\leq R_{0,\mu}(\omega)<\infty\).
		By Fubini--Tonelli theorem,
		\begin{equation}\label{eq:Fubini--Tonelli theorem}
			\begin{aligned}
				&\int_0^\infty\int_0^\infty
				p_{\mu,\lambda_1}(t)p_{\mu,\varrho}(r)
				\mathcal N(\theta_{-(t+r)}\omega)\,dr\,dt=
				\int_0^\infty
				(p_{\mu,\lambda_1}*p_{\mu,\varrho})(q)
				\mathcal N(\theta_{-q}\omega)\,dq.
			\end{aligned}
		\end{equation}
		Substituting \eqref{eq:attractor-state-radius} and
		\eqref{eq:Fubini--Tonelli theorem} into
		\eqref{eq:attractor-profile-radius},
		Lemma~\ref{lem:attractor-kernel-estimates}\textup{(ii)} and \eqref{eq:attractor-kernels-L1} yield
		\begin{equation}\label{eq:attractor-radius-expansion}
			\begin{aligned}
				R_{1,\mu}(\omega)
				={}&
				\frac{\|f(0)\|_X}{\lambda_1}
				+\frac{L\|f(0)\|_X}{\varrho}
				\|p_{\mu,\lambda_1}\|_{L^1(\mathbb R_+)}
				+\mathcal N(\omega)\\
				&+L\int_0^\infty
				\bigl[
				p_{\mu,\lambda_1}
				+L(p_{\mu,\lambda_1}*p_{\mu,\varrho})
				\bigr](t)
				\mathcal N(\theta_{-t}\omega)\,dt\\
				={}&
				R_{0,\mu}(\omega)
				+\frac{L\|f(0)\|_X}{\varrho}
				\left(
				\|p_{\mu,\lambda_1}\|_{L^1(\mathbb R_+)}
				-\frac1{\lambda_1}
				\right)
				\leq R_{0,\mu}(\omega).
			\end{aligned}
		\end{equation}
		Hence \(0\leq R_{1,\mu}(\omega)\leq R_{0,\mu}(\omega)<\infty\).
		For \(i=0,1\),	measurability of \(R_{i,\mu}\) follows from joint measurability of
		\((t,\omega)\mapsto\mathcal N(\theta_{-t}\omega)\) and measurability
		of \(\mathcal O_{1,\mu}\). 
		Since \(\mathcal O_{1,\mu}\in\mathcal F\), both zero
		extensions are measurable on \(\mathcal O\).
		
		For \(|s|\leq N\),  by
		\eqref{eq:attractor-good-set-preliminary}, we have
		\begin{equation*}
			\begin{aligned}
				0\leq
				p_{\mu,\varrho}(t)\mathcal N(\theta_{s-t}\omega)
				\leq
				p_{\mu,\varrho}(t)\mathcal N_N(\theta_{-t}\omega)\quad \text{and}\quad 
				\int_0^\infty
				p_{\mu,\varrho}(t)\mathcal N_N(\theta_{-t}\omega)\,dt
				<\infty.
			\end{aligned}
		\end{equation*}
		Since \(N\) is arbitrary, 	 continuity of
		\(s\mapsto\mathcal N(\theta_s\omega)\) in Lemma~\ref{lem:attractor-noise-majorant}  and dominated convergence in
		\eqref{eq:attractor-state-radius} give continuity of
		\(s\mapsto R_{0,\mu}(\theta_s\omega)\).
		By invariance of \(\mathcal O_{1,\mu}\) and
		\eqref{eq:attractor-radius-expansion}, the same
		convergence holds for \(R_{1,\mu}\). On \(\mathcal O\setminus\mathcal O_{1,\mu}\), invariance
		of the complement gives
		\(
		R_{i,\mu}(\theta_s\omega)=0\) for
		\( i=0,1.
		\)
		\(\sup_{|s|\leq N}R_{i,\mu}(\theta_s\omega)\) can therefore
		be taken over rational \(s\) and are measurable. 
		Finally, \eqref{eq:attractor-state-radius},
		\eqref{eq:attractor-radius-expansion},
		\eqref{eq:radius-weighted-noise-moment}  and \(\mathbb E\mathcal N_N<\infty\) in Lemma \ref{lem:attractor-noise-majorant} yield
		\[
		\begin{aligned}
			\mathbb E\sup_{|s|\leq N}R_{i,\mu}(\theta_s\omega)
			\leq
			\mathbb E\sup_{|s|\leq N}R_{0,\mu}(\theta_s\omega)
			&\leq
			\frac{\|f(0)\|_X}{\varrho}
			+\mathbb E\mathcal N_N
			+L\mathbb E\int_0^\infty
			p_{\mu,\varrho}(t)\mathcal N_N(\theta_{-t}\omega)\,dt\\
			&=
			\frac{\|f(0)\|_X}{\varrho}
			+\left(
			1+L\|p_{\mu,\varrho}\|_{L^1(\mathbb R_+)}
			\right)\mathbb E\mathcal N_N
			<\infty, \quad i=0,1.
		\end{aligned}
		\]	This proves \eqref{eq:attractor-radius-local-moments}.
	\end{proof}

	For the rest of this section, we restrict the metric dynamical system to
	\(\mathcal O_{1,\mu}\) with its trace \(\sigma\)-algebra
	\(\mathcal F_{1,\mu}:=\{F\cap\mathcal O_{1,\mu}:F\in\mathcal F\}\)
	and \(\mathbb P\) denotes the restricted probability measure.
	By invariance, the restriction of \(\Gamma_\mu\) remains a continuous
	random dynamical system on \(\mathcal Y\).
	For \(s,t\in\mathbb R\) with \(t\geq s\),
	\eqref{eq:mild-stochastic} and \eqref{eq:young-shift} give
	\begin{equation}\label{eq:attractor-absolute-mild}
		\begin{aligned}
			u_\mu(t-s,\theta_s\omega,x)
			=P_\mu(t-s)x
			+\int_s^tP_\mu(t-q)f(u_\mu(q-s,\theta_s\omega,x))\,dq
			+\int_s^tP_\mu(t-q)G\,d\omega(q).
		\end{aligned}
	\end{equation}
	\begin{lemma}\label{lem:attractor-state-comparison}
		For \(s,t\in\mathbb R\) with \(t\geq s\),
		\(x,y\in X\), and \(\omega\in\mathcal O_{1,\mu}\), we have
		\begin{align}
			&\|u_\mu(t-s,\theta_s\omega,x)\|_X
			\leq
			p_{\mu,\varrho}(t-s)\|x\|_X+R_{0,\mu}(\theta_t\omega),
			\label{eq:attractor-state-comparison}\\
			&\bigl\|u_\mu(t-s,\theta_s\omega,x)
			-u_\mu(t-s,\theta_s\omega,y)\bigr\|_X
			\leq
			p_{\mu,\varrho}(t-s)\|x-y\|_X.
			\label{eq:attractor-state-difference}
		\end{align}
	\end{lemma}
	
	\begin{proof}
		Fix \(s,t\in\mathbb R\) with \(s\leq t\) and
		\(\omega\in\mathcal O_{1,\mu}\). 
		By \eqref{eq:young-shift}, \eqref{eq:attractor-past-noise},
		and \eqref{eq:attractor-past-noise-bound},
		\[
		\left\|\int_s^{s+r}P_\mu(s+r-a)G\,d\omega(a)\right\|_X
		=\|S_{\mu,r}(\theta_{s+r}\omega)(0)\|_X
		\leq\mathcal N(\theta_{s+r}\omega).
		\]
		For \(0\leq r\leq t-s\), set
		\(
		h(r):=\|u_\mu(r,\theta_s\omega,x)\|_X
		\)
		and
		\[
		g(r):=
		p_{\mu,\lambda_1}(r)\|x\|_X
		+\|f(0)\|_X\int_0^r p_{\mu,\lambda_1}(q)\,dq
		+\mathcal N(\theta_{s+r}\omega).
		\]
		Hence \eqref{eq:attractor-absolute-mild},
		\eqref{eq:P-exact-norm}, and \eqref{eq:f-lipschitz} give
		\begin{equation}\label{eq:absolute-mild}
			h(r)\leq g(r)
			+L\int_0^r p_{\mu,\lambda_1}(r-q)h(q)\,dq.
		\end{equation}
		Equation	\eqref{eq:absolute-mild}		and 
		Lemma~\ref{lem:attractor-kernel-estimates}\textup{(ii)} yields
		\begin{equation}\label{eq:h-pvarrho}
			\begin{aligned}
				h
				&=h+L\bigl[
				p_{\mu,\varrho}-p_{\mu,\lambda_1}
				-L(p_{\mu,\varrho}*p_{\mu,\lambda_1})
				\bigr]*h\\&=(h-Lp_{\mu,\lambda_1}*h)
				+Lp_\varrho*(h-Lp_{\mu,\lambda_1}*h)
				\leq g+Lp_{\mu,\varrho}*g.
			\end{aligned}
		\end{equation}
		Substituting the definition of \(g\) into \eqref{eq:h-pvarrho}
		and using    Lemma~\ref{lem:attractor-kernel-estimates}\textup{(ii)} gives
		\[
		\begin{aligned}
			h(r)
			&\leq g(r)+L(p_{\mu,\varrho}*g)(r)\\
			&=p_{\mu,\varrho}(r)\|x\|_X
			+\|f(0)\|_X\int_0^r p_{\mu,\varrho}(a)\,da
			+\mathcal N(\theta_{s+r}\omega)
			+L\int_0^r p_{\mu,\varrho}(a)
			\mathcal N(\theta_{s+r-a}\omega)\,da. 
		\end{aligned}
		\]
		Taking \(r=t-s\) and using \eqref{eq:attractor-kernels-L1},
		\eqref{eq:attractor-state-radius}, and
		\(\theta_{t-a}\omega=\theta_{-a}(\theta_t\omega)\), we obtain
		\[
		\|u_\mu(t-s,\theta_s\omega,x)\|_X
		\leq p_{\mu,\varrho}(t-s)\|x\|_X
		+R_{0,\mu}(\theta_t\omega).
		\]
		This proves \eqref{eq:attractor-state-comparison}.
		
		Finally, subtracting \eqref{eq:attractor-absolute-mild} for
		\(x\) and \(y\) cancels the stochastic integral.
		By \eqref{eq:P-exact-norm} and \eqref{eq:f-lipschitz},
		\[
		\begin{aligned}
			&\|u_\mu(r,\theta_s\omega,x)-u_\mu(r,\theta_s\omega,y)\|_X\\
			&\quad
			-L\int_0^r p_{\mu,\lambda_1}(r-q)
			\|u_\mu(q,\theta_s\omega,x)-u_\mu(q,\theta_s\omega,y)\|_X\,dq
			\leq p_{\mu,\lambda_1}(r)\|x-y\|_X.
		\end{aligned}
		\]
		Similarly, 
		\eqref{eq:absolute-mild}, \eqref{eq:h-pvarrho} gives
		\[
		\begin{aligned}
			&\|u_\mu(r,\theta_s\omega,x)-u_\mu(r,\theta_s\omega,y)\|_X
			\leq
			\bigl[p_{\mu,\lambda_1}
			+L(p_{\mu,\varrho}*p_{\mu,\lambda_1})\bigr](r)\|x-y\|_X
			=p_{\mu,\varrho}(r)\|x-y\|_X.
		\end{aligned}
		\]
		Taking \(r=t-s\) proves \eqref{eq:attractor-state-difference}.
	\end{proof}

	Define the random closed ball
	\begin{equation}\label{eq:attractor-absorbing-ball}
		B_{R,\mu}(\omega):=
		\left\{
		\Upsilon\in\mathcal Y:
		\|\Upsilon\|_{\mathcal Y}\leq2R_{1,\mu}(\omega)+1
		\right\}.
	\end{equation}
	The following Lemma   shows that
	\(B_{R,\mu}\) pullback absorbs every member of
	\(\mathfrak D_{\widetilde B}\).
	\begin{lemma}\label{prop:attractor-profile-bounds}
		For every \(T\geq0\), \(\omega\in\mathcal O_{1,\mu}\),
		and \(x,y\in X\),  we have
		\begin{align}
			&\|\Gamma_\mu(T,\theta_{-T}\omega,\Psi_{\mu,x})\|_{\mathcal Y}
			\leq
			2p_{\mu,\varrho}(T)\|x\|_X+2R_{1,\mu}(\omega),
			\label{eq:attractor-profile-bound}\\
			&\|\Gamma_\mu(T,\theta_{-T}\omega,\Psi_{\mu,x})
			-\Gamma_\mu(T,\theta_{-T}\omega,\Psi_{\mu,y})\|_{\mathcal Y}
			\leq
			2p_{\mu,\varrho}(T)\|x-y\|_X.
			\label{eq:attractor-profile-difference}
		\end{align}
		Moreover, for every \(B\in\mathfrak D_{\mu,\varrho}\)
		and  \(\omega\in\mathcal O_{1,\mu}\),
		there exists \(T_B(\omega)\geq0\) such that
		\[
		\Gamma_\mu(T,\theta_{-T}\omega,
		\widetilde B(\theta_{-T}\omega))
		\subseteq B_{R,\mu}(\omega),
		\qquad T\geq T_B(\omega).
		\]
	\end{lemma}
	
	\begin{proof}
		Fix \(T\geq0\) and \(\omega\in\mathcal O_{1,\mu}\).
		\eqref{eq:Psi-x}, \eqref{eq:lifted-map} and \eqref{eq:attractor-past-noise},
		\eqref{eq:lifted-map} give,
		for \(r\geq0\),
		\begin{equation}\label{eq:attractor-profile-decomposition}
			\begin{aligned}
				\Gamma_\mu(T,\theta_{-T}\omega,\Psi_{\mu,x})(r)
				&=P_\mu(T+r)x+S_{\mu,T}(\omega)(r)
				+\int_0^T P_\mu(r+a)
				f\bigl(u_\mu(T-a,\theta_{-T}\omega,x)\bigr)\,da.
			\end{aligned}
		\end{equation}
		By \eqref{eq:P-exact-norm} and the complete monotonicity in
		Lemma~\ref{lem:pnu-properties},
		\[
		\|P_\mu(r+a)\|_{\mathcal L(X)}
		=p_{\mu,\lambda_1}(r+a)
		\leq p_{\mu,\lambda_1}(a),
		\qquad r,a\geq0.
		\]
		Subtracting \eqref{eq:attractor-profile-decomposition} for \(x\)
		and \(y\), and using \eqref{eq:f-lipschitz}, \eqref{eq:attractor-state-difference}
		and Lemma~\ref{lem:attractor-kernel-estimates}\textup{(ii)}, we obtain
		\begin{equation}\label{eq:attractor-bound1}
			\begin{aligned}
				&\|\Gamma_\mu(T,\theta_{-T}\omega,\Psi_{\mu,x})
				-\Gamma_\mu(T,\theta_{-T}\omega,\Psi_{\mu,y})\|_\infty\\
				&\quad\leq
				\left[p_{\mu,\lambda_1}(T)
				+L\int_0^T p_{\mu,\lambda_1}(a)
				p_{\mu,\varrho}(T-a)\,da\right]\|x-y\|_X
				=p_{\mu,\varrho}(T)\|x-y\|_X,
			\end{aligned}
		\end{equation}
		Taking \(x=0\) in \eqref{eq:attractor-profile-decomposition}
		and applying \eqref{eq:attractor-state-comparison} with
		\(s=-T\) and \(t=-a\),
		we  find
		\begin{equation}\label{eq:attractor-bound2}
			\begin{aligned}
				\|\Gamma_\mu(T,\theta_{-T}\omega,\Psi_{\mu,0})\|_\infty
				\leq \mathcal N(\omega)
				+\int_0^T p_{\mu,\lambda_1}(a)
				\bigl[\|f(0)\|_X+LR_{0,\mu}(\theta_{-a}\omega)\bigr]\,da
				\leq R_{1,\mu}(\omega),
			\end{aligned}
		\end{equation}
		by \eqref{eq:f-lipschitz}, \eqref{eq:attractor-past-noise-bound}, \eqref{eq:attractor-kernels-L1} and
		\eqref{eq:attractor-profile-radius}.
		Since \(\Gamma_\mu(T,\theta_{-T}\omega,\Psi_{\mu,x})\in\mathcal Y\)
		is continuous, \eqref{eq:attractor-bound1}, \eqref{eq:attractor-bound2}, 
		and \eqref{eq:C_b} give
		\[
		\begin{aligned}
			\|\Gamma_\mu(T,\theta_{-T}\omega,\Psi_{\mu,x})\|_{\mathcal Y}
			&\leq
			2\|\Gamma_\mu(T,\theta_{-T}\omega,\Psi_{\mu,x})
			-\Gamma_\mu(T,\theta_{-T}\omega,\Psi_{\mu,0})\|_\infty\\
			&\quad +2\|\Gamma_\mu(T,\theta_{-T}\omega,\Psi_{\mu,0})\|_\infty
			\leq 2p_{\mu,\varrho}(T)\|x\|_X+2R_{1,\mu}(\omega),
		\end{aligned}
		\]
		and
		\[
		\begin{aligned}
			&\|\Gamma_\mu(T,\theta_{-T}\omega,\Psi_{\mu,x})
			-\Gamma_\mu(T,\theta_{-T}\omega,\Psi_{\mu,y})\|_{\mathcal Y}\\
			&\quad\leq
			2\|\Gamma_\mu(T,\theta_{-T}\omega,\Psi_{\mu,x})
			-\Gamma_\mu(T,\theta_{-T}\omega,\Psi_{\mu,y})\|_\infty
			\leq 2p_{\mu,\varrho}(T)\|x-y\|_X.
		\end{aligned}
		\]
		This proves \eqref{eq:attractor-profile-bound} and
		\eqref{eq:attractor-profile-difference}.
		
		For every \(B\in\mathfrak D_{\mu,\varrho}\) and
		\(\omega\in\mathcal O_{1,\mu}\) satisfying
		\eqref{eq:admissible-random-initial-set},
		there exists \(T_B(\omega)\geq0\) such that
		\(
		2p_{\mu,\varrho}(T)r_B(\theta_{-T}\omega)\leq1\)
		for \(T\geq T_B(\omega).
		\)
		For these \(T\), \eqref{eq:attractor-profile-bound} yields
		\[
		\begin{aligned}
			\sup_{x\in B(\theta_{-T}\omega)}
			\|\Gamma_\mu(T,\theta_{-T}\omega,\Psi_{\mu,x})\|_{\mathcal Y}
			&\leq 2p_{\mu,\varrho}(T)r_B(\theta_{-T}\omega)
			+2R_{1,\mu}(\omega)
			\leq 2R_{1,\mu}(\omega)+1.
		\end{aligned}
		\]
		Hence, by \eqref{eq:widetilde B} and
		\eqref{eq:attractor-absorbing-ball},
		\[
		\Gamma_\mu\bigl(T,\theta_{-T}\omega,
		\widetilde B(\theta_{-T}\omega)\bigr)
		\subseteq B_{R,\mu}(\omega),
		\qquad T\geq T_B(\omega).
		\]
		
		Finally, for every fixed \(\Phi^*\in\mathcal Y\),
		\begin{equation*}
			\operatorname{dist}_{\mathcal Y}(\{\Phi^*\},B_{R,\mu}(\omega))
			:=
			\inf_{\Upsilon\in B_{R,\mu} (\omega)}
			\|\Phi^*-\Upsilon\|_{\mathcal Y}
			=
			\max\left\{
			\|\Phi^*\|_{\mathcal Y}-(2R_{1,\mu}(\omega)+1),\,0
			\right\}.
		\end{equation*}
		The right-hand side  in above formula is measurable  by
		Lemma~\ref{lem:attractor-invariant-good-set}.
		Thus \(B_{R,\mu}\) is a random closed ball.
	\end{proof}

	\subsection{Pullback asymptotic compactness }
	\label{subsec:attractor-compactness}
	\begin{lemma}\label{lem:attractor-history-compactness}
		Let \(\mathcal G:[0,\infty)\to[0,\infty)\) be measurable and satisfy
		\begin{equation}\label{eq:attractor-history-envelope}
			K_{\mu,\mathcal G}:=\int_0^\infty p_{\mu,\lambda_1}(t) \mathcal G(t)\,dt<\infty.
		\end{equation}
		For every strongly measurable \(F:[0,\infty)\to X\) with
		\(\|F(t)\|_X\leq \mathcal G(t)\), define
		\[
		\mathscr{F}_{\mu,F}(r)=\int_0^\infty P_\mu(r+t)F(t)\,dt.
		\qquad r\geq0.
		\]
		Then, \(\mathscr{F}_{\mu,F}\in C_b([0,\infty);X)\subset\mathcal Y\) and the
		family of all such \(\mathscr{F}_{\mu,F}\) is relatively compact in
		\(\mathcal Y\).
	\end{lemma}

	\begin{proof}
		For each \(r\geq0\), the map \(t\mapsto P_\mu(r+t)F(t)\) is strongly
		measurable by \eqref{eq:operator-action-continuity} and the strong
		measurability of \(F\). By \eqref{eq:P-exact-norm} and the monotonicity
		of \(p_{\mu,\lambda_1}\) in Lemma~\ref{lem:pnu-properties},
		\[
		\|P_\mu(r+t)F(t)\|_X
		\leq p_{\mu,\lambda_1}(r+t)\mathcal G(t)
		\leq p_{\mu,\lambda_1}(t)\mathcal G(t).
		\]
		Hence the Bochner integral exists by
		\eqref{eq:attractor-history-envelope}, and
		\begin{equation}\label{eq:attractor-history-uniform-bound}
			\sup_F\sup_{r\geq0}\|\mathscr F_{\mu,F}(r)\|_X
			\leq K_{\mu,\mathcal G}.
		\end{equation}
		
		We first prove relative compactness of the values in \(X\).
		By \eqref{eq:attractor-scalar-integrability},
		\(\lambda_n p_{\mu,\lambda_n}(s)\leq s^{-1}\) for \(s>0\).
		Thus \eqref{eq:P-spectral-definition} gives
		\[
		\|P_\mu(s)v\|_{X_1}^2
		=\sum_{n=1}^{\infty}
		\lambda_n^2p_{\mu,\lambda_n}(s)^2|v_n|^2
		\leq s^{-2}\|v\|_X^2.
		\]
		Choose  \(s_0>0\)  such that \(|s-s_0|<s_0/2\). Then, we have
		\[
		\lambda_n^2|p_{\mu,\lambda_n}(s)-p_{\mu,\lambda_n}(s_0)|^2
		|v_n|^2\leq 9s_0^{-2}|v_n|^2.
		\]
		Since \(p_{\mu,\lambda_n}\) is continuous and
		\(\sum_{n=1}^{\infty}|v_n|^2<\infty\), dominated convergence gives 	\[
		\|P_\mu(s)v-P_\mu(s_0)v\|_{X_1}^2
		=
		\sum_{n=1}^{\infty}
		\lambda_n^2
		|p_{\mu,\lambda_n}(s)-p_{\mu,\lambda_n}(s_0)|^2|v_n|^2
		\to 0
		\qquad\text{as }s\to s_0.
		\]
		Thus \(P_\mu(\cdot)v\in C((0,\infty);X_1)\) for every \(v\in X\).
		
		Fix \(0<\varepsilon<M<\infty\) and define
		\[
		\mathscr F_{\mu,F}^{\varepsilon,M}(r)
		:=\int_\varepsilon^M P_\mu(r+t)F(t)\,dt .
		\]
		Choose \(X\)-valued simple functions \(F_k\to F\) almost everywhere.
		By the continuity of \(P_\mu(\cdot)v\) , \(t\mapsto P_\mu(r+t)F_k(t)\) is
		strongly measurable in \(X_1\) on \([\varepsilon,M]\). Moreover,
		\[
		\|P_\mu(r+t)(F_k(t)-F(t))\|_{X_1}
		\leq \varepsilon^{-1}\|F_k(t)-F(t)\|_X\to0
		\quad\text{for almost every }t\in[\varepsilon,M].
		\]
		Hence \(t\mapsto P_\mu(r+t)F(t)\) is strongly measurable in \(X_1\).
		By \eqref{eq-probability measure},
		\(p_{\mu,\lambda_1}(t)\geq p_{\mu,\lambda_1}(M)>0\) for \(0\leq t\leq M\). Consequently,
		\[
		\begin{aligned}
			\int_\varepsilon^M\|P_\mu(r+t)F(t)\|_{X_1}\,dt
			&\leq \frac1\varepsilon\int_\varepsilon^M\mathcal G(t)\,dt
			\leq \frac{1}{\varepsilon p_{\mu,\lambda_1}(M)}
			\int_\varepsilon^M p_{\mu,\lambda_1}(t)\mathcal G(t)\,dt
			\leq \frac{K_{\mu,\mathcal G}}{\varepsilon p_{\mu,\lambda_1}(M)} .
		\end{aligned}
		\]
		Thus \(\mathscr V_{\varepsilon,M}:=
		\{\mathscr F_{\mu,F}^{\varepsilon,M}(r):F\text{ as above},\ r\geq0\}\)
		is bounded in \(X_1\) and relatively compact in \(X\) by the compact
		embedding \(X_1\hookrightarrow X\) in Section~\ref{sec:preliminaries}.
		
		For the full integrals, \eqref{eq:attractor-history-envelope} gives
		\[
		\begin{aligned}
			&\sup_F\sup_{r\geq0}
			\|\mathscr F_{\mu,F}(r)-\mathscr F_{\mu,F}^{\varepsilon,M}(r)\|_X
			\leq
			\int_0^\varepsilon p_{\mu,\lambda_1}(t)\mathcal G(t)\,dt
			+\int_M^\infty p_{\mu,\lambda_1}(t)\mathcal G(t)\,dt
			\to0
		\end{aligned}
		\]
		as \(\varepsilon\to0\) and \(M\to\infty\).
		Given \(\alpha>0\), choose \(\varepsilon,M\) such that
		\(
		\sup_F\sup_{r\geq0}
		\|\mathscr{F}_{\mu,F}(r)-\mathscr{F}_{\mu,F}^{\varepsilon,M}(r)\|_X
		<\frac\alpha2.
		\)
		Since \(\mathscr V_{\varepsilon,M}\) is relatively compact,
		it provides finitely many points \(y_1,\ldots,y_m\in X\)
		such that
		\(
		\mathscr V_{\varepsilon,M}
		\subseteq
		\bigcup_{\ell=1}^m B_X(y_\ell,\alpha/2).
		\)
		Choose \(\ell\) with
		\(\|\mathscr{F}_{\mu,F}^{\varepsilon,M}(r)-y_\ell\|_X<\alpha/2\).
		Then
		\[
		\begin{aligned}
			\|\mathscr{F}_{\mu,F}(r)-y_\ell\|_X
			\leq
			\|\mathscr{F}_{\mu,F}(r)-\mathscr{F}_{\mu,F}^{\varepsilon,M}(r)\|_X
			+
			\|\mathscr{F}_{\mu,F}^{\varepsilon,M}(r)-y_\ell\|_X
			<\alpha.
		\end{aligned}
		\]	Therefore, since \(\alpha>0\) is arbitrary, it follows that
		\(
		\mathscr V
		:=
		\{\mathscr{F}_{\mu,F}(r):F\text{ as above},\ r\geq0\}\subseteq
		\bigcup_{\ell=1}^m B_X(y_\ell,\alpha)
		\)
		is totally bounded in \(X\).
		Completeness of \(X\) then implies that  \(\mathscr V\)  is
		relatively compact in \(X\).

		We next prove equicontinuity. By
		Lemma~\ref{lem:attractor-kernel-estimates}\textup{(iii),(iv)},
		for \(t>0\) and \(r,h\geq0\),
		\[
		\|P_\mu(r+t+h)-P_\mu(r+t)\|_{\mathcal L(X)}
		\leq\|P_\mu(t+h)-P_\mu(t)\|_{\mathcal L(X)}
		\leq p_{\mu,\lambda_1}(t),
		\]
		and \(\|P_\mu(t+h)-P_\mu(t)\|_{\mathcal L(X)}\to0\) as \(h\to0\).
		As \(h\to0\), Dominated convergence and \eqref{eq:attractor-history-envelope} yield
		\begin{equation*}
			\begin{aligned}
				&\sup_F\sup_{r\geq0}
				\|\mathscr F_{\mu,F}(r+h)-\mathscr F_{\mu,F}(r)\|_X
				\leq\int_0^\infty
				\|P_\mu(t+h)-P_\mu(t)\|_{\mathcal L(X)}\mathcal G(t)\,dt
				\to0.
			\end{aligned}
		\end{equation*}
		This proves equicontinuity on \([0,\infty)\).
		By \eqref{eq:attractor-history-uniform-bound} and \eqref{eq:C_b},
		\(\mathscr F_{\mu,F}\in C_b([0,\infty);X)\subset\mathcal Y\).
		
		Finally, by the Arzel\`a--Ascoli theorem on each \([0,N]\) and a
		diagonal argument, every sequence \(\{\mathscr F_j\}_{j\ge1}\) in
		the family has a subsequence \(\{\mathscr F_{j_k}\}_{k\ge1}\)  and a limit
		\(\mathscr F_*\in C([0,\infty);X)\). Thus,
		\[
		\sup_{0\leq r\leq N}\|\mathscr F_j(r)-\mathscr F_*(r)\|_X
		\to0,\qquad N\in\mathbb N .
		\]
		Passing to the pointwise limit in
		\eqref{eq:attractor-history-uniform-bound} gives
		\(\|\mathscr F_*\|_\infty\leq K_{\mu,\mathcal G}\), so
		\(\mathscr F_*\in\mathcal Y\) by \eqref{eq:C_b}.
		For each \(J\in\mathbb N\), \eqref{eq:profile-norm} gives
		\[
		\begin{aligned}
			\|\mathscr F_j-\mathscr F_*\|_{\mathcal Y}
			&\leq 2\sup_{0\leq r\leq J}\|\mathscr F_j(r)-\mathscr F_*(r)\|_X
			+2^{1-J}K_{\mu,\mathcal G}.
		\end{aligned}
		\]
		Thus \(\limsup_{j\to\infty}\|\mathscr F_j-\mathscr F_*\|_{\mathcal Y}
		\leq 2^{1-J}K_{\mu,\mathcal G}\). Letting \(J\to\infty\) proves
		convergence in \(\mathcal Y\) and hence the 	\(\{\mathscr{F}_{\mu,F}\}\)   is relatively compact
		in \(\mathcal Y\).
	\end{proof}
	
	\begin{lemma}
		\label{prop:attractor-physical-asymptotic-compactness}
		Let \(B\in\mathfrak D_{\mu,\varrho}\), and fix
		\(\omega\in\mathcal O_{1,\mu}\) such that
		\eqref{eq:admissible-random-initial-set} holds.
		If \(T_j\to\infty\) and
		\(x_j\in B(\theta_{-T_j}\omega)\) for \(j\geq1\),
		then the sequence
		\(
		\{\Gamma_\mu(T_j,\theta_{-T_j}\omega,\Psi_{x_j})\}_{j\geq1}
		\)
		has a convergent subsequence in \(\mathcal Y\).
	\end{lemma}

	\begin{proof}
		For each \(j\geq1\), define
		\[
		F_j(a):=
		\begin{cases}
			f(u_\mu (T_j-a,\theta_{-T_j}\omega,0)), & 0\leq a\leq T_j,\\
			0, & a>T_j.
		\end{cases}
		\]
		Each \(F_j\) is strongly measurable by the continuity of \(u\) and \(f\).
		Taking \(x=0\) in \eqref{eq:attractor-profile-decomposition} gives
		\begin{equation}\label{eq:D-S}
			\Gamma_\mu(T_j,\theta_{-T_j}\omega,\Psi_0)
			=\mathscr F_{\mu,F_j}+S_{\mu,T_j}(\omega).
		\end{equation}
		
		By \eqref{eq:f-lipschitz} and \eqref{eq:attractor-state-comparison}, for \(0\leq a\leq T_j\),
		\[
		\|F_j(a)\|_X
		\leq \|f(0)\|_X+L\|u_\mu (T_j-a,\theta_{-T_j}\omega,0)\|_X
		\leq \|f(0)\|_X+LR_{0,\mu}(\theta_{-a}\omega)
		=:\mathcal G^*(a).
		\]
		The same bound holds for \(a>T_j\) since \(F_j(a)=0\).
		By \eqref{eq:attractor-scalar-integrability},
		\eqref{eq:attractor-profile-radius}, and
		Lemma~\ref{lem:attractor-invariant-good-set},
		\[
		\int_0^\infty p_{\mu,\lambda_1}(a)\mathcal G^*(a)\,da
		\leq \frac{\|f(0)\|_X}{\lambda_1}
		+L\int_0^\infty p_{\mu,\lambda_1}(a)
		R_{0,\mu}(\theta_{-a}\omega)\,da
		\leq R_{1,\mu}(\omega)<\infty.
		\]
		Thus \(\mathcal G^*\) satisfies \eqref{eq:attractor-history-envelope}.
		By Lemma~\ref{lem:attractor-history-compactness}, there exist
		a subsequence \(\{{j_k}\}_{k\ge1}\)  and \(\mathscr F^*\in\mathcal Y\) such that
		\[
		\mathscr F_{\mu,F_{j_k}}\to\mathscr F^*
		\qquad\text{in }\mathcal Y.
		\]
		Since \(T_j\to\infty\), \eqref{eq:past-noise-weighted-convergence}
		also gives \(S_{\mu,T_j}(\omega)\to W_\mu(\omega)\) in \(\mathcal Y\).
		
		Finally, by \eqref{eq:admissible-random-initial-set}, \eqref{eq:D-S}, \eqref{eq:attractor-profile-difference},
		and \(x_{j_k}\in B(\theta_{-T_{j_k}}\omega)\), we obtain
		\[
		\begin{aligned}
			&\bigl\|
			\Gamma_\mu(T_{j_k},\theta_{-T_{j_k}}\omega,\Psi_{x_{j_k}})
			-(\mathscr F^*+W_\mu(\omega))
			\bigr\|_{\mathcal Y}\\
			&\quad\leq
			2p_{\mu,\varrho}(T_{j_k})r_B(\theta_{-T_{j_k}}\omega)
			+\|\mathscr F_{\mu,F_{j_k}}-\mathscr F^*\|_{\mathcal Y}
			+\|S_{\mu,T_{j_k}}(\omega)-W_\mu(\omega)\|_{\mathcal Y}
			\to0.
		\end{aligned}
		\]
	\end{proof}

	\begin{lemma}\label{lem:attractor-zero-orbit-continuity}
		The mapping
		\((T,\omega)\mapsto\Gamma_\mu(T,\theta_{-T}\omega,\Psi_0)\)
		is jointly measurable. For each \(\omega\in\mathcal O_{1,\mu}\),
		\(T\mapsto \Gamma_\mu(T,\theta_{-T}\omega,\Psi_0)\)
		is continuous from \([0,\infty)\) into \(\mathcal Y\).
		For every
		integer \(n\geq1\), the set
		\begin{equation}\label{eq:attractor-compact-tail}
			\varkappa_{\mu,n}(\omega):=
			\overline{
				\{\Gamma_\mu(T,\theta_{-T}\omega,\Psi_0):T\geq n\}
			}^{\,\mathcal Y}
		\end{equation}
		is nonempty and compact.
	\end{lemma}
	
	\begin{proof}
		By Lemma~\ref{prop:lifted-map-measurability} and the
		joint measurability of the Wiener shifts, the mapping
		\(
		(T,\omega)\mapsto
		\Gamma_\mu(T,\theta_{-T}\omega,\Psi_0)
		\)
		is jointly measurable.
		Fix \(T^*>0\) and \(\omega\in\mathcal O_{1,\mu}\).
		For \(0<q\leq h\leq1\) and \(|s|\leq T^*\),
		\eqref{eq:local Holder} gives
		\(
		[\theta_s\omega]_{\gamma;[0,h]}
		\leq[\omega]_{\gamma;[-T^*,T^*+1]}.
		\)
		Hence \eqref{eq:two-parameter-bound} yield
		\begin{equation}\label{eq:Z_mu-Ch}
			\begin{aligned}
				\sup_{r\geq0}\|Z_\mu(q,r,\theta_s\omega)\|_X
				\leq C\|G\|_{\mathcal L(X,X_\delta)}
				[\theta_s\omega]_{\gamma;[0,q]}
				\bigl(q^\gamma+q^{\rho+\gamma}\bigr)\leq Ch^\gamma,
				\qquad |s|\leq T^*.
			\end{aligned}
		\end{equation}		
		Using \eqref{eq:P-exact-norm},
		(\ref{ass:nonlinearity}), and
		\eqref{eq:Z_mu-Ch} in \eqref{eq:mild-stochastic}
		with zero initial value, we obtain
		\[
		\|u_\mu (q,\theta_s\omega;0)\|_X
		\leq\|f(0)\|_Xh
		+L\int_0^q\|u_\mu (a,\theta_s\omega;0)\|_X\,da
		+Ch^\gamma,
		\qquad 0\leq q\leq h.
		\]
		Using Gronwall's inequality gives
		\[
		\sup_{0\leq q\leq h}\|u_\mu (q,\theta_s\omega;0)\|_X
		\leq\bigl(\|f(0)\|_Xh+Ch^\gamma\bigr)e^{Lh}.
		\]
		Since \(\Psi_0\equiv0\), it follows from \eqref{eq:lifted-map},
		\eqref{eq:P-exact-norm}, \eqref{eq:f-lipschitz},
		\eqref{eq:Z_mu-Ch}, and \eqref{eq:C_b} that
		
		\begin{equation}\label{eq:Gamma_mu-Ch}
			\begin{aligned}
				\sup_{|s|\leq T^*}
				\|\Gamma_\mu(h,\theta_s\omega,\Psi_0)\|_{\mathcal Y}
				&\leq2h\left[
				\|f(0)\|_X
				+L\bigl(\|f(0)\|_Xh+Ch^\gamma\bigr)e^{Lh}
				\right]+2Ch^\gamma\to0.
			\end{aligned}
		\end{equation}
		as \(h\to0\).
		Let \(0\leq T\leq T'\leq T^*\) with \(h=T'-T\leq1\).
		Since \(\theta_h\theta_{-T'}\omega=\theta_{-T}\omega\),
		\eqref{eq:lifted-cocycle} gives
		\[
		\Gamma_\mu(T',\theta_{-T'}\omega,\Psi_0)
		=\Gamma_\mu\bigl(T,\theta_{-T}\omega,
		\Gamma_\mu(h,\theta_{-T'}\omega,\Psi_0)\bigr).
		\]
		Thus, as \(h\to 0\), by \eqref{eq:lifted-map-lipschitz} and \eqref{eq:Gamma_mu-Ch},
		\[
		\begin{aligned}
			&\|\Gamma_\mu(T',\theta_{-T'}\omega,\Psi_0)
			-\Gamma_\mu(T,\theta_{-T}\omega,\Psi_0)\|_{\mathcal Y}
			\leq C_{T^*}
			\|\Gamma_\mu(h,\theta_{-T'}\omega,\Psi_0)\|_{\mathcal Y}
			\to0
			.
		\end{aligned}
		\]
		Since \(T^*>0\) is arbitrary, the mapping
		\(T\mapsto \Gamma_\mu(T,\theta_{-T}\omega,\Psi_0)\)
		is continuous from \([0,\infty)\) into \(\mathcal Y\).
		
		Finally, fix \(n\geq1\) and let \(T_j\geq n\).
		If \(\{T_{j}\}_{j\ge1}\)  is bounded, there exist a subsequence
		\(\{T_{j_k}\}_{k\ge1}\) and \(T_*\geq n\) such that
		\(T_{j_k}\to T_*\).
		By the continuity established above,
		\[
		\Gamma_\mu(T_{j_k},\theta_{-T_{j_k}}\omega,\Psi_0)
		\to 
		\Gamma_\mu(T_*,\theta_{-T_*}\omega,\Psi_0)
		\qquad\text{in }\mathcal Y.
		\]
		If  \(\{T_{j}\}_{j\ge1}\) is unbounded, choose a subsequence
		\(\{T_{j_k}\}_{k\ge1}\) tending to infinity.
		The constant random set \(B_0(\omega)=\{0\}\) belongs to
		\(\mathfrak D_{\mu,\varrho}\) and satisfies
		\eqref{eq:admissible-random-initial-set}, since its radius
		is identically zero.
		Applying
		Lemma~\ref{prop:attractor-physical-asymptotic-compactness}
		with \(B=B_0\) and zero initial values, we obtain every sequence in 
		\(
		\{\Gamma_\mu(T,\theta_{-T}\omega,\Psi_0):T\geq n\}
		\)
		has a subsequence converging in \(\mathcal Y\).
		Hence this set is relatively compact, and its closure
		\(\varkappa_{\mu,n}(\omega)\) is compact.
		It is nonempty since 
		\(\Gamma_\mu(n,\theta_{-n}\omega,\Psi_0)\in \{\Gamma_\mu(T,\theta_{-T}\omega,\Psi_0):T\geq n\}\).
	\end{proof}

	\subsection{The minimal random pullback attractor}
	\label{subsec:attractor-construction}
	\begin{theorem}\label{thm:generalized-random-pullback-attractor}
		Fix \(\mu\geq0\). Suppose that
		(\ref{ass:kernel}), (\ref{ass:kernel-holder}),
		(\ref{ass:noise-coefficient}), and (\ref{ass:nonlinearity})
		hold, with \(L<\lambda_1\), and set
		\(\varrho:=\lambda_1-L>0\).
		Let \(\Gamma_\mu\) be defined by
		\eqref{eq:lifted-map}, and define
		\begin{equation}\label{eq:attractor-definition}
			\mathcal A_\mu(\omega):=
			\bigcap_{n=1}^{\infty}
			\overline{
				\{\Gamma_\mu(T,\theta_{-T}\omega,\Psi_0):T\geq n\}
			}^{\,\mathcal Y},
			\qquad \omega\in\mathcal O_{1,\mu}.
		\end{equation}
		Then \(\mathcal A_\mu\) is a nonempty random compact set
		in \(\mathcal Y\) and is a minimal random pullback
		attractor for \(\Gamma_\mu\) relative to
		\(\mathfrak D_{\widetilde B}\).
		It has the following properties.

		\smallskip
		\noindent\textup{(i)} For all \(s\geq0\) and \(\omega\in\mathcal O_{1,\mu}\),
		\begin{equation}\label{eq:attractor-strict-invariance}
			\Gamma_\mu(s,\omega,\mathcal A_\mu(\omega))
			=\mathcal A_\mu(\theta_s\omega).
		\end{equation}
		
		\smallskip
		\noindent\textup{(ii)} For every \(B\in\mathfrak D_{\mu,\varrho}\) and
		every \(\omega\in\mathcal O_{1,\mu}\) satisfying
		\eqref{eq:admissible-random-initial-set},
		\begin{equation*}
			\lim_{T\to\infty}\operatorname{dist}_{\mathcal Y}
			\bigl(\Gamma_\mu(T,\theta_{-T}\omega,\widetilde B(\theta_{-T}\omega)),
			\mathcal A_\mu(\omega)\bigr)=0.
		\end{equation*}
		
		\smallskip
		\noindent\textup{(iii)} For every random pullback attractor
		\(\mathcal K\) for \(\Gamma_\mu\) relative to \(\mathfrak D_{\widetilde B}\),
		we have \(\mathcal A_\mu(\omega)\subseteq\mathcal K(\omega)\)
		almost surely.
		
		\smallskip
		\noindent\textup{(iv)} For every \(\omega\in\mathcal O_{1,\mu}\),
		\begin{equation}\label{eq:attractor-size-bound}
			\sup_{\mathscr T\in\mathcal A_\mu(\omega)}
			\|\mathscr T\|_{\mathcal Y}\leq 2R_{1,\mu}(\omega).
		\end{equation}
	\end{theorem}
	
	\begin{proof}
		
		Fix \(\omega\in\mathcal O_{1,\mu}\).
		By Lemma~\ref{lem:attractor-zero-orbit-continuity} and
		\eqref{eq:attractor-compact-tail}, the sets
		\(\varkappa_{\mu,n}(\omega)\) are nonempty, compact, and decreasing
		in \(n\). Hence their intersection \(\mathcal A_\mu(\omega)\) in
		\eqref{eq:attractor-definition} is nonempty and compact.
		
		If \(T_j\to\infty\), then \(T_j\geq n\) eventually for every \(n\).
		Since \(
		\Gamma_\mu(T_j,\theta_{-T_j}\omega,\Psi_0)
		\to \mathscr T\) and \(\varkappa_{\mu,n}(\omega)\) is closed, 	it follows that \(\mathscr T\in\varkappa_{\mu,n}(\omega)\).
		Conversely, if \(\mathscr T\in\mathcal A_\mu(\omega)\), then
		\(\mathscr T\in\varkappa_{\mu,j}(\omega)\) for every \(j\), so we
		may choose \(T_j\geq j\) such that
		\(
		\|\Gamma_\mu(T_j,\theta_{-T_j}\omega,\Psi_0)
		-\mathscr T\|_{\mathcal Y}<1/j.
		\)
		Consequently,
		\begin{equation}\label{eq:attractor-sequential-characterization}
			\mathcal A_\mu(\omega)=
			\left\{
			\mathscr{T}\in\mathcal Y:
			\begin{array}{l}
				\text{there exist }T_j\to\infty\text{ such that}\\
				\Gamma_\mu(T_j,\theta_{-T_j}\omega,\Psi_0)
				\to\mathscr{T}\text{ in }\mathcal Y
			\end{array}
			\right\}.
		\end{equation}
		To prove measurability, fix \(\Psi\in\mathcal Y\). 
		By \eqref{eq:attractor-compact-tail} and continuity in \(T\)
		from Lemma~\ref{lem:attractor-zero-orbit-continuity},
		\[
		d_n(\omega)
		:=\operatorname{dist}_{\mathcal Y}
		(\{\Psi\},\varkappa_{\mu,n}(\omega))
		= \inf_{q\in\mathbb Q\cap[n,\infty)}
		\|\Psi-\Gamma_\mu(q,\theta_{-q}\omega,\Psi_0)\|_{\mathcal Y}.
		\]
		Lemma~\ref{lem:attractor-zero-orbit-continuity} gives measurability
		of \(\omega\mapsto \Gamma_\mu(q,\theta_{-q}\omega,\Psi_0)\) for each fixed \(q\), so \(d_n\) is
		measurable. 
		Moreover, we have
		\[
		d_n(\omega)\leq d_{n+1}(\omega)
		\leq\operatorname{dist}_{\mathcal Y}
		(\{\Psi\},\mathcal A_\mu(\omega)).
		\]
		
		By compactness, choose \(\mathscr P_n\in\varkappa_{\mu,n}(\omega)\)
		such that \(\|\Psi-\mathscr P_n\|_{\mathcal Y}=d_n(\omega)\).
		Compactness of \(\varkappa_{\mu,1}(\omega)\)
		gives a  subsequence 	\(\mathscr P_{n_k}\to\mathscr P\) and the  limit \(\mathscr P\) belongs to  \(\mathcal A_\mu(\omega)\).  This gives
		\begin{equation*}
			\operatorname{dist}_{\mathcal Y}(\{\Psi\},\mathcal A_\mu(\omega))
			\leq \|\Psi-\mathscr P\|_{\mathcal Y}
			=\lim_{n\to\infty}d_n(\omega).
		\end{equation*}
		Thus, we have 
		\begin{equation}\label{eq:attractor-measurable-distance}
			\operatorname{dist}_{\mathcal Y}(\{\Psi\},\mathcal A_\mu(\omega))
			=\lim_{n\to\infty}d_n(\omega).
		\end{equation}
		Since each \(d_n\) is measurable, this identity shows that
		the distance to \(\mathcal A_\mu(\omega)\) is measurable
		for every fixed \(\Psi\in\mathcal Y\).
		Thus, by Definition~\ref{def:attractor-random-closed-set},
		\(\mathcal A_\mu\) is a nonempty random compact set on
		\((\mathcal O_{1,\mu},\mathcal F_{1,\mu},\mathbb P)\).

		\textup{(i)}
		Fix \(s\geq0\). By Lemma~\ref{lem:attractor-invariant-good-set},
		\(\theta_s\omega\in\mathcal O_{1,\mu}\). The  
		\eqref{eq:lifted-cocycle} gives
		\begin{equation}\label{eq:attractor-zero-orbit-shift}
			\Gamma_\mu\bigl(s,\omega,\Gamma_\mu(T,\theta_{-T}\omega,\Psi_0)\bigr)
			=\Gamma_\mu\bigl(T+s,\theta_{-(T+s)}(\theta_s\omega),\Psi_0\bigr).
		\end{equation}
		Using Lemma~\ref{prop:lifted-map-continuity},
		\eqref{eq:attractor-zero-orbit-shift}, and
		\eqref{eq:attractor-sequential-characterization},
		we obtain
		\[
		\Gamma_\mu(s,\omega,\mathscr T)
		=\lim_{j\to\infty}
		\Gamma_\mu\bigl(
		s,\omega,
		\Gamma_\mu(T_j,\theta_{-T_j}\omega,\Psi_0)
		\bigr)
		=
		\lim_{j\to\infty}
		\Gamma_\mu\bigl(
		T_j+s,\theta_{-(T_j+s)}(\theta_s\omega),\Psi_0
		\bigr)
		\in\mathcal A_\mu(\theta_s\omega).
		\]
		Hence
		\(
		\Gamma_\mu(s,\omega,\mathcal A_\mu(\omega))
		\subseteq\mathcal A_\mu(\theta_s\omega).
		\)
		
		Conversely, for \(\widetilde{\mathscr T}\in\mathcal A_\mu(\theta_s\omega)\),
		choose \(T_j\to\infty\), \(T_j\geq s+1\), such that
		\(\Gamma_\mu(T_j,\theta_{s-T_j}\omega,\Psi_0)\to\widetilde{\mathscr T}\).
		Since \(T_j-s\to\infty\), compactness of \(\varkappa_{\mu,1}(\omega)\)
		and \eqref{eq:attractor-sequential-characterization} give a subsequence
		such that
		\[
		\Gamma_\mu(T_{j_k}-s,\theta_{s-T_{j_k}}\omega,\Psi_0)
		\to\mathscr T\in\mathcal A_\mu(\omega).
		\]
		Applying \(\Gamma_\mu(s,\omega,\cdot)\) and using
		\eqref{eq:attractor-zero-orbit-shift} gives
		\(\Gamma_\mu(s,\omega,\mathscr T)=\widetilde{\mathscr T}\).
		This proves \eqref{eq:attractor-strict-invariance}.
		
		\textup{(ii)}
		Fix \(\omega\in\mathcal O_{1,\mu}\). We first show that
		\begin{equation}\label{eq:attractor-zero-orbit-attraction}
			\operatorname{dist}_{\mathcal Y}
			\bigl(\{\Gamma_\mu(T,\theta_{-T}\omega,\Psi_0)\},
			\mathcal A_\mu(\omega)\bigr)
			\to0, \quad\text{as}
			\qquad T\to\infty.
		\end{equation}
		Suppose that
		\eqref{eq:attractor-zero-orbit-attraction} does not hold.
		There exists \(\varepsilon>0\)and  \(j\geq1\) such that
		\begin{equation*}
			\operatorname{dist}_{\mathcal Y}
			\bigl(
			\{\Gamma_\mu(T_j,\theta_{-T_j}\omega,\Psi_0)\},
			\mathcal A_\mu(\omega)
			\bigr)
			\geq\varepsilon, \quad T_j\geq j.
		\end{equation*}
		Then \(T_j\to\infty\) and \(T_j\geq1\) for every \(j\).
		By  
		\eqref{eq:attractor-compact-tail},
		\(
		\Gamma_\mu(T_j,\theta_{-T_j}\omega,\Psi_0)
		\in
		\varkappa_{\mu,1}(\omega)
		\), which is compact by
		Lemma~\ref{lem:attractor-zero-orbit-continuity}.
		Hence a subsequence converges in \(\mathcal Y\) to some \(\mathscr T\).
		Since \(T_{j_k}\to\infty\), \eqref{eq:attractor-sequential-characterization}
		gives \(\mathscr T\in\mathcal A_\mu(\omega)\). Thus
		\begin{equation*}
			\begin{aligned}
				\varepsilon
				\leq \operatorname{dist}_{\mathcal Y}
				\bigl(\{\Gamma_\mu(T_{j_k},\theta_{-T_{j_k}}\omega,\Psi_0)\},
				\mathcal A_\mu(\omega)\bigr)
				\leq
				\|\Gamma_\mu(T_{j_k},\theta_{-T_{j_k}}\omega,\Psi_0)
				-\mathscr T\|_{\mathcal Y}
				\to0,
			\end{aligned}
		\end{equation*}
		a contradiction. 
		This proves \eqref{eq:attractor-zero-orbit-attraction}.

		Now let \(B\in\mathfrak D_{\mu,\varrho}\) and suppose that
		\eqref{eq:admissible-random-initial-set} holds at \(\omega\).
		Using \eqref{eq:attractor-profile-difference} with \(y=0\) and applying
		\eqref{eq:widetilde B},	\eqref{eq:admissible-random-initial-set}
		and	\eqref{eq:attractor-zero-orbit-attraction}, we obtain
		\[
		\begin{aligned}
			&\operatorname{dist}_{\mathcal Y}
			\bigl(
			\Gamma_\mu(T,\theta_{-T}\omega,
			\widetilde B(\theta_{-T}\omega)),
			\mathcal A_\mu(\omega)
			\bigr)\\
			&\leq
			\sup_{x\in B(\theta_{-T}\omega)}	\|\Gamma_\mu(T,\theta_{-T}\omega,\Psi_{\mu,x})
			-\Gamma_\mu(T,\theta_{-T}\omega,\Psi_0)\|_{\mathcal Y}
			+
			\operatorname{dist}_{\mathcal Y}
			\bigl(
			\{\Gamma_\mu(T,\theta_{-T}\omega,\Psi_0)\},
			\mathcal A_\mu(\omega)
			\bigr)\\
			&\leq
			2p_{\mu,\varrho}(T)r_B(\theta_{-T}\omega)+
			\operatorname{dist}_{\mathcal Y}
			\bigl(
			\{\Gamma_\mu(T,\theta_{-T}\omega,\Psi_0)\},
			\mathcal A_\mu(\omega)
			\bigr)\to0, \quad T\to\infty.
		\end{aligned}
		\]

		\textup{(iii)}
		Let \(\mathcal K\) be a  random pullback attractor relative to
		\(\mathfrak D_{\widetilde B}\). By \eqref{eq:attractor-initial-family},
		the constant family \(\{\Psi_0\}\) belongs to \(\mathfrak D_{\widetilde B}\),
		the pullback attraction property of \(\mathcal K\) gives
		\[
		\operatorname{dist}_{\mathcal Y}
		\bigl(
		\{\Gamma_\mu(T,\theta_{-T}\omega,\Psi_0)\},
		\mathcal K(\omega)
		\bigr)
		\to 0, \quad T\to\infty.
		\]
		Fix such an \(\omega\in\mathcal O_{1,\mu}\).
		For any \(\mathscr T\in\mathcal A_\mu(\omega)\), choose \(T_j\to\infty\)
		as in \eqref{eq:attractor-sequential-characterization}.
		By continuity of the distance to \(\mathcal K(\omega)\),
		\[
		\operatorname{dist}_{\mathcal Y}(\{\mathscr T\},\mathcal K(\omega))
		=\lim_{j\to\infty}\operatorname{dist}_{\mathcal Y}
		\bigl(\{\Gamma_\mu(T_j,\theta_{-T_j}\omega,\Psi_0)\},\mathcal K(\omega)\bigr)
		=0.
		\]
		Closedness of \(\mathcal K(\omega)\) gives
		\(\mathcal A_\mu(\omega)\subseteq\mathcal K(\omega)\).
		If \(\mathcal K\) is also minimal, its minimality gives
		\(
		\mathcal K(\omega)\subseteq\mathcal A_\mu(\omega)\)
		almost surely.
		So the minimal attractor is unique almost surely.

		\textup{(iv)}
		By \eqref{eq:attractor-profile-bound} with \(x=0\),
		\[
		\|\Gamma_\mu(T,\theta_{-T}\omega,\Psi_0)\|_{\mathcal Y}
		\leq 2R_{1,\mu}(\omega),\qquad T\geq0.
		\]
		Passing to the limit in \eqref{eq:attractor-sequential-characterization}
		proves \eqref{eq:attractor-size-bound}.
		
	\end{proof}
	
	\begin{remark}\label{rem:attractor-strong-minimality}
		Property \textup{(iii)} extends to random closed sets
		that pullback attract every member of
		\(\mathfrak D_{\widetilde B}\).
		More precisely, if \(\mathcal C\) is such a set, then
		\(
		\mathcal A_\mu(\omega)\subseteq\mathcal C(\omega)
		\) almost surely.
		The same proof applies to \(\mathcal C\), without
		requiring compactness or strict invariance.
		Moreover,  extending \textup{(iv)}
		by \eqref{eq:attractor-radius-local-moments}
		and
		combining this with \eqref{eq:attractor-size-bound},
		we obtain
		\(
		\mathbb E\!\left[
		\sup_{\mathscr T\in\mathcal A_\mu(\omega)}
		\|\mathscr T\|_{\mathcal Y}
		\right]
		\leq2\mathbb E R_{1,\mu}<\infty.
		\)
	\end{remark}

	\begin{corollary}\label{cor:attractor-tempered}
		There exists a measurable invariant set
		\(\mathcal O_{*,\mu}\subseteq\mathcal O_{1,\mu}\) with
		\(\mathbb P(\mathcal O_{*,\mu})=1\) such that
		\(1+R_{1,\mu}\) is tempered; that is,
		\begin{equation}\label{eq:attractor-radius-tempered}
			\lim_{|s|\to\infty}
			e^{-\varepsilon|s|}
			\bigl(1+R_{1,\mu}(\theta_s\omega)\bigr)=0,
			\qquad
			\varepsilon>0,\quad \omega\in\mathcal O_{*,\mu}.
		\end{equation}
		The random ball \(B_{R,\mu}\) in
		\eqref{eq:attractor-absorbing-ball} has a tempered
		radius, and the attractor \(\mathcal A_\mu\) admits a
		tempered random bound.
	\end{corollary}

	\begin{proof}

		For \(\omega\in\mathcal O_{1,\mu}\), define
		\[
		J(\omega):=\sup_{0\le s\le1}
		\bigl(1+R_{1,\mu}(\theta_s\omega)\bigr)
		=\sup_{s\in\mathbb Q\cap[0,1]}
		\bigl(1+R_{1,\mu}(\theta_s\omega)\bigr)
		\]
		by  
		Lemma~\ref{lem:attractor-invariant-good-set}.
		Thus \(J\) is measurable, and
		\eqref{eq:attractor-radius-local-moments} gives \(\mathbb EJ<\infty\).
		Define
		\[
		\begin{aligned}
			\mathcal O_{*,\mu}
			&:=\left\{\omega\in\mathcal O_{1,\mu}:
			\lim_{|s|\to\infty}e^{-\varepsilon|s|}
			\bigl(1+R_{1,\mu}(\theta_s\omega)\bigr)=0
			\text{ for every }\varepsilon>0\right\}.
		\end{aligned}
		\]
		For \(n\in\mathbb Z\), we have
		\(
		J(\theta_n\omega)
		=\sup_{n\le s\le n+1}\bigl(1+R_{1,\mu}(\theta_s\omega)\bigr).
		\) Since \(\bigl||s|-|n|\bigr|\le1\) on \([n,n+1]\),  we have
		\[
		e^{-q(|n|+1)}J(\theta_n\omega)
		\le
		\sup_{n\le s\le n+1}
		e^{-q|s|}\bigl(1+R_{1,\mu}(\theta_s\omega)\bigr)
		\le
		e^{-q(|n|-1)}J(\theta_n\omega)
		\]
		for \(q>0\). Together with \(e^{-\varepsilon|s|}\le e^{-q|s|}\) for
		\(0<q<\varepsilon\), this gives
		\[
		\mathcal O_{*,\mu}
		=
		\bigcap_{q\in\mathbb Q\cap(0,\infty)}
		\left\{\omega\in\mathcal O_{1,\mu}:
		\lim_{|n|\to\infty}e^{-q|n|}J(\theta_n\omega)=0\right\}.
		\]	
		Since \(\omega\mapsto J(\theta_n\omega)\) is measurable,
		each convergence event in the countable intersection is measurable.
		Hence \(\mathcal O_{*,\mu}\) is measurable.

		For each rational \(q>0\), Markov's inequality and the invariance
		of \(\mathbb P\) give
		\[
		\sum_{n\in\mathbb Z}
		\mathbb P\bigl(J(\theta_n\omega)>e^{q|n|/2}\bigr)
		\le \mathbb E J\sum_{n\in\mathbb Z}e^{-q|n|/2}<\infty.
		\]
		By the Borel--Cantelli lemma, almost surely,
		\(e^{-q|n|}J(\theta_n\omega)\le e^{-q|n|/2}\) for all sufficiently
		large \(|n|\). Thus each event in the countable intersection above
		has probability one, and \(\mathbb P(\mathcal O_{*,\mu})=1\).

		For \(\omega\in\mathcal O_{*,\mu}\), \(\tau\in\mathbb R\), and
		\(\varepsilon>0\), the shift identity gives
		\[
		e^{-\varepsilon|s|}
		\bigl(1+R_{1,\mu}(\theta_s\theta_\tau\omega)\bigr)
		\leq e^{\varepsilon|\tau|}e^{-\varepsilon|s+\tau|}
		\bigl(1+R_{1,\mu}(\theta_{s+\tau}\omega)\bigr)
		\to0,
		\]
		as \(|s|\to\infty\). Since \(\mathcal O_{1,\mu}\) is invariant by
		Lemma~\ref{lem:attractor-invariant-good-set}, this shows that
		\(\theta_\tau\mathcal O_{*,\mu}\subseteq\mathcal O_{*,\mu}\).
		Applying the same inclusion with \(-\tau\) gives the reverse inclusion.
		Thus \(\theta_\tau\mathcal O_{*,\mu}=\mathcal O_{*,\mu}\)
		.

		Finally, 
		by \eqref{eq:attractor-radius-tempered}, 
		\[
		\begin{aligned}
			0
			\leq e^{-\varepsilon|s|}
			\bigl(2R_{1,\mu}(\theta_s\omega)+1\bigr)
			\leq 2e^{-\varepsilon|s|}
			\bigl(1+R_{1,\mu}(\theta_s\omega)\bigr)
			\to 0
			\qquad\text{as }|s|\to\infty.
		\end{aligned}
		\]
		Thus, the radius of
		\(B_{R,\mu}(\omega)\) in
		\eqref{eq:attractor-absorbing-ball}
		is tempered.
		Moreover,  apply \eqref{eq:attractor-size-bound}
		at \(\theta_s\omega\), giving
		\[
		\begin{aligned}
			0\leq e^{-\varepsilon|s|}
			\sup_{\mathscr T\in\mathcal A_\mu(\theta_s\omega)}
			\|\mathscr T\|_{\mathcal Y}
			\leq 2e^{-\varepsilon|s|}
			R_{1,\mu}(\theta_s\omega)
			\to0
			\qquad\text{as }|s|\to\infty.
		\end{aligned}
		\]
		Thus \(\mathcal A_\mu\) admits a tempered random bound.	
	\end{proof}

	\begin{remark}\label{rem:attractor-scope}
		For \(\mu>0\), temperedness alone does not provide the required
		convergence under the general kernel assumptions.
		If the additional estimate
		\(
		p_{\mu,\varrho}(T)\leq Ce^{-cT}
		\)
		holds for some constants \(C,c>0\) and \(T\geq0\), then every such
		tempered random set satisfies
		\eqref{eq:admissible-random-initial-set}.
		When \(\mu=0\),  
		\(p_{0,\varrho}(T)=e^{-\varrho T}\),  every nonempty
		closed bounded tempered random set belongs to
		\(\mathfrak D_{0,\varrho}\).
	\end{remark}

	\begin{lemma}\label{prop:attractor-physical-projection}
		Define
		\begin{equation*}
			\mathscr A_\mu(\omega)
			:=\{\mathscr T(0):\mathscr T\in\mathcal A_\mu(\omega)\}.
		\end{equation*}		
		Then \(\mathscr{A}_\mu\) is a nonempty random compact set
		in \(X\). For every \(B\in\mathfrak D_{\mu,\varrho}\),
		\[
		\sup_{x\in B(\theta_{-T}\omega)}
		\inf_{y\in\mathscr{A}_\mu(\omega)}
		\|u(T,\theta_{-T}\omega;x)-y\|_X
		\to 0
		\quad\text{as }T\to\infty,
		\qquad \mathbb P\text{-almost surely}.
		\]
	\end{lemma}
	
	\begin{proof}
		By \eqref{eq:profile-norm}, the evaluation map
		\(\iota_0:\mathcal Y\to X\), \(\iota_0(\mathscr T)=\mathscr T(0)\),
		satisfies
		\[
		\|\iota_0(\mathscr T_1)-\iota_0(\mathscr T_2)\|_X
		\leq \|\mathscr T_1-\mathscr T_2\|_{\mathcal Y}.
		\]
		Hence \(\mathscr A_\mu(\omega)=\iota_0(\mathcal A_\mu(\omega))\)
		is nonempty and compact by
		Theorem~\ref{thm:generalized-random-pullback-attractor}.

		Since \(\mathcal A_\mu\) is a random compact set in the separable Banach
		space \(\mathcal Y\), there exist measurable maps
		\(\mathscr T_j:\mathcal O_{1,\mu}\to\mathcal Y\), \(j\ge1\) such that
		\(\{\mathscr T_j(\omega)\}_{j\ge1}\) is a dense subset of
		\(\mathcal A_\mu(\omega)\) by   the Castaing representation theorem.
		For each \(x\in X\), continuity of \(\iota_0\) gives
		\[
		\operatorname{dist}_X(\{x\},\mathscr A_\mu(\omega))
		=\inf_{j\geq1}\|x-\iota_0(\mathscr T_j(\omega))\|_X.
		\]
		The right side of the above formula is measurable, so \(\mathscr A_\mu\) is a random
		compact set in \(X\).
		
		Finally, fix \(B\in\mathfrak D_{\mu,\varrho}\).
		By \eqref{eq:gamma-evaluation-at-zero}  and
		\eqref{eq:Psi-x}, we have
		\(
		u(T,\theta_{-T}\omega;x)
		=\iota_0\bigl(\Gamma_\mu(T,\theta_{-T}\omega,\Psi_{\mu,x})\bigr).
		\)
		Consequently, for every \(x\in B(\theta_{-T}\omega)\), we have
		\[
		\begin{aligned}
			\inf_{y\in\mathscr{A}_\mu(\omega)}
			\|u(T,\theta_{-T}\omega;x)-y\|_X
			\leq\inf_{\mathscr{T}\in\mathcal A_\mu(\omega)}
			\|\Gamma_\mu(T,\theta_{-T}\omega,\Psi_{\mu,x})-\mathscr{T}\|_{\mathcal Y}.
		\end{aligned}
		\]
		By \ref{thm:generalized-random-pullback-attractor}(ii), as \(T\to\infty \),
		we obtain
		\[
		\begin{aligned}
			\sup_{x\in B(\theta_{-T}\omega)}
			\inf_{y\in\mathscr{A}_\mu(\omega)}
			\|u(T,\theta_{-T}\omega;x)-y\|_X&\leq
			\operatorname{dist}_{\mathcal Y}
			\bigl(
			\Gamma_\mu(T,\theta_{-T}\omega,
			\widetilde B(\theta_{-T}\omega)),
			\mathcal A_\mu(\omega)
			\bigr)
			\to 0.
		\end{aligned}
		\]
	\end{proof}
	\begin{remark}\label{rem:attractor}
		For \(\mu>0\), the cocycle \(\Gamma_\mu\) acts on \(\mathcal Y\)
		and has the random pullback attractor \(\mathcal A_\mu\).
		For \(\Psi_{\mu,x}=P_\mu(\cdot)x\), evaluation at \(r=0\) recovers
		the solution of \eqref{e1.1}.
		Thus \(\mathscr A_\mu(\omega)
		=\{\Psi(0):\Psi\in\mathcal A_\mu(\omega)\}\) is the part of
		\(\mathcal A_\mu(\omega)\) observable in \(X\).
	\end{remark}

	\section{Convergence of random attractors as \(\mu\to0\)}
	\label{subsec:vanishing-memory}
	We compare the random pullback attractors \(\mathcal A_\mu\) in the
	common space \(\mathcal Y\) as \(\mu\to0\).   Throughout this section,
	suppose that (\ref{ass:kernel}), (\ref{ass:kernel-holder}),
	(\ref{ass:noise-coefficient}), and (\ref{ass:nonlinearity}) hold,
	with \(L<\lambda_1\) and \(\varrho:=\lambda_1-L\).
	At \(\mu=0\), \eqref{e1.1} becomes
	\begin{equation}\label{eq:6.1}
		\begin{cases}
			\mathrm du(t)=(A u+f(u))\,\mathrm dt
			+G\,\mathrm dB^H(t), & t>0,\  \\ 
			u(0)=x\in X.
		\end{cases}
	\end{equation}
	By \eqref{eq:pnu-laplace} and \eqref{eq:P-spectral-definition},
	\(p_{0,\nu}(t)=e^{-\nu t}\) and \(P_0(t)=e^{tA}\). In particular,
	\begin{equation}\label{eq:semigroup property}
		P_0(r+a)=P_0(r)P_0(a),\qquad r,a\geq0.
	\end{equation}
	For \(\mu=0\), \(r\geq0\), \(x\in X\), and \(\mathscr T\in\mathcal Y\), by \eqref{eq:Psi-x}, \eqref{eq:profile-norm},
	\eqref{eq:Psi-xy}, and \eqref{eq:Psi-x-bound}, we have 
	\begin{equation*}
		\Psi_{0,x}(r)=P_0(r)x,	\quad   \Psi_{0,x}(0)=x,
	\end{equation*}
	and
	\begin{equation}\label{eq:vm-evaluation}
		\|\mathscr T(0)\|_X\leq\|\mathscr T\|_{\mathcal Y},\quad
		\|\Psi_{0,x}\|_{\mathcal Y}\leq2\|x\|_X.
	\end{equation}
	
	By Theorem~\ref{thm:generalized-random-pullback-attractor} with
	\(\mu=0\), \(\mathcal A_0\) is the random pullback attractor
	for \(\Gamma_0\) in \(\mathcal Y\) relative to
	\(\mathfrak D_{\widetilde B}\).
	We now show that
	\(
	\mathscr A_0(\omega)
	=\{\mathscr T(0):\mathscr T\in\mathcal A_0(\omega)\},
	\)
	defined in Lemma~\ref{prop:attractor-physical-projection}, is the
	tempered random pullback attractor of \eqref{eq:6.1} on \(X\). Equations  \eqref{eq:semigroup property}, \eqref{eq:lifted-map},
	\eqref{eq:mild-stochastic} and \eqref{eq:Psi-x} give
	\begin{equation}\label{eq:vm-semigroup-lift}
		\begin{aligned}
			\Gamma_0(t,\omega,\Psi_{0,x})(r)
			&=
			P_0(r)\Bigl[
			P_0(t)x
			+\int_0^t P_0(t-s)f(u_0(s,\omega,x))\,ds
			+Z_0(t,0,\omega)
			\Bigr]\\
			&=
			P_0(r)u_0(t,\omega,x)
			=
			\Psi_{0,u_0(t,\omega,x)}(r).
		\end{aligned}
	\end{equation}
	Fix \(\omega\in\mathcal O_{1,0}\) and
	\(\mathscr T\in\mathcal A_0(\omega)\). Let  \(y_j:=u_0(T_j,\theta_{-T_j}\omega,0)\).
	By \eqref{eq:attractor-sequential-characterization}
	and \eqref{eq:vm-semigroup-lift}, there exist \(T_j\to\infty\)
	such that 
	\(
	\Psi_{0,y_j}
	=\Gamma_0(T_j,\theta_{-T_j}\omega,\Psi_0)
	\to\mathscr T\)
	in \(\mathcal Y
	\).
	By \eqref{eq:vm-evaluation} and  \(y_j\to\mathscr T(0)\) in \(X\),
	\[
	\|\Psi_{0,y_j}-\Psi_{0,\mathscr T(0)}\|_{\mathcal Y}
	\leq 2\|y_j-\mathscr T(0)\|_X\to0.
	\]
	Thus \(\mathscr T=\Psi_{0,\mathscr T(0)}\), and the definition of
	\(\mathscr A_0(\omega)\) gives
	\begin{equation}\label{eq:vm-limit-attractor-lift}
		\mathcal A_0(\omega)
		=\{\Psi_{0,x}:x\in\mathscr A_0(\omega)\}.
	\end{equation}
	Using 	\eqref{eq:attractor-strict-invariance},
	\eqref{eq:vm-semigroup-lift} and \eqref{eq:vm-limit-attractor-lift}
	yields
	\begin{equation}\label{eq:A0-invariance}
		\mathscr A_0(\theta_t\omega)
		=\{u_0(t,\omega,y):y\in\mathscr A_0(\omega)\},
		\qquad t\geq0.
	\end{equation}
	By Lemma~\ref{prop:attractor-physical-projection},
	Remark~\ref{rem:attractor-scope},  Corollary~\ref{cor:attractor-tempered}  and the strict invariance
	\eqref{eq:A0-invariance}, \(\mathscr A_0\) is the tempered random
	pullback attractor of \eqref{eq:6.1} on \(X\).

	To estimate integrals over the entire past uniformly for small \(\mu\),
	we first construct a common integrable bound for the  kernels.

	\begin{lemma}\label{lem:vm-kernel}
		There exist \(\widetilde{\mu}>0\) and a nonnegative function
		\(\mathscr G\in C([0,\infty))\cap L^1(0,\infty)\) such that
		\begin{equation}\label{eq:vm-envelope}
			0\leq p_{\mu,\nu}(t)\leq \mathscr G(t),
			\qquad 0\leq\mu\leq\widetilde{\mu},\quad \nu\geq\varrho,
			\quad t\geq0.
		\end{equation}
		Moreover, as \(\mu\to0\),
		\begin{equation}\label{eq:vm-global-kernel}
			\sup_{t\geq0}\|P_\mu(t)-P_0(t)\|_{\mathcal L(X)}
			\leq2\mu\widehat k(\lambda_1)
			\to 0,
		\end{equation}
		and
		\begin{equation}\label{eq:vm-derivative}
			\sup_{t\geq\varepsilon}
			\|P_\mu'(t)-P_0'(t)\|_{\mathcal L(X)}
			\to 0
			\qquad\text{for every  } \varepsilon>0.
		\end{equation}
	\end{lemma}

	\begin{proof}
		Fix \(\nu>0\) and set \(e_\nu(t):=e^{-\nu t}\) and
		\(g_\nu:=e_\nu*k\). Since
		\(g_\nu(t)=e^{-\nu t}\int_0^t e^{\nu s}k(s)\,ds\),
		\(g_\nu\) is locally absolutely continuous, \(g_\nu(0)=0\), and
		\(h_\nu(t):=g_\nu'(t)=k(t)-\nu g_\nu(t)\) for  \(t>0\).
		By (\ref{ass:kernel}), \(k\) is nonnegative and nonincreasing.
		Hence, for \(t>0\),
		\[
		\nu g_\nu(t)
		\geq k(t)\int_0^t \nu e^{-\nu(t-s)}\,ds
		=(1-e^{-\nu t})k(t),
		\]
		so \(h_\nu(t)=k(t)-\nu g_\nu(t)\leq e^{-\nu t}k(t)\).
		The triangle inequality therefore gives
		\[
		|h_\nu(t)|
		\leq e^{-\nu t}k(t)
		+\bigl|e^{-\nu t}k(t)-h_\nu(t)\bigr|
		=2e^{-\nu t}k(t)-h_\nu(t).
		\]
		Using	\(g_\nu(0)=0\leq g_\nu(R)\), we have
		\[
		\int_0^R |h_\nu(t)|\,dt
		\leq 2\int_0^R e^{-\nu t}k(t)\,dt-g_\nu(R)
		\leq 2\widehat k(\nu).
		\]
		Letting \(R\to\infty\), we obtain
		\begin{equation}\label{eq:vm-h-integrability}
			\|h_\nu\|_{L^1(0,\infty)}\leq 2\widehat k(\nu).
		\end{equation}
		By \eqref{eq:mhat-def} and \eqref{eq:pnu-laplace}, we have
		\(
		\widehat h_\nu(\lambda)
		=\frac{\lambda\widehat k(\lambda)}{\lambda+\nu},\)
		and 	\((1+\mu\widehat h_\nu(\lambda))\widehat p_{\mu,\nu}(\lambda)
		=\frac{1}{\lambda+\nu}=\widehat e_\nu(\lambda).
		\)
		Inverting the Laplace transform gives
		\begin{equation}\label{eq:vm-perturbation-identity}
			p_{\mu,\nu}=e_\nu-\mu h_\nu*p_{\mu,\nu}.
		\end{equation}
		
		Choose \(0<\widetilde\mu\leq1\) such that
		\(\widetilde\mu\|h_\varrho\|_{L^1}<1\), and set
		\begin{equation*}
			\mathscr G
			:=e_\varrho+\sum_{j=1}^{\infty}
			\widetilde\mu^j |h_\varrho|^{*j}*e_\varrho,
		\end{equation*}
		where \(|h_\varrho|^{*j}\) denotes the \(j\)-fold convolution of
		\(|h_\varrho|\) with itself.
		For each \(j\),	\(\widetilde\mu^j |h_\varrho|^{*j}*e_\varrho\) is continuous. Since for \(v\in L^1(0,\infty)\), we have
		\(
		(v*e_\varrho)(t)
		=e^{-\varrho t}\int_0^t e^{\varrho s}v(s)\,ds.
		\)
		Young's inequality gives
		\[
		\bigl\|\widetilde\mu^j |h_\varrho|^{*j}*e_\varrho\bigr\|_{L^p}
		\leq
		\bigl(\widetilde\mu\|h_\varrho\|_{L^1}\bigr)^j
		\|e_\varrho\|_{L^p},
		\qquad p\in\{1,\infty\}.
		\]
		Hence the series converges uniformly and in \(L^1(0,\infty)\), with
		\(
		\|\mathscr G\|_{L^1}
		\leq
		\frac{1}{\varrho\bigl(1-\widetilde\mu\|h_\varrho\|_{L^1}\bigr)}.
		\)
		For \(0\leq\mu\leq\widetilde\mu\),
		\eqref{eq:vm-perturbation-identity} gives
		\(
		p_{\mu,\varrho}
		\leq e_\varrho+\widetilde\mu |h_\varrho|*p_{\mu,\varrho}.
		\)
		Iterating this inequality \(m\geq1\) times and using
		\eqref{eq:pnu-basic-bound}, we obtain
		\[
		p_{\mu,\varrho}(t)
		\leq \mathscr G(t)
		+\widetilde\mu^m
		\bigl(|h_\varrho|^{*m}*p_{\mu,\varrho}\bigr)(t)
		\leq \mathscr G(t)
		+\bigl(\widetilde\mu\|h_\varrho\|_{L^1}\bigr)^m.
		\]
		Since \(\widetilde\mu\|h_\varrho\|_{L^1}<1\), letting \(m\to\infty\)
		gives \(p_{\mu,\varrho}(t)\leq\mathscr G(t)\).
		The monotonicity of \(p_{\mu,\nu}\) in \(\nu\) from
		Lemma~\ref{lem:pnu-properties} then yields \eqref{eq:vm-envelope}.

		By \eqref{eq:vm-perturbation-identity}, \eqref{eq:pnu-basic-bound},
		and \eqref{eq:vm-h-integrability},
		\[
		\|p_{\mu,\nu}-e_\nu\|_{L^\infty}
		\leq\mu\|h_\nu\|_{L^1}\|p_{\mu,\nu}\|_{L^\infty}
		\leq2\mu\widehat k(\nu).
		\]
		Taking \(\nu=\lambda_n\) and using
		\(\widehat k(\lambda_n)\leq\widehat k(\lambda_1)\) in
		\eqref{eq:P-spectral-definition} gives \eqref{eq:vm-global-kernel}.

		Finally, by Taylor's formula, \eqref{eq:two-derivative},
		\eqref{eq:P-spectral-definition}, and \eqref{eq:P-derivative-spectral},
		for \(\varepsilon,h>0\),
		\begin{equation}\label{eq:PmuP0}
			\begin{aligned}
				&\sup_{t\geq\varepsilon}
				\left\|
				\frac{P_\mu(t+h)-P_\mu(t)}{h}-P_\mu'(t)
				\right\|_{\mathcal L(X)}
				\leq
				\frac{h}{2}
				\sup_{\substack{n\geq1\\ s\geq\varepsilon}}
				|p_{\mu,\lambda_n}''(s)|
				\leq \frac{2h}{e^2\varepsilon^2}.
			\end{aligned}
		\end{equation}	
		Applying \eqref{eq:PmuP0} to \(P_\mu\) and \(P_0\) and using
		\eqref{eq:vm-global-kernel}, we obtain
		\[
		\begin{aligned}
			\sup_{t\geq\varepsilon}\|P_\mu'(t)-P_0'(t)\|_{\mathcal L(X)}
			&\leq\frac{4h}{e^2\varepsilon^2}
			+\frac{2}{h}\sup_{t\geq0}\|P_\mu(t)-P_0(t)\|_{\mathcal L(X)}
			\leq\frac{4h}{e^2\varepsilon^2}
			+\frac{4\mu\widehat k(\lambda_1)}{h}.
		\end{aligned}
		\]
		Taking \(h=\sqrt{\mu}\) for \(\mu>0\) proves \eqref{eq:vm-derivative}.
	\end{proof}		
	
	\begin{lemma}\label{lem:vm-common-bounds}
		There exists a measurable invariant set
		\(\mathcal O_{**}\subseteq\mathcal O_0\) with
		\(\mathbb P(\mathcal O_{**})=1\). For 
		\(\omega\in\mathcal O_{**}\) and  \(0\leq\mu\leq\widetilde\mu\),
		the attractor \(\mathcal A_\mu(\omega)\) is well defined
		by \eqref{eq:attractor-definition}.
		Define
		\begin{equation}\label{eq:vm-state-radius}
			\overline R_0(\omega)
			:=\frac{\|f(0)\|_X}{\varrho}+\mathcal N(\omega)
			+L\int_0^\infty \mathscr G(t)\mathcal N(\theta_{-t}\omega)\,dt.
		\end{equation}
		Then \(\overline R_0\) is finite, measurable and
		\(
		R_{0,\mu}(\omega)\leq\overline R_0(\omega)\). The map
		\(s\to \overline R_0(\theta_s\omega)\) is continuous and 
		\(1+\overline R_0\) is tempered.
		For all  \(\mu\), \(t\geq s\), and \(x\in X\),
		\begin{equation}\label{eq:vm-state-bound}
			\|u_\mu(t-s,\theta_s\omega;x)\|_X
			\leq p_{\mu,\varrho}(t-s)\|x\|_X
			+\overline R_0(\theta_t\omega).
		\end{equation}
		Moreover,
		\begin{equation}\label{eq:vm-attractor-bound}
			\sup_{\mathscr{T}\in\mathcal A_\mu(\omega)}\|\mathscr{T}(0)\|_X
			\leq\overline R_0(\omega).
		\end{equation}
		Furthermore,
		\begin{equation}\label{eq:vm-drift-integrability}
			\int_0^\infty \mathscr G(t)
			\bigl(\|f(0)\|_X+L\overline R_0(\theta_{-t}\omega)\bigr)\,dt<\infty.
		\end{equation}
	\end{lemma}

	\begin{proof}
		With \(\mathcal N_N\) as in \eqref{eq:attractor-noise-N}, define
		\[
		\mathcal O_{2}:=
		\left\{\omega\in\mathcal O_0:
		\int_0^\infty h(a)\mathcal N_N(\theta_{-a}\omega)\,da<\infty
		\ \text{for every }h\in\{\mathscr G,\mathscr G*\mathscr G\},
		\ N\in\mathbb N\right\}.
		\]
		By Lemma~\ref{lem:vm-kernel} and Young's inequality,
		\(\mathscr G\) and \(\mathscr G*\mathscr G\) are nonnegative and integrable.
		Thus the proof of Lemma~\ref{lem:attractor-invariant-good-set},
		applied to   \(\mathscr G\) and \(\mathscr G*\mathscr G\) in \eqref{eq:attractor-good-set-preliminary},
		gives that \(\mathcal O_{2}\) is measurable, invariant,
		and 	\(\mathbb P(\mathcal O_{2})=1\). For every \(0\leq\mu\leq\widetilde\mu\),
		\eqref{eq:vm-envelope} gives
		\[
		p_{\mu,\varrho},\ p_{\mu,\lambda_1}\leq\mathscr G,
		\qquad
		p_{\mu,\lambda_1}*p_{\mu,\varrho}\leq\mathscr G*\mathscr G.
		\]
		Hence \(\mathcal O_{2}\subseteq\mathcal O_{1,\mu}\)
		by \eqref{eq:attractor-good-set-preliminary}, and
		Theorem~\ref{thm:generalized-random-pullback-attractor}
		gives \(\mathcal A_\mu(\omega)\) as in
		\eqref{eq:attractor-definition} for every
		\(\omega\in\mathcal O_{2}\).
		
		Define \(\overline R_0\) by \eqref{eq:vm-state-radius} on
		\(\mathcal O_{2}\) and by zero on its complement.
		This defines a finite measurable random variable by the definition
		of \(\mathcal O_{2}\) and joint measurability of the shifts.
		For \(\omega\in\mathcal O_{2}\) and \(|s|\leq N\),
		\eqref{eq:attractor-noise-N} gives
		\[
		\mathscr G(a)\mathcal N(\theta_{s-a}\omega)
		\leq \mathscr G(a)\mathcal N_N(\theta_{-a}\omega),
		\]
		whose right-hand side is integrable in \(a\).
		Thus Lemma~\ref{lem:attractor-noise-majorant} and dominated convergence
		give continuity of \(s\mapsto\overline R_0(\theta_s\omega)\).
		Moreover, Fubini--Tonelli theorem, shift invariance of \(\mathbb P\),
		and \(\mathbb E\mathcal N_N<\infty\) in Lemma \ref{lem:attractor-noise-majorant} yield
		\begin{equation}\label{eq:vm-common}
			\mathbb E\sup_{|s|\leq N}\overline R_0(\theta_s\omega)
			\leq \frac{\|f(0)\|_X}{\varrho}
			+\bigl(1+L\|\mathscr G\|_{L^1}\bigr)\mathbb E\mathcal N_N
			<\infty.
		\end{equation}
		Applying the argument of Corollary~\ref{cor:attractor-tempered}
		to \(1+\overline R_0\), using \eqref{eq:vm-common} with \(N=1\),
		shows that the set
		\begin{equation}\label{eq:vm-common-tempered-set}
			\mathcal O_{**}:=
			\left\{\omega\in\mathcal O_{2}:
			\lim_{|s|\to\infty}e^{-\varepsilon|s|}
			\bigl(1+\overline R_0(\theta_s\omega)\bigr)=0
			\ \text{for every }\varepsilon>0\right\}
		\end{equation}
		is measurable, invariant, and of full measure.
		
		Fix \(\omega\in\mathcal O_{**}\) and
		\(0\leq\mu\leq\widetilde\mu\).
		Comparing \eqref{eq:attractor-state-radius} with
		\eqref{eq:vm-state-radius} and using \eqref{eq:vm-envelope}
		gives \(R_{0,\mu}\leq\overline R_0\) on \(\mathcal O_{**}\).
		Hence \eqref{eq:attractor-state-comparison} yields
		\eqref{eq:vm-state-bound}. Using \eqref{eq:gamma-evaluation-at-zero}, we obtain
		\begin{equation}\label{eq:gamma-R}
			\|\Gamma_\mu(T,\theta_{-T}\omega,\Psi_0)(0)\|_X
			=\|u_\mu(T,\theta_{-T}\omega;0)\|_X
			\leq\overline R_0(\omega),\qquad T\geq0.
		\end{equation}
		By \eqref{eq:attractor-sequential-characterization} and \(
		\Gamma_\mu(T_j,\theta_{-T_j}\omega,\Psi_0)(0)
		\to\mathscr T(0)
		\) in \(X\),
		passing to the limit in  \eqref{eq:gamma-R} gives
		\eqref{eq:vm-attractor-bound} as \(\mathscr T\) was arbitrary.
		
		Finally, substituting \eqref{eq:vm-state-radius}, applying
		Fubini--Tonelli theorem and using
		\(\omega\in\mathcal O_{2}\) and
		\(\mathcal N\leq\mathcal N_1\) in \eqref{eq:attractor-noise-N}, we obtain
		\[
		\begin{aligned}
			&\int_0^\infty \mathscr G(t)
			\bigl(\|f(0)\|_X+L\overline R_0(\theta_{-t}\omega)\bigr)\,dt\\
			&\quad=
			\left(1+\frac{L}{\varrho}\right)
			\|f(0)\|_X\|\mathscr G\|_{L^1}
			+L\int_0^\infty
			\bigl[\mathscr G(t)+L(\mathscr G*\mathscr G)(t)\bigr]
			\mathcal N(\theta_{-t}\omega)\,dt
			<\infty.
		\end{aligned}
		\]
		This proves \eqref{eq:vm-drift-integrability}.
	\end{proof}

	\begin{lemma}\label{lem:vm-noise}
		Let \(W_\mu\) be as defined in
		Lemma~\ref{lem:attractor-past-noise-convolution}.
		For every \(\omega\in\mathcal O_0\) and \(T>0\), we have
		\begin{align}
			&\sup_{|s|\leq T}
			\|W_\mu(\theta_s\omega)-W_0(\theta_s\omega)\|_\infty
			\to 0, \quad\text{as}\quad \mu\to0;
			\label{eq:vm-past-noise-convergence}\\
			&\sup_{0\leq t\leq T}\sup_{r\geq0}
			\|Z_\mu(t,r,\omega)-Z_0(t,r,\omega)\|_X
			\to0, \quad\text{as}\quad \mu\to0.
			\label{eq:vm-finite-noise-convergence}
		\end{align}
		Moreover,
		\begin{equation}\label{eq:vm-noise-semigroup}
			W_0(\omega)(r)=P_0(r)W_0(\omega)(0),
			\qquad r\geq0.
		\end{equation}
	\end{lemma}

	\begin{proof}
		Fix \(\omega\in\mathcal O_0\) and \(T>0\). Since \(\mathcal O_0\) is
		invariant, \eqref{eq:attractor-noise-ibp} applies to \(\theta_s\omega\).
		As the past interval expands to \((-\infty,0]\), its boundary term tends
		to zero uniformly in \(r\geq0\) and \(|s|\leq T\) by
		\eqref{eq:attractor-operator-derivative} and \eqref{eq:noise-shift-sup-tail}.
		Thus \eqref{eq:attractor-past-noise-limit}
		gives, for \(|s|\leq T\) and \(r\geq0\),
		\begin{equation}\label{eq:vm-past-noise-ibp}
			W_\mu(\theta_s\omega)(r)
			=\int_0^\infty P_\mu'(r+a)G
			\bigl(\omega(s-a)-\omega(s)\bigr)\,da.
		\end{equation}
		Indeed, \eqref{eq:shifted-noise-increment-bound} and
		Lemma~\ref{lem:attractor-kernel-estimates}\textup{(iii)} give
		\begin{equation}\label{eq:bp}
			b^*(a)\sup_{|s|\leq T}\|\omega(s-a)-\omega(s)\|_X
			\leq C_{T,\omega}
			\begin{cases}
				a^{\gamma-1}, & 0<a\leq1,\\
				a^{\zeta-2}, & a>1,
			\end{cases}
		\end{equation}
		which is integrable since \(\gamma>0\) and \(\zeta<1\).
		By \eqref{eq:attractor-operator-derivative},
		\begin{equation}\label{eq:i-bound}
			\sup_{v\geq a}\|P_\mu'(v)-P_0'(v)\|_{\mathcal L(X)}
			\leq 2b^*(a),\qquad a>0.
		\end{equation}
		The left-hand side tends to zero for each \(a>0\) by
		\eqref{eq:vm-derivative}. Hence \eqref{eq:vm-past-noise-ibp} and
		dominated convergence, using \eqref{eq:bp} and \eqref{eq:i-bound}, give
		\[
		\begin{aligned}
			&\sup_{|s|\leq T}
			\|W_\mu(\theta_s\omega)-W_0(\theta_s\omega)\|_\infty\\
			&\quad\leq \|G\|_{\mathcal L(X)}
			\int_0^\infty
			\sup_{v\geq a}\|P_\mu'(v)-P_0'(v)\|_{\mathcal L(X)}
			\sup_{|s|\leq T}\|\omega(s-a)-\omega(s)\|_X\,da
			\to0 .
		\end{aligned}
		\]
		
		By \eqref{eq:young-shift} and \eqref{eq:attractor-past-noise},
		splitting the integral at zero yields
		\[
		S_{\mu,\tau+t}(\theta_t\omega)(r)
		=S_{\mu,\tau}(\omega)(t+r)+Z_\mu(t,r,\omega).
		\]
		Letting \(\tau\to\infty\) by
		\eqref{eq:attractor-past-noise-limit} gives
		\(
		Z_\mu(t,r,\omega)
		=W_\mu(\theta_t\omega)(r)-W_\mu(\omega)(t+r).
		\)
		Hence, by \eqref{eq:vm-past-noise-convergence},
		\[
		\sup_{0\leq t\leq T}\sup_{r\geq0}
		\|Z_\mu(t,r,\omega)-Z_0(t,r,\omega)\|_X
		\leq 2\sup_{|s|\leq T}
		\|W_\mu(\theta_s\omega)-W_0(\theta_s\omega)\|_\infty
		\to0.
		\]
		Finally, \eqref{eq:semigroup property} implies
		\(
		P_0'(r+a)=P_0(r)P_0'(a)\)
		for \(r\geq0\) and \(a>0\).
		Substituting   into
		\eqref{eq:vm-past-noise-ibp}, we obtain
		\[
		W_0(\omega)(r)
		=
		P_0(r)\int_0^\infty P_0'(a)G\omega(-a)\,da
		=
		P_0(r)W_0(\omega)(0),
		\]
		which proves \eqref{eq:vm-noise-semigroup}.
	\end{proof}
	
	For the remaining attractor results, we restrict the metric dynamical
	system to \(\mathcal O_{**}\), equipped with the trace \(\sigma\)-algebra
	and the restricted probability measure.
	
	\begin{lemma}\label{lem:vm-history-comparison}
		For every \(\omega\in\mathcal O_{**}\),
		\begin{equation}\label{eq:vm-profile-defect}
			\sup_{\mathscr{T}\in\mathcal A_\mu(\omega)}
			\sup_{r\geq0}\|\mathscr{T}(r)-\Psi_{0,\mathscr T(0)}(r)\|_X
			\to 0\qquad\text{as }\mu\to0.
		\end{equation}
		Let  \(u_0(t,\omega,x)\) for the
		mild solution of \eqref{e1.1} with \(\mu=0\).
		Then, for every \(T>0\),
		\begin{equation}\label{eq:vm-restart-comparison}
			\sup_{\mathscr{T}\in\mathcal A_\mu(\omega)}
			\sup_{0\leq t\leq T}
			\|\Gamma_\mu(t,\omega,\mathscr{T})(0)
			-u_0(t,\omega,\mathscr{T}(0))\|_X
			\to 0\qquad\text{as }\mu\to0.
		\end{equation}
	\end{lemma}

	\begin{proof}
		Fix \(\omega\in\mathcal O_{**}\).
		For \(0\leq\mu\leq\widetilde\mu\), define
		\[
		D_\mu(a):=
		\sup_{r\geq0}
		\|P_\mu(r+a)-P_0(r)P_\mu(a)\|_{\mathcal L(X)},
		\qquad a\geq0.
		\]
		
		The semigroup property \eqref{eq:semigroup property} gives
		\[
		P_\mu(r+a)-P_0(r)P_\mu(a)
		=(P_\mu-P_0)(r+a)-P_0(r)(P_\mu-P_0)(a).
		\]
		Hence \eqref{eq:P-exact-norm}, Lemma~\ref{lem:pnu-properties},
		\eqref{eq:vm-envelope}, and \eqref{eq:vm-global-kernel} give
		\begin{equation}\label{eq:D-G}
			\begin{aligned}
				D_\mu(a)\leq 2p_{\mu,\lambda_1}(a)\leq2\mathscr G(a)
				\quad 	\text{and}\quad 
				D_\mu(a)\leq 2\sup_{v\geq0}
				\|P_\mu(v)-P_0(v)\|_{\mathcal L(X)}
				\leq4\mu\widehat k(\lambda_1).
			\end{aligned}
		\end{equation}
		For \(a>0\), the supremum defining \(D_\mu(a)\) may be taken over
		nonnegative rational \(r\) by Lemma~\ref{lem:attractor-kernel-estimates}
		\textup{(iii)}, so \(D_\mu\) is measurable.
		By \eqref{eq:D-G}, \eqref{eq:vm-drift-integrability}, and
		dominated convergence,
		\begin{equation}\label{eq:vm-drift-defect}
			J_\mu(\omega):=\int_0^\infty D_\mu(a)
			\bigl(\|f(0)\|_X+L\overline R_0(\theta_{-a}\omega)\bigr)\,da
			\to0, \quad \mu\to 0.
		\end{equation}
		For \(0\leq a\leq\tau\), applying \eqref{eq:vm-state-bound} with
		\((s,t,x)=(-\tau,-a,0)\) and using \eqref{eq:f-lipschitz} gives
		\begin{equation}\label{eq:f-bound}
			\|f(u_\mu(\tau-a,\theta_{-\tau}\omega,0))\|_X
			\leq \|f(0)\|_X+L\overline R_0(\theta_{-a}\omega).
		\end{equation}
		By \eqref{eq:attractor-profile-decomposition} with \(x=0\),
		\[
		\begin{aligned}
			&\Gamma_\mu(\tau,\theta_{-\tau}\omega,\Psi_0)(r)
			-P_0(r)\Gamma_\mu(\tau,\theta_{-\tau}\omega,\Psi_0)(0)\\
			&\quad=\int_0^\tau
			[P_\mu(r+a)-P_0(r)P_\mu(a)]
			f(u_\mu(\tau-a,\theta_{-\tau}\omega,0))\,da
			+S_{\mu,\tau}(\omega)(r)-P_0(r)S_{\mu,\tau}(\omega)(0).
		\end{aligned}
		\]
		By \eqref{eq:f-bound} and \eqref{eq:vm-drift-defect}, for every
		\(\tau,r\geq0\),
		\begin{equation}\label{eq:f-bound1}
			\begin{aligned}
				&\|\Gamma_\mu(\tau,\theta_{-\tau}\omega,\Psi_0)(r)
				-P_0(r)\Gamma_\mu(\tau,\theta_{-\tau}\omega,\Psi_0)(0)\|_X
				\leq J_\mu(\omega)
				+\|S_{\mu,\tau}(\omega)(r)-P_0(r)S_{\mu,\tau}(\omega)(0)\|_X .
			\end{aligned}
		\end{equation}
		Fix \(\mathscr T\in\mathcal A_\mu(\omega)\).
		By \eqref{eq:attractor-sequential-characterization}, there exists
		\(\tau_j\to\infty\) such that
		\(\Gamma_\mu(\tau_j,\theta_{-\tau_j}\omega,\Psi_0)\to\mathscr T\)
		in \(\mathcal Y\).
		For each fixed \(r\geq0\), as \(j\to\infty\),   \eqref{eq:profile-local-control} gives
		\(
		\Gamma_\mu(\tau_j,\theta_{-\tau_j}\omega,\Psi_0)(r)
		\to \mathscr T(r)\) and
		\(	\Gamma_\mu(\tau_j,\theta_{-\tau_j}\omega,\Psi_0)(0)
		\to \mathscr T(0)\)
		in \(X\).
		Taking \(\tau=\tau_j\) in \eqref{eq:f-bound1} and using
		\eqref{eq:attractor-past-noise-limit}, we obtain
		\[
		\|\mathscr T(r)-P_0(r)\mathscr T(0)\|_X
		\leq J_\mu(\omega)
		+\|W_\mu(\omega)(r)-P_0(r)W_\mu(\omega)(0)\|_X .
		\]
		Taking the suprema over \(r\geq0\) and
		\(\mathscr T\in\mathcal A_\mu(\omega)\), and using
		\eqref{eq:vm-noise-semigroup} and \eqref{eq:P-exact-norm}, gives
		\[
		\begin{aligned}
			&\sup_{\mathscr T\in\mathcal A_\mu(\omega)}\sup_{r\geq0}
			\|\mathscr T(r)-P_0(r)\mathscr T(0)\|_X
			\leq J_\mu(\omega)
			+2\|W_\mu(\omega)-W_0(\omega)\|_\infty
			\to0,\quad\mu\to0,
		\end{aligned}
		\]
		by \eqref{eq:vm-drift-defect} and \eqref{eq:vm-past-noise-convergence}.
		Since \(\Psi_{0,\mathscr T(0)}(r)=P_0(r)\mathscr T(0)\) by
		\eqref{eq:Psi-x}, this proves \eqref{eq:vm-profile-defect}.

		Fix \(T>0\). By \eqref{eq:vm-state-bound} at \(\mu=0\),
		\eqref{eq:vm-attractor-bound}, and Lemma~\ref{lem:vm-common-bounds},
		\begin{equation}\label{eq:u0-finiteness}
			\sup_{\mathscr T\in\mathcal A_\mu(\omega)}\sup_{0\leq s\leq T}
			\|u_0(s,\omega,\mathscr T(0))\|_X
			\leq \overline R_0(\omega)
			+\sup_{0\leq s\leq T}\overline R_0(\theta_s\omega)<\infty.
		\end{equation}
		For \(a\geq0\) and \(w_1,w_2\in X\), \eqref{eq:f-lipschitz},
		\eqref{eq:P-exact-norm}, and \eqref{eq:vm-global-kernel} give
		\[
		\begin{aligned}
			\|P_\mu(a)f(w_1)-P_0(a)f(w_2)\|_X
			&\leq \|P_\mu(a)[f(w_1)-f(w_2)]\|_X
			+\|(P_\mu(a)-P_0(a))f(w_2)\|_X\\
			&\leq L\|w_1-w_2\|_X+2\mu\widehat k(\lambda_1)\|f(w_2)\|_X .
		\end{aligned}
		\]
		Subtract \eqref{eq:mild-stochastic} with \(\mu=0\) and
		\(x=\mathscr T(0)\) from \eqref{eq:uPsi} with \(\Psi=\mathscr T\).
		Using \eqref{eq:gamma-evaluation-at-zero}, the preceding estimate,
		and Gronwall's inequality, we obtain
		\[
		\begin{aligned}
			&\sup_{0\leq t\leq T}
			\|\Gamma_\mu(t,\omega,\mathscr T)(0)
			-u_0(t,\omega,\mathscr T(0))\|_X\\
			&\quad\leq e^{LT}\Bigl[
			\sup_{r\geq0}\|\mathscr T(r)-P_0(r)\mathscr T(0)\|_X
			+\sup_{0\leq s\leq T}\|Z_\mu(s,0,\omega)-Z_0(s,0,\omega)\|_X\\
			&\qquad\qquad
			+2\mu\widehat k(\lambda_1)
			\int_0^T\|f(u_0(s,\omega,\mathscr T(0)))\|_X\,ds\Bigr].
		\end{aligned}
		\]
		The \(\int_0^T\|f(u_0(s,\omega;\mathscr T(0)))\|_X\) is uniformly bounded by \eqref{eq:f-lipschitz}
		and \eqref{eq:u0-finiteness}. As \(\mu\to0\), taking the supremum over
		\(\mathscr T\in\mathcal A_\mu(\omega)\) and using
		\eqref{eq:vm-profile-defect} and \eqref{eq:vm-finite-noise-convergence}
		proves \eqref{eq:vm-restart-comparison}.
	\end{proof}			
	
	\begin{definition}\label{def:vm-upper-semicontinuity}\cite[Definition~7.1]{banaskiewicz_kalita}
		The family of random attractors \(\{\mathcal A_\mu\}_{\mu\geq0}\)
		is pathwise upper semicontinuous at \(\mu=0\) in \(\mathcal Y\) if
		\[
		\lim_{\mu\to0}
		\operatorname{dist}_{\mathcal Y}
		\bigl(\mathcal A_\mu(\omega),\mathcal A_0(\omega)\bigr)=0
		\qquad\text{for every }\omega\in\mathcal O_{**}.
		\]
	\end{definition}

	\begin{theorem}\label{thm:vm-profile-convergence}
		Suppose that (\ref{ass:kernel}),
		(\ref{ass:kernel-holder}), (\ref{ass:noise-coefficient}), and
		(\ref{ass:nonlinearity}) hold, with \(L<\lambda_1\).
		For \(0\leq\mu\leq\widetilde\mu\), let \(\mathcal A_\mu\) be the minimal random pullback attractor
		constructed in Theorem~\ref{thm:generalized-random-pullback-attractor}.
		Then these attractors are upper semicontinuous at \(\mu=0\)
		in \(\mathcal Y\). More precisely,
		\begin{equation}\label{eq:vm-profile-convergence}
			\lim_{\mu\to0}
			\operatorname{dist}_{\mathcal Y}
			\bigl(\mathcal A_\mu(\omega),\mathcal A_0(\omega)\bigr)=0,
			\qquad\omega\in\mathcal O_{**}.
		\end{equation}
	\end{theorem}
	
	\begin{proof}

		Fix \(\omega\in\mathcal O_{**}\), \(T>0\),
		\(0\le\mu\le\widetilde\mu\), and
		\(\mathscr T\in\mathcal A_\mu(\omega)\).
		By \eqref{eq:attractor-strict-invariance}, there exists
		\(\widetilde\Xi\in\mathcal A_\mu(\theta_{-T}\omega)\) such that
		\begin{equation}\label{eq:P-G-L}
			\mathscr T=\Gamma_\mu(T,\theta_{-T}\omega,\widetilde\Xi).
		\end{equation}
		Since \(p_{0,\varrho}(T)=e^{-\varrho T}\),
		\eqref{eq:attractor-state-difference} gives
		\begin{equation}\label{eq:vm-limit-contraction}
			\|u_0(T,\theta_{-T}\omega,x)-u_0(T,\theta_{-T}\omega,y)\|_X
			\leq e^{-\varrho T}\|x-y\|_X,\qquad x,y\in X.
		\end{equation}
		Choose \(y\in\mathscr A_0(\theta_{-T}\omega)\).
		By \eqref{eq:A0-invariance} and \eqref{eq:vm-attractor-bound},
		\[
		u_0(T,\theta_{-T}\omega,y)\in\mathscr A_0(\omega)
		\quad\text{and}\quad 
		\|\widetilde\Xi(0)-y\|_X
		\leq 2\overline R_0(\theta_{-T}\omega).
		\]
		Hence \eqref{eq:P-G-L} and \eqref{eq:vm-limit-contraction} imply
		\begin{equation}\label{eq:inf-R0}
			\begin{aligned}
				\inf_{\widetilde{y}\in\mathscr A_0(\omega)}\|\mathscr T(0)-\widetilde{y}\|_X
				&\leq \|\mathscr T(0)-u_0(T,\theta_{-T}\omega,y)\|_X\\
				&\leq \|\Gamma_\mu(T,\theta_{-T}\omega,\widetilde\Xi)(0)
				-u_0(T,\theta_{-T}\omega,\widetilde\Xi(0))\|_X
				+2e^{-\varrho T}\overline R_0(\theta_{-T}\omega).
			\end{aligned}
		\end{equation}
		Since \(\Psi_{0,\mathscr T(0)}(0)=\mathscr T(0)\),
		\eqref{eq:profile-norm}, \eqref{eq:vm-evaluation}, and
		\eqref{eq:vm-limit-attractor-lift} give
		\begin{equation}\label{eq:inf-bound}
			\begin{aligned}
				\inf_{\mathscr U\in\mathcal A_0(\omega)}
				\|\mathscr T-\mathscr U\|_{\mathcal Y}
				&\leq \|\mathscr T-\Psi_{0,\mathscr T(0)}\|_{\mathcal Y}
				+\inf_{\widetilde y\in\mathscr A_0(\omega)}
				\|\Psi_{0,\mathscr T(0)-\widetilde{y}}\|_{\mathcal Y}\\
				&\leq \sup_{r\geq0}\|\mathscr T(r)-P_0(r)\mathscr T(0)\|_X
				+2\inf_{\widetilde y\in\mathscr A_0(\omega)}\|\mathscr T(0)-\widetilde{y}\|_X .
			\end{aligned}
		\end{equation}
		
		For \(T>0\) and \(0\leq\mu\leq\widetilde\mu\), set
		\[
		\begin{aligned}
			\mathcal E_\mu(T,\omega)
			&:= \sup_{\mathscr T\in\mathcal A_\mu(\omega)}\sup_{r\geq0}
			\|\mathscr T(r)-P_0(r)\mathscr T(0)\|_X\\
			&\quad
			+2\sup_{\Xi\in\mathcal A_\mu(\theta_{-T}\omega)}
			\|\Gamma_\mu(T,\theta_{-T}\omega,\Xi)(0)
			-u_0(T,\theta_{-T}\omega,\Xi(0))\|_X .
		\end{aligned}
		\]
		Substituting \eqref{eq:inf-R0} into \eqref{eq:inf-bound} and
		taking the supremum over \(\mathscr T\in\mathcal A_\mu(\omega)\),
		we obtain
		\begin{equation}\label{eq:vm-upper-semicontinuity-bound}
			\operatorname{dist}_{\mathcal Y}
			\bigl(\mathcal A_\mu(\omega),\mathcal A_0(\omega)\bigr)
			\leq \mathcal E_\mu(T,\omega)
			+4e^{-\varrho T}\overline R_0(\theta_{-T}\omega).
		\end{equation}
		
		Given \(\varepsilon>0\), choose \(T>0\) such that
		\(4e^{-\varrho T}\overline R_0(\theta_{-T}\omega)<\varepsilon/2\)
		by \eqref{eq:vm-common-tempered-set}.
		For this fixed \(T\), \eqref{eq:vm-profile-defect} at \(\omega\)
		and \eqref{eq:vm-restart-comparison} at
		\(\theta_{-T}\omega\in\mathcal O_{**}\) give
		\(\mathcal E_\mu(T,\omega)<\varepsilon/2\) for all sufficiently
		small \(\mu>0\). 
		Thus
		\[
		\operatorname{dist}_{\mathcal Y}
		\bigl(\mathcal A_\mu(\omega),\mathcal A_0(\omega)\bigr)
		<\varepsilon,
		\]
		which proves \eqref{eq:vm-profile-convergence}.
		
	\end{proof}

	\section*{Competing interests}
	
	The authors declare that they have no competing interests.
	
	\section*{Acknowledgements}
	This work is supported by Natural Science Foundation of Guangxi (No. 2026GXNSFAA00640023)

\end{document}